\documentclass[12pt]{amsart}
\usepackage[utf8]{inputenc}
\usepackage[scale=0.80,centering]{geometry}

\usepackage{eye}
\eyesetup{interval width=12mm,level distance=1.5mm,first level=3.5mm}

\usepackage{bbm}
\usepackage{amssymb, mathrsfs}
\usepackage{amsmath}
\usepackage[all]{xy}
\usepackage{tikz}
\usepackage{amsfonts}
\usepackage{latexsym}
\usepackage{graphicx}
\usepackage{enumerate}
\usepackage{amscd,amssymb,amsmath,amsbsy,amsthm}
	\definecolor{linkred}{rgb}{0.7,0.2,0.2}
	\definecolor{linkblue}{rgb}{0,0.2,0.6}
	\definecolor{linkgreen}{rgb}{0,0.6,0.2}
\usepackage[colorlinks,plainpages,backref,
    linkcolor=linkblue,
    citecolor=linkgreen,
    urlcolor=linkred]{hyperref}

\newtheorem{Th}{Theorem}[section]
\newtheorem{Lemma}[Th]{Lemma}

\newtheorem{Coro}[Th]{Corollary}
\newtheorem{Prop}[Th]{Proposition}
\newtheorem{Eg}[Th]{Example}
\newtheorem{Rmk}[Th]{Remark}

\allowdisplaybreaks

\newtheorem{Thm}{Theorem}[section]

\def\semi{\propto}

\def\id{\operatorname{id}}
\def\pt{\operatorname{pt}}

\def\calO{\mathcal{O}}
\def\pt{\mathsf{pt}}

\begin{document}

\title{Schubert Calculus in the Eye of the Storm}

\author{Neil J.Y. Fan}
 \address{Department of Mathematics, Sichuan University, Chengdu, Sichuan 610065, P.R. China}
\email{fan@scu.edu.cn}

\author{Peter L. Guo}
\address{Center for Combinatorics, LPMC, Nankai University, Tianjin 300071, P.R. China}
\email{lguo@nankai.edu.cn}

\author{Changzheng Li}
 \address{School of Mathematics, Sun Yat-sen University, Guangzhou 510275, P.R. China}
\email{lichangzh@mail.sysu.edu.cn}

\author{Rui Xiong}
\address{Department of Mathematics, University of Michigan, Ann Arbor, MI 48109, USA}
\email{xiongrui@umich.edu}

\thanks{}

\maketitle

\begin{abstract}
   We introduce new combinatorial objects,  called   storm diagrams, to study the torus-equivariant quantum $K$-theory  of complete flag varieties   $GL_n(\mathbb{C})/B$. 
   We consider the quantum product of an arbitrary Schubert class with a
modified Hirzebruch $\lambda_y$-class of a tautological subbundle. We show that its Schubert expansion coefficients are given by
the total weights of storm diagrams. As applications, we establish     Pieri rules for the torus-equivariant quantum $K$-theory of partial flag varieties $GL_n(\mathbb{C})/P$. Our storm diagrams provide the first combinatorial model for a manifestly positive   multiplication formula for partial flag varieties beyond the Chevalley case, with restrictions imposed on only one of the two Schubert classes.
\end{abstract}

\setcounter{tocdepth}{1}

\section{Introduction}

\subsection{Background}
Let $G=GL_n(\mathbb{C})$,  $B$ be the Borel subgroup  consisting of upper triangular matrices of $G$, and $T\subset B$ be the diagonal torus. The  equivariant quantum $K$-theory $QK_T(G/B)$ of the complete flag variety $G/B$ is a deformation of the equivariant $K$-theory $K_T(G/B)$, and has very rich structures connecting  geometry, representation theory and combinatorics. For example,  Kato \cite{Kato1} proved that, after localization,  $QK_T(G/B)$ is isomorphic to the equivariant $K$-homology of the affine Grassmannian ${\rm Gr}_{SL_n}$. This establishes  a conjecture of Lam, the third named author,  Mihalcea and Shimozono \cite{LLMS} (see also \cite{IIM}), giving a $K$-theoretic analogue of  Peterson's isomorphism \cite{Pet97, LaSh} for quantum cohomology. 
   
   The   Schubert classes  $\mathcal{O}_u$, indexed by permutations $u\in  S_n$,   form an additive basis of $QK_T(G/B)$ over $K_T({\rm pt})[q_1, \cdots, q_{n-1}]$. The Schubert  structure coefficients $N_{u, v}^{w, d}$ in the quantum product  $$\mathcal{O}_u \star\mathcal{O}_v=\sum_{w\in S_n, d\in H_2(G/B, \mathbb{Z})} N_{u, v}^{w, d}\mathcal{O}_w {q}^d$$
are defined by a combination of  genus-zero equivariant $K$-theoretic Gromov--Witten invariants. They recover the multiplication in   $K_T(G/B)$ by taking the degree $d=0$ part, and recover the equivariant quantum cohomology ring $QH_T(G/B)$  
by specifying the coefficients appropriately \cite{BM11}. In both theories, the Schubert structure coefficients are known to satisfy positivity properties. For $QK_T(G/B)$, the corresponding  positivity property was  conjectured    in \cite[Conjecture 2.2]{BCMP}, as stated in  \eqref{eqn:conjPos}. This conjecture motivates the search for manifestly positive formulas for the sign-normalized
structure coefficients.
Till now, there is only  one known case, namely the Chevalley formula \cite{Kato8, LNS} which concerns  the case   $v=s_k$ and plays an important role in the study of $QK_T(G/B)$ \cite{LeMa, MNSI, MNSII}.  
In this paper, we establish a (positive) Pieri rule for $QK_T(G/B)$, describing the quantum product $\mathcal{O}_u\star \mathcal{O}_v$ for arbitrary $u\in S_n$ and  $v=s_is_{i+1}\cdots s_k$ with $1\leq i\leq k<n$. Our approach combines 
$$ \mbox{Combinatorics (storm diagrams)}  +  \mbox{Repr.  Theory (Weyl group action)}  + \mbox{Geometry (localization)}.$$

Even in the classical specialization $K_T(G/B)$, our formula is the first positive one for the multiplication  with an arbitrary  Schubert class $\mathcal{O}_u$, after  the Chevalley formula  obtained in \cite{FuLa, PiRa, GrRa} more than two decades ago.
Other known  positive formulae in $K_T(G/B)$ include those  
for the localization   coefficients $N_{u, w}^w$ \cite{KoKu, McN, GrKr},   products of Schubert classes pulled back from a common  complex Grassmannian by Pechenik and Yong \cite{PeYo} (see also \cite{WZJ}) or a two-step flag variety by Knutson and Zinn-Justin \cite{KZJ1}, and products indexed by permutations with separated descents by Knutson and Zinn-Justin  \cite{KZJ3}.

More generally, we establish Pieri rules in the equivariant quantum $K$-theory $QK_T(G/P)$ of partial flag varieties with $G=GL_n(\mathbb{C})$. For comparison, the equivariant quantum multiplication in projective spaces can be computed from the results of Buch and Mihalcea  \cite{BM11},  while a Littlewood--Richardson rule  for  the non-equivariant quantum $K$-theory of the two-step flag variety $F\ell_{1, n-1; n}$ was established by Xu \cite{Xu}. A Pieri rule for $QK(G/B)$ was given by Naito and Sagaki \cite{NaSa}. In    non-equivariant $K$-theory, Buch \cite{Buch} established a Littlewood--Richardson rule for complex Grassmannians, and Lenart and Sottile  \cite{LeSo} proved a Pieri rule for flag varieties $G/P$.

\subsection{Storm diagram model and main result}
\label{subsec:main-results}
\begin{figure}[hhh]
\begin{tikzpicture}[scale=0.8] \draw[fill](0,0)circle(0.5);\draw[lightgray,ultra thick](0,0)circle(3.0);\node at(1.08,2.80){\(\color{blue}{\scriptstyle 12}\)};\node at(0.47,2.96){\(\color{blue}{\scriptstyle 10}\)};\node at(2.33,1.89){\(\color{blue}{\scriptstyle 10}\)};\node at(1.89,2.33){\(\color{red}{\scriptstyle 5}\)};\node at(2.96,0.47){\(\color{red}{\scriptstyle 1}\)};\node at(2.80,1.08){\(\color{red}{\scriptstyle 1}\)};\node at(2.80,-1.08){\(\color{blue}{\scriptstyle 11}\)};\node at(2.96,-0.47){\(\color{blue}{\scriptstyle 11}\)};\node at(1.89,-2.33){\(\color{blue}{\scriptstyle 9}\)};\node at(2.33,-1.89){\(\color{red}{\scriptstyle 6}\)};\node at(0.47,-2.96){\(\color{red}{\scriptstyle 2}\)};\node at(1.08,-2.80){\(\color{red}{\scriptstyle 4}\)};\node at(-1.08,-2.80){\(\color{blue}{\scriptstyle 7}\)};\node at(-0.47,-2.96){\(\color{red}{\scriptstyle 2}\)};\node at(-2.33,-1.89){\(\color{red}{\scriptstyle 6}\)};\node at(-1.89,-2.33){\(\color{blue}{\scriptstyle 7}\)};\node at(-2.96,-0.47){\(\color{blue}{\scriptstyle 8}\)};\node at(-2.80,-1.08){\(\color{blue}{\scriptstyle 8}\)};\node at(-2.80,1.08){\(\color{red}{\scriptstyle 3}\)};\node at(-2.96,0.47){\(\color{blue}{\scriptstyle 9}\)};\node at(-1.89,2.33){\(\color{red}{\scriptstyle 5}\)};\node at(-2.33,1.89){\(\color{red}{\scriptstyle 3}\)};\node at(-0.47,2.96){\(\color{red}{\scriptstyle 4}\)};\node at(-1.08,2.80){\(\color{blue}{\scriptstyle 12}\)};\draw[thick](0.00,5.30)to(0.00,0.70);\draw[dashed](2.65,4.59)to(0.35,0.61);\draw[dashed](4.59,2.65)to(0.61,0.35);\draw[dashed](5.30,0.00)to(0.70,0.00);\draw[dashed](4.59,-2.65)to(0.61,-0.35);\draw[dashed](2.65,-4.59)to(0.35,-0.61);\draw[dashed](0.00,-5.30)to(0.00,-0.70);\draw[dashed](-2.65,-4.59)to(-0.35,-0.61);\draw[dashed](-4.59,-2.65)to(-0.61,-0.35);\draw[dashed](-5.30,0.00)to(-0.70,0.00);\draw[dashed](-4.59,2.65)to(-0.61,0.35);\draw[dashed](-2.65,4.59)to(-0.35,0.61);\draw[blue,thick](-0.90,2.33)arc(111:69:2.50);\draw[blue,thick](1.00,2.61)to(0.90,2.33);\draw[blue,thick](-1.00,2.61)to(-0.90,2.33);\draw[blue,thick](0.30,1.88)arc(81:39:1.90);\draw[blue,thick](2.18,1.76)to(1.48,1.20);\draw[blue,thick](0.44,2.77)to(0.30,1.88);\draw[red,thick](3.27,1.25)arc(21:9:3.50);\draw[red,thick](3.16,0.50)to(3.46,0.55);\draw[red,thick](2.99,1.15)to(3.27,1.25);\draw[blue,thick](2.17,-0.34)arc(-9:-21:2.20);\draw[blue,thick](2.61,-1.00)to(2.05,-0.79);\draw[blue,thick](2.77,-0.44)to(2.17,-0.34);\draw[blue,thick](-1.58,0.25)arc(171:-51:1.60);\draw[blue,thick](1.76,-2.18)to(1.01,-1.24);\draw[blue,thick](-2.77,0.44)to(-1.58,0.25);\draw[red,thick](-0.59,-3.75)arc(-99:-81:3.80);\draw[red,thick](0.50,-3.16)to(0.59,-3.75);\draw[red,thick](-0.50,-3.16)to(-0.59,-3.75);\draw[blue,thick](-0.63,-0.78)arc(231:-111:1.00);\draw[blue,thick](-1.00,-2.61)to(-0.36,-0.93);\draw[blue,thick](-1.76,-2.18)to(-0.63,-0.78);\draw[red,thick](3.89,-3.15)arc(-399:-141:5.00);\draw[red,thick](-2.49,-2.01)to(-3.89,-3.15);\draw[red,thick](2.49,-2.01)to(3.89,-3.15);\draw[blue,thick](-1.21,-0.47)arc(-159:-171:1.30);\draw[blue,thick](-2.77,-0.44)to(-1.28,-0.20);\draw[blue,thick](-2.61,-1.00)to(-1.21,-0.47);\draw[red,thick](-3.19,2.58)arc(-219:-201:4.10);\draw[red,thick](-2.99,1.15)to(-3.83,1.47);\draw[red,thick](-2.49,2.01)to(-3.19,2.58);\draw[red,thick](2.96,3.65)arc(-309:-231:4.70);\draw[red,thick](-2.01,2.49)to(-2.96,3.65);\draw[red,thick](2.01,2.49)to(2.96,3.65);\draw[red,thick](1.58,-4.11)arc(-429:-261:4.40);\draw[red,thick](-0.50,3.16)to(-0.69,4.35);\draw[red,thick](1.15,-2.99)to(1.58,-4.11);\draw[red,very thick,->](0.0000,4.4000)to +(0.001000,0.000000);\draw[red,very thick,->](0.0000,4.7000)to +(0.001000,0.000000);\draw[red,very thick,->](0.0000,5.0000)to +(0.001000,0.000000);\draw[blue,very thick,->](0.0000,1.0000)to +(-0.001000,0.000000);\draw[blue,very thick,->](0.0000,1.6000)to +(-0.001000,0.000000);\draw[blue,very thick,->](0.0000,2.5000)to +(-0.001000,0.000000);\draw[red,very thick,->](2.2000,3.8105)to +(0.000866,-0.000500);\draw[red,very thick,->](2.3500,4.0703)to +(0.000866,-0.000500);\draw[red,very thick,->](2.5000,4.3301)to +(0.000866,-0.000500);\draw[blue,very thick,->](0.5000,0.8660)to +(-0.000866,0.000500);\draw[blue,very thick,->](0.8000,1.3856)to +(-0.000866,0.000500);\draw[blue,very thick,->](0.9500,1.6454)to +(-0.000866,0.000500);\draw[red,very thick,->](3.8105,2.2000)to +(0.000500,-0.000866);\draw[red,very thick,->](4.3301,2.5000)to +(0.000500,-0.000866);\draw[blue,very thick,->](0.8660,0.5000)to +(-0.000500,0.000866);\draw[blue,very thick,->](1.3856,0.8000)to +(-0.000500,0.000866);\draw[red,very thick,->](4.4000,0.0000)to +(0.000000,-0.001000);\draw[red,very thick,->](5.0000,0.0000)to +(0.000000,-0.001000);\draw[blue,very thick,->](1.0000,0.0000)to +(0.000000,0.001000);\draw[blue,very thick,->](1.6000,0.0000)to +(0.000000,0.001000);\draw[red,very thick,->](3.8105,-2.2000)to +(-0.000500,-0.000866);\draw[red,very thick,->](4.3301,-2.5000)to +(-0.000500,-0.000866);\draw[blue,very thick,->](0.8660,-0.5000)to +(0.000500,0.000866);\draw[blue,very thick,->](1.3856,-0.8000)to +(0.000500,0.000866);\draw[red,very thick,->](2.2000,-3.8105)to +(-0.000866,-0.000500);\draw[blue,very thick,->](0.5000,-0.8660)to +(0.000866,0.000500);\draw[red,very thick,->](0.0000,-3.8000)to +(-0.001000,0.000000);\draw[blue,very thick,->](0.0000,-1.0000)to +(0.001000,0.000000);\draw[red,very thick,->](-4.3301,-2.5000)to +(-0.000500,0.000866);\draw[blue,very thick,->](-0.8660,-0.5000)to +(0.000500,-0.000866);\draw[red,very thick,->](-5.0000,0.0000)to +(0.000000,0.001000);\draw[blue,very thick,->](-1.0000,0.0000)to +(0.000000,-0.001000);\draw[red,very thick,->](-3.5507,2.0500)to +(0.000500,0.000866);\draw[red,very thick,->](-4.3301,2.5000)to +(0.000500,0.000866);\draw[blue,very thick,->](-0.8660,0.5000)to +(-0.000500,-0.000866);\draw[blue,very thick,->](-1.3856,0.8000)to +(-0.000500,-0.000866);\draw[red,very thick,->](-2.3500,4.0703)to +(0.000866,0.000500);\draw[red,very thick,->](-2.5000,4.3301)to +(0.000866,0.000500);\draw[blue,very thick,->](-0.5000,0.8660)to +(-0.000866,-0.000500);\draw[blue,very thick,->](-0.8000,1.3856)to +(-0.000866,-0.000500);\draw[red,very thick,->](0.0000,4.4000)to +(0.001000,0.000000);\draw[red,very thick,->](0.0000,4.7000)to +(0.001000,0.000000);\draw[red,very thick,->](0.0000,5.0000)to +(0.001000,0.000000);\draw[blue,very thick,->](0.0000,1.0000)to +(-0.001000,0.000000);\draw[blue,very thick,->](0.0000,1.6000)to +(-0.001000,0.000000);\draw[blue,very thick,->](0.0000,2.5000)to +(-0.001000,0.000000);
\node at(1.42,5.31){\(\scriptstyle \beta\)};\node at(3.89,3.89){\(\scriptstyle \beta\)};\node at(5.31,1.42){\(\scriptstyle y+t_{3}\)};\node at(5.31,-1.42){\(\scriptstyle 1\)};\node at(3.89,-3.89){\(\scriptstyle \beta\)};\node at(1.42,-5.31){\(\scriptstyle 1+\beta t_{6}\)};\node at(-1.42,-5.31){\(\scriptstyle 1\)};\node at(-3.89,-3.89){\(\scriptstyle 1+\beta t_{8}\)};\node at(-5.31,-1.42){\(\scriptstyle 1\)};\node at(-5.31,1.42){\(\scriptstyle 1+\beta t_{10}\)};\node at(-3.89,3.89){\(\scriptstyle 1-\beta y\)};\node at(-1.42,5.31){\(\scriptstyle 1+\beta t_{12}\)};\node at(0.00,5.60){\(\scriptstyle q_{\{4,5,6\},\{7,9,12\}}\)};
\end{tikzpicture}\vspace{-.4cm}
\caption{A storm diagram}\label{fig:eye612}
\end{figure}
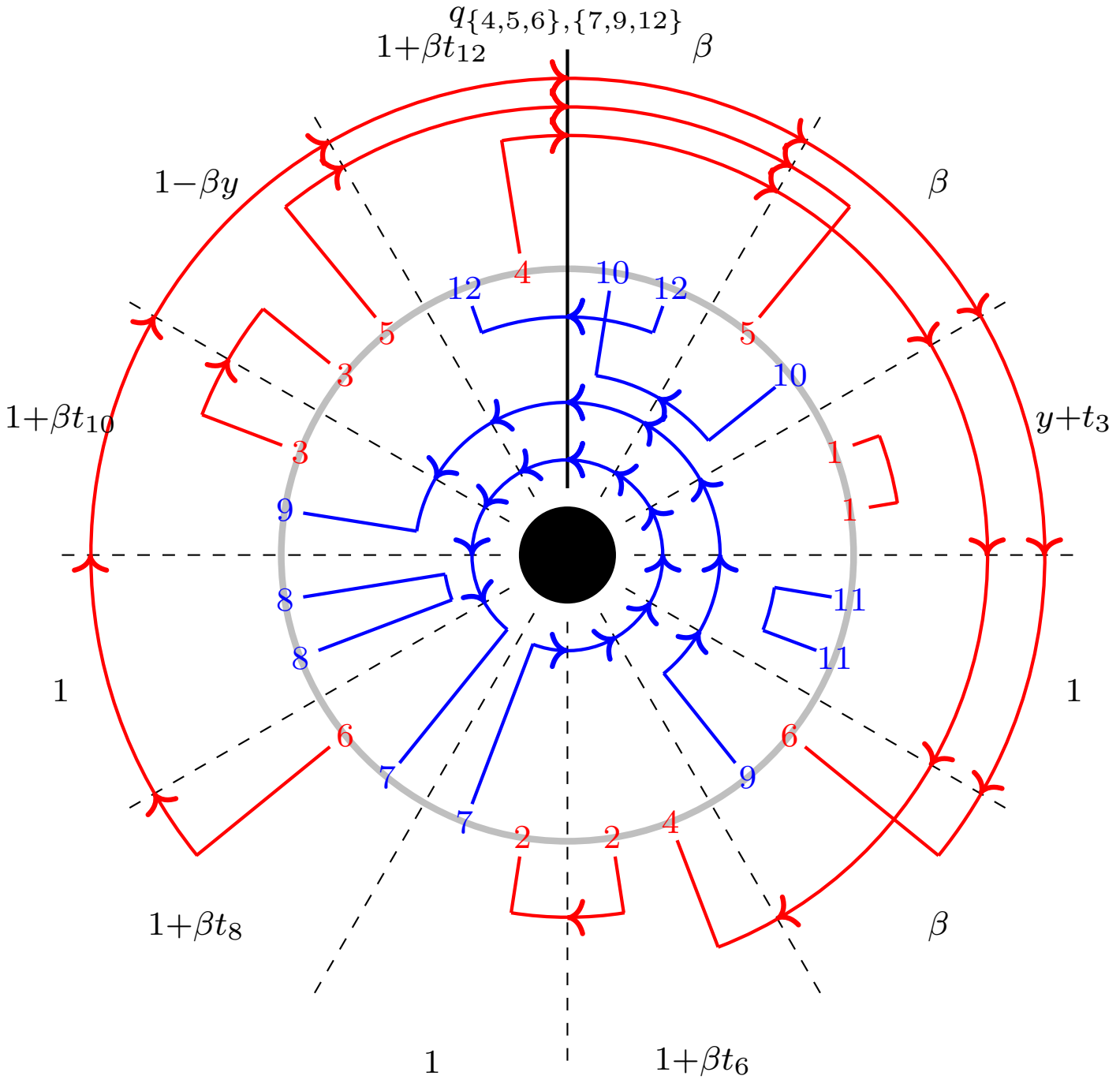

Let $k\in\{1, \cdots, n\}$ and  $u,w\in S_n$. We construct a diagram $\mathcal{D}_k(u, w)$ on the cylinder
$(\mathbb R/n\mathbb Z)\times\mathbb R$,  as follows.
\begin{itemize}
\item Divide the horizontal circle into $n$ intervals. In interval $i$,
place $w^{-1}(i)$ on the left and $u^{-1}(i)$ on the right in the clockwise order.
\item If $u(j)=w(j)$, join the two occurrences
of $j$ inside their common interval, 
at height $j$ if $1\leq j\leq k$, and at depth $n+1-j$ if $k<j\leq n$.
\item For $1\leq a\leq k$ with $u(a)\neq w(a)$, draw a path from the $a$ in interval $u(a)$ to the $a$ in interval $w(a)$, at height $a$, running to the right; 
\item For $k<b\leq n$ with $u(b)\neq w(b)$, draw a path from the $b$ in interval $u(b)$ to
the $b$ in interval $w(b)$, at depth $n+1-b$, running to the left.
\end{itemize}
In particular, no path makes a full circle.  See Figure \ref{fig:eye612} for an example. We call  $\mathcal{D}_k(u, w)$ a \emph{storm diagram}, inspired by its resemblance
to the eye of a storm. 

\vspace{-0.1cm}

\begin{figure}[hhhh]
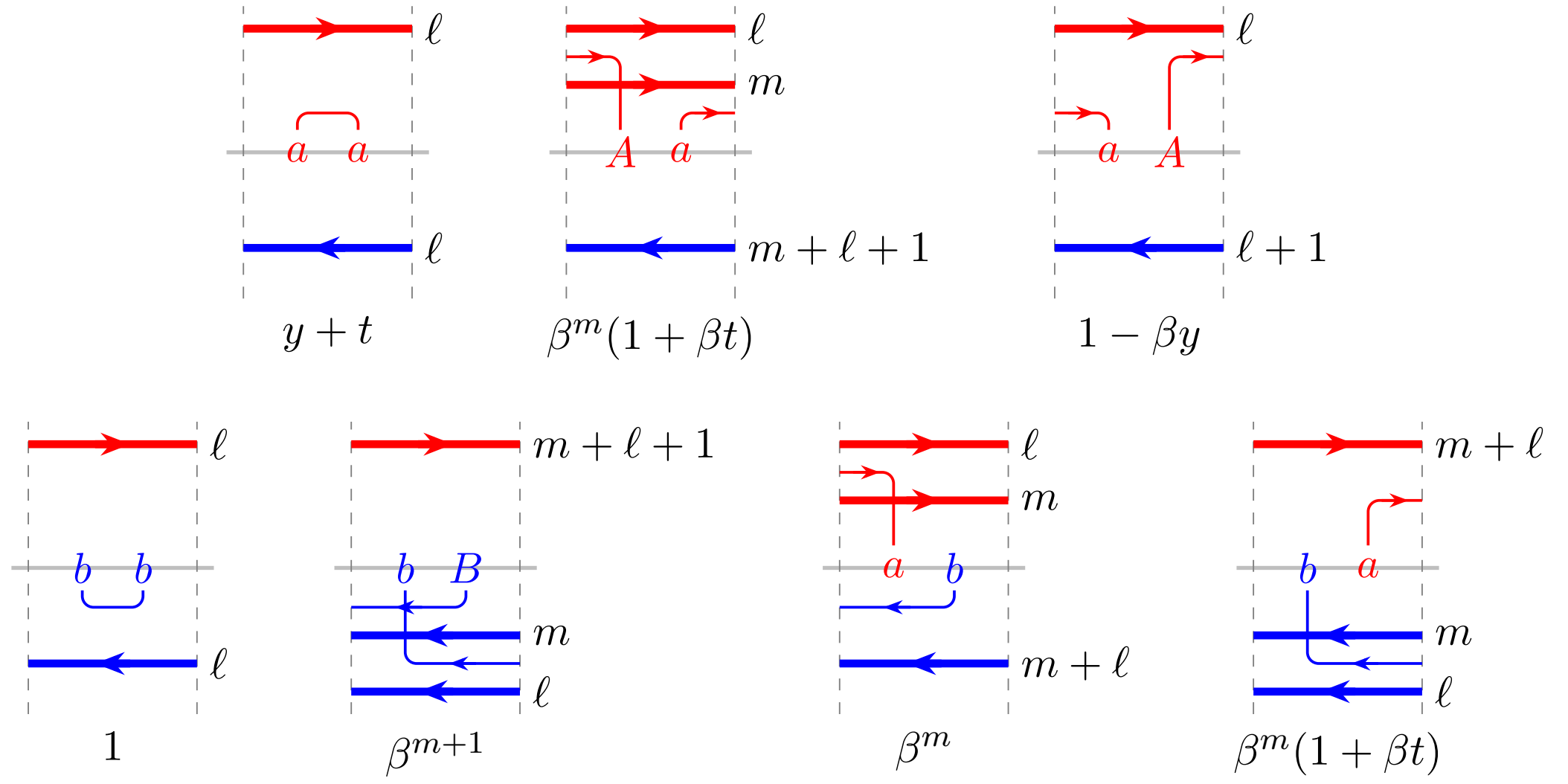

\centering
$\begin{matrix}
\begin{eye}[interval width=15mm,level distance=2.5mm]{4}{8}{1}
  \wind{a}{a}{1}{1}{1}
  \wind[thick]{}{\ell}{4}{.5}{1.5}
  \wind[thick]{}{\ell}{6}{.5}{1.5}
  \label{1}{\(y+t\)}; 
\end{eye}\qquad 
\begin{eye}[interval width=15mm,level distance=2.5mm]{4}{8}{1}
  \wind{}{A}{3}{.5}{1}
  \wind{a}{}{1}{1}{1.5}
  \wind[thick]{}{m}{2}{.5}{1.5}
  \wind[thick]{}{\ell}{4}{.5}{1.5}
  \wind[thick]{m+\ell+1}{}{6}{1.5}{.5}
  \label{1}{\(\beta^m(1+\beta t)\)}; 
\end{eye}\qquad 
\begin{eye}[interval width=15mm,level distance=2.5mm]{4}{8}{1}
  \wind{}{a}{1}{.5}{1}
  \wind{A}{}{3}{1}{1.5}
  \wind[thick]{}{\ell}{4}{.5}{1.5}
  \wind[thick]{\ell+1}{}{6}{1.5}{.5}
  \label{1}{\(1-\beta y\)}; 
\end{eye}\\\\
\begin{eye}[interval width=15mm,level distance=2.5mm]{4}{8}{1}
  \wind{b}{b}{8}{1}{1}
  \wind[thick]{}{\ell}{4}{.5}{1.5}
  \wind[thick]{}{\ell}{6}{.5}{1.5}
  \label{1}{\(1\)}; 
\end{eye}\qquad 
\begin{eye}[interval width=15mm,level distance=2.5mm]{4}{8}{1}
  \wind[thick]{}{m+\ell+1}{4}{.5}{1.5}
  \wind[thick]{\ell}{}{5}{1.5}{.5}
  \wind[thick]{m}{}{7}{1.5}{.5}
  \wind{}{b}{6}{1.5}{1}
  \wind{B}{}{8}{1}{.5}
  \label{1}{\(\beta^{m+1}\)}; 
\end{eye}\qquad 
\begin{eye}[interval width=15mm,level distance=2.5mm]{4}{8}{1}
  \wind{}{a}{3}{.5}{1}
  \wind[thick]{}{m}{2}{.5}{1.5}
  \wind[thick]{}{\ell}{4}{.5}{1.5}
  \wind{b}{}{8}{1}{.5}
  \wind[thick]{m+\ell}{}{6}{1.5}{.5}
  \label{1}{\(\beta^m\)}; 
\end{eye}\qquad 
\begin{eye}[interval width=15mm,level distance=2.5mm]{4}{8}{1}
  \wind[thick]{}{m+\ell}{4}{.5}{1.5}
  \wind[thick]{\ell}{}{5}{1.5}{.5}
  \wind[thick]{m}{}{7}{1.5}{.5}
  \wind{}{b}{6}{1.5}{1}
  \wind{a}{}{2}{1}{1.5}
  \label{1}{\(\beta^{m}(1+\beta t)\)}; 
\end{eye}
\end{matrix}
$
\caption{Local configurations with nonzero weights \,\,($a<A\leq k<b<B$)}
\label{fig:local-configurations}
\end{figure}

For each $1\leq i\leq n$, we define the local weight $c^{(i)}_{u,w}(t_i,y)$ of the storm diagram $\mathcal{D}_k(u, w)$ at  interval $i$ to be the corresponding quantity in Figure \ref{fig:local-configurations}  with $t=t_i$ if the local configuration is given by one of the  seven pictures therein, or zero otherwise. See also Table \ref{tab:local-rules} for a precise description. 
In those pictures, $a<A\leq k<b<B$ whenever these labels occur.
The label beside a thick arrow gives the number of parallel paths,
where $m,\ell\geq0$.  
We retain $\beta$ in the pictures and set $\beta=-1$ throughout unless otherwise stated.

Let $I^{(i)}=I^{(i)}_k(u, w)$ (resp. $J^{(i)}=J^{(i)}_k(u, w)$) denote the set of labels $p$ of paths  crossing the right wall of interval $i$ with $p\leq k$ (resp. $p>k$). In particular, we have  
$$
 I:=I^{(n)}=\{a\le k:w(a)<u(a)\},\qquad
 J:=J^{(n)}=\{b>k:u(b)<w(b)\}.
$$
Define the Laurent monomial $q_{IJ}:= {\prod_{a\in I}q_aq_{a+1}\cdots q_{n-1}\over \prod_{b\in J}q_bq_{b+1}\cdots q_{n-1}} \in \mathbb{Z}[q_1^{\pm 1}, \cdots, q_{n-1}^{\pm 1}]$, where empty products are understood to be 1; in particular,   $q_{IJ}=1$ if $I=J=\varnothing$. We then define the total weight of $\mathcal{D}_k(u, w)$ to be 
\begin{equation}\label{eqn:totalweight}
    F_{u, w}^q:=q_{IJ}\prod_{i=1}^nc_{u, w}^{(i)}(t_i, y).
\end{equation}
We call the right wall of interval $i$  \emph{balanced} 
if $|I^{(i)}|=|J^{(i)}|$. A storm diagram is \emph{balanced} if every  wall   is balanced. Note that all the seven local configurations in Figure \ref{fig:local-configurations} have balanced right walls. It follows that $F_{u, w}^q$ is nonzero only if $\mathcal{D}_k(u, w)$ is a balanced diagram. Furthermore when this holds, we have $|I|=|J|$ and $q_{IJ}\in \mathbb{Z}[q_1, \cdots, q_{n-1}]$.

Note $G/B=\{V_1\subset\cdots\subset V_{n-1}\subset \mathbb{C}^n: \dim V_j=j\}=:F\ell_n$.  Let $\mathcal{V}_k$ denote the $k$-th tautological subbundle over $G/B$, whose fiber at $V_\bullet$ is given by the vector subspace $V_k$. 
Consider the modified  Hirzebruch $\lambda_y$-class $$\lambda_y(\mathcal{V}_k)
:=\sum_{j=0}^k (-1)^j(1+y)^{k-j}[\Lambda^j\mathcal{V}_k]
\in K_T(G/B)[y].$$ 
We obtain the following main result.

\begin{Thm}\label{thm:main-formula}
Let $k\in \{1, \cdots, n\}$ and  $u\in S_n$. Write 
\begin{equation}\label{eq:intro-expansion}
\calO_u\star\lambda_y(\mathcal V_k)
=\sum_{w\in S_n}c^q_{u,w}(t,y)  \calO_w
\quad\text{in }QK_T(G/B)[y].
\end{equation}
Then we have \begin{equation}\label{eq:intro-main-formula}
c^q_{u,w}(t,y)
=F_{u, w}^q=q_{IJ}\prod_{i=1}^n c^{(i)}_{u,w}(t_i,y).
\end{equation}
\end{Thm}

To prove the above theorem, we show that the coefficients $c_{u, w}^q(t, y)$ satisfy certain recursive relationships as in Propositions \ref{thm:clsrec} and \ref{qmtindonw} and initial conditions as in Propositions \ref{prop:clasinitial} and \ref{prop:quaninitial}, by using the left Weyl group action and the left Demazure operators \cite{MNS} on $QK_T(G/B)$, together with one specific Demazure operator in Theorem \ref{thm:full-flag-Demazure} defined with the help of the $K$--theoretic Peterson isomorphism. Then we achieve the  goal, by showing that the total weights $F_{u, w}^q$ of the storm diagrams satisfy the same recursive relations and initial conditions in a purely combinatorial way.

\begin{Eg}

For $n=12$ and $k=6$, we consider
$u=[3,6,10,12,11,8,7,9,5,2,4,1]$ and $ w=[3,7,11,6,2,5,8,9,10,1,4,12]$ in one-line notation.
We   draw the storm diagram in Figure \ref{fig:eye612}, where  
the intervals are numbered clockwise from the marked cut at the top. 

After cutting the storm diagram  open along the wall between intervals 12 and 1, we obtain Figure \ref{fig:linediagram}. For each
label $i$, the path starts in interval $u(i)$ and ends in interval $w(i)$.

\begin{figure}[ht]
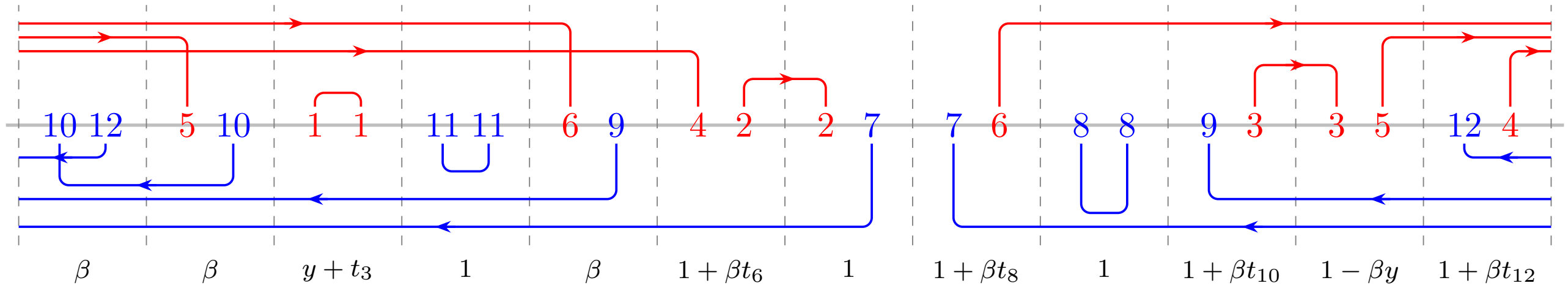

\begin{eye}[interval width=0.08\linewidth]{6}{12}{12}
\wind{1}{1}{1}{3}{3}
\wind{2}{2}{2}{6}{7}
\wind{3}{3}{3}{10}{11}
\wind{4}{}{4}{12}{12.5}
\wind{}{4}{4}{.5}{6}
\wind{5}{}{5}{11}{12.5}
\wind{}{5}{5}{.5}{2}
\wind{6}{}{6}{8}{12.5}
\wind{}{6}{6}{.5}{5}
\wind{7}{}{7}{7}{0.5}
\wind{}{7}{7}{12.5}{8}
\wind{8}{8}{8}{9}{9}
\wind{9}{}{9}{5}{0.5}
\wind{}{9}{9}{12.5}{10}
\wind{10}{10}{10}{2}{1}
\wind{11}{11}{11}{4}{4}
\wind{12}{}{12}{1}{0.5}
\wind{}{12}{12}{12.5}{12}
\label{1}{\tiny $\beta$}
\label{2}{\tiny $\beta$}
\label{3}{\tiny $y+t_3$}
\label{4}{\tiny $1$}
\label{5}{\tiny $\beta$}
\label{6}{\tiny $1+\beta t_6$}
\label{7}{\tiny $1$}
\label{8}{\tiny $1+\beta t_8$}
\label{9}{\tiny $1$}
\label{10}{\tiny $1+\beta t_{10}$}
\label{11}{\tiny $1-\beta y$}
\label{12}{\tiny $1+\beta t_{12}$}
\end{eye}
\caption{The storm diagram of
Figure \ref{fig:eye612}, cut open between intervals $n$ and 1}\label{fig:linediagram}
\end{figure}
Here the 1st interval is the 5th local configuration  above with $m=0$ and $\ell=2$, so that the resulting weight is $\beta$;  paths $p$ with $p$ in $I=\{4,5,6\}$ or $J=\{7,9,12\}$ cross the cut.
Thus 
$q_{IJ}=q_4q_5^2q_6^3 q_7^2q_8^2 q_9q_{10}q_{11}$ and 
\begin{align*}
c^q_{u,w}(t,y)
={}&-q_{IJ}(y+t_3)(1+y)(1-t_6)(1-t_8)(1-t_{10})(1-t_{12}).
\end{align*}
\end{Eg}

\subsection{Applications}
As one main application of  Theorem \ref{thm:main-formula}, we obtain a (manifestly positive) Pieri rule in $QK_T(G/B)$. 
Recently,  Gorbounov,   Korff and Mihalcea also used Hirzebruch $\lambda_y$-classes to study the equivariant quantum $K$-theory $QK_T(Gr(m, n))$ of complex Grassmannians  \cite[Theorem 1.1]{GKM25}.   
Building on this work, they obtained explicit combinatorial formulae
for the Schubert structure  coefficients using partition functions
and cylindric Grothendieck polynomials \cite[Theorems~6.2, 6.3, and~7.1]{GKM}, while these  formulae are not cancellation-free. Note that  
the positivity conjecture was proved for the non-equivariant  case $QK(Gr(m, n))$ in \cite{BCMP22}.

In our case, a key point is that the total weight of the storm diagram $\mathcal{D}_k(u, w)$, after taking  specializations of the form   $y=[\mathcal O_{\mathbb{P}^{i-1}}(1)]-1\in K_{T}(\mathbb{P}^{i-1})$, can be related with the localization type coefficients $$N_{s_{r-b+1}s_{r-b+2}\cdots s_r, w_A}^{w_A}=\mathcal{O}_{ s_{r-b+1}s_{r-b+2}\cdots s_r}|_{w_A}=:(-1)^b\psi_b(A).$$
Here $w_A\in S_n$ is a Grassmannian permutation for $Gr(r, n)$, and $A=\{w_A(1), w_A(2), \cdots, w_A(r)\}$. In particular, there have been various known manifestly positive formulae of $\psi_b(A)$  \cite{KoKu, McN, GrKr}. For example,  one such positive formula is given in  Remark  \ref{rmk:loccoef}.

Consider    the following sets depending on $u, w$ and $k$:
$$
\mathsf F=\{j\in u([k]):u^{-1}(j)=w^{-1}(j)\}, \qquad  \mathsf U=\{j\in u([k]):u^{-1}(j)<w^{-1}(j)\}.
$$    
Note  that   $\deg q_j=2$ for all $j$, and we denote by $\{\alpha_1, \cdots, \alpha_{n-1}\}$ the standard base of simple roots of $G$.
Using Theorem \ref{thm:main-formula}, we establish the following Pieri rule in \textbf{Theorem   \ref{thm:equivPieri}}.

\begin{Th}\label{thm:equivPieriINTR}
Let $u,w\in S_n$, $1\leq i\leq k<n$, and $d\in H_2(G/B, \mathbb{Z})$. Then   $N_{u, s_{i}s_{i+1}\cdots s_k}^{w,d}=0$ unless  $c_{u, w}^q(t, y)\neq 0$ and $q^d=q_{IJ}$.  
When both conditions hold, we have
\begin{equation*} 
 (-1)^{\ell(u)+k-i+1-\ell(w)}
 N_{u,s_{i}s_{i+1}\cdots s_k}^{w,d}
 =\binom{|\mathsf U|-1}{i-1-|\mathsf F|}
+ \sum_{\substack{S\subseteq\mathsf U~:~|S|\ge i-|\mathsf F|}}
\psi_{|\mathsf F|+|S|-i+1}({\mathsf F\cup S}).
\end{equation*}
In particular, $ (-1)^{\ell(u)+k-i+1+\ell(w)+\deg q^d}
 N_{u,s_{i}s_{i+1}\cdots s_k}^{w,d}\in \mathbb{Z}_{\ge 0}[e^{-\alpha_1}-1, \cdots, e^{-\alpha_{n-1}}-1]$.
\end{Th}
\noindent We learned through personal communication that   Naito, Nguyen and Sagaki have independently obtained a related Pieri rule for $QK_T(G/B)$ \cite{NNS}. Their approach combines the general Chevalley formula for semi-infinite flag manifolds \cite[Theorem~33]{LNS} with the non-equivariant quantum Pieri rule \cite[Theorem~2.13]{NaSa}.

 Seidel \cite{Sei} constructed certain   classes in quantum
cohomology.  
The analogue for the quantum $K$-theory of      
(co)minuscule Grassmannians was studied by Buch, Chaput and Perrin    \cite{BCP}, and was further extended to the equivariant setting   \cite{BCP0} (see also \cite{GKM25} for complex Grassmannians).   The analogue for the equivariant quantum $K$-theory of flag varieties of    general Lie type was  proved by Ikeda, Kouno and Naito   \cite{IKN}. In type $A$, we  provide an alternative proof of the equivariant Seidel product for $QK_T(G/B)$  in  \textbf{Corollary \ref{cor:equivSeidel}} using Theorem \ref{thm:main-formula}. 
A typical application of the Seidel product in quantum Schubert calculus is of the  form:  
 $$\mbox{quantum Pieri rule } = \mbox{ classical Pieri rule } + \mbox{ Seidel action}.$$
    This has been shown for the quantum $K$-theory of complex Grassmannians in \cite{LLSY} and minuscule Grassmannians of classical Lie type in \cite{BCP}. It was also discussed for $QK(F\ell_n)$   
          in \cite[Conjecture 3.1 and Problem 3.2]{LiSo}, which can be resolved by the non-equivariant specialization of Theorem \ref{thm:parttofullINT}.

Consider the partial flag variety $G/P=\{V_{n_1}\subset \cdots\subset V_{n_r}\subset \mathbb{C}^n: \dim V_{n_j}=n_j\}=:F\ell_{n_1, \cdots, n_r; n}$. The natural projection $\pi: G/B=F\ell_n\longrightarrow F\ell_{n_1, \cdots, n_r; n}=G/P$ induces a surjective ring homomorphism 
 $\Phi: QK_T(G/B)\to QK_T(G/P)$ by Kato \cite[Theorem 2.19]{Kato3}. Note that
 the ring homomorphism $\Phi$ does not induce a positive Pieri rule for $QK_T(G/P)$ directly. Instead, 
by applying the equivariant Seidel product  twice, we obtain the following     key property, and thereby establish a positive Pieri rule in \textbf{Theorem \ref{thm:equiPierPara}}. Here $\hat u, \hat w\in S_n$ and $\delta_k(u, w)\in H_2(G/B, \mathbb{Z})$ are associated elements depending only on $u, w, k$ (see \textbf{Theorem \ref{thm:parttofull}} for precise descriptions of all the relevant notations).  

\begin{Th}\label{thm:parttofullINT} 
Let $k\in \{n_1, \cdots, n_r\}$, $u, w\in W^P$, and  $1\leq i\leq k$. In $QK_T(G/P)$, we have 
 $$ \sum_{d\in H_2(G/P, \mathbb{Z})}\bar N_{u, s_is_{i+1}\cdots s_k}^{w, d}Q^d= N_{\hat u, s_is_{i+1}\cdots s_k}^{\hat w, \delta_k(u, w)} Q^{\pi_*(\delta_k(u, w))}. $$ 
\end{Th}
\noindent The above theorem may be viewed as a $K$-theoretic analogue of Peterson's comparison formula \cite{Wood, HL15} for quantum cohomology, in the case of Pieri-type coefficients.

\subsection*{History and organization of the paper}
This project was started in 2024. While preparing \cite{FGX}, the first,
second, and last authors observed that the configurations in each row
of the puzzle rule in that paper correspond bijectively to those in a special case of the equivariant $K$-theoretic Pieri rule.
We then asked whether the general Pieri rule could be described in a similar way. 
Motivated by the quantum Pieri rule for Grassmannians \cite{Pos05-2}, we further conjectured that placing the same model on a cylinder would give a quantum extension.
Together with explorations in various computational examples by SageMath, this led us to introduce the combinatorial model of storm diagrams with   the local weights defined by Figure \ref{fig:local-configurations}.

To prove Theorem \ref{thm:main-formula}, we   noticed  the key Proposition \ref{prop:two-column},  which would imply the desired recursions for the total weights of storm diagrams and could be proved   combinatorially.   SageMath computations for $(k,n)=(5,10)$ led us to essentially all the cases listed in Appendix~\ref{app:two-column}.
We believed that the list was exhaustive, but did not have enough patience to confirm it.
If we had  infinite patience and never become bored, we could have finished this paper before the AI era. Now, with AI’s help, we finally can.

The paper is organized as follows. In Section 2, we review the necessary background. In Section 3, we show that the structure coefficients $c^q_{u, w}(t, y)$ satisfy certain recursive relations and initial conditions. 
In Section 4, we prove Theorem \ref{thm:main-formula}, by showing that the total weights $F_{u, w}^q$ of storm diagrams   satisfy the same   recursive relations and initial conditions. In Section 5, we establish Pieri rules for $QK_T(G/P)$ as main applications. 

\subsection*{AI Disclosure}

In early 2026, Gemini completed the case classification for
Proposition~\ref{prop:two-column}. ChatGPT 6 subsequently polished
the proof, now presented in Appendix~\ref{app:two-column}.
It also wrote the \LaTeX{} package used to retypeset all diagrams
except Figure~\ref{fig:eye612}, and assisted in establishing
Lemma~\ref{lemma:Olambda} and proofreading the manuscript.

\subsection*{Acknowledgements}
The authors would like to thank Leonardo C. Mihalcea and Satoshi Naito for helpful discussions.
Fan is supported by the National Key R\&D
Program of China (No. 2025YFA1017702) and the National Natural Science Foundation
of China (No. 12471314).   Guo is supported by the National Natural Science Foundation of China (No.
12371329)  and the Fundamental Research Funds for the Central Universities (No.
63263094).

\section{Preliminaries}

In this section, we review the necessary background.
\subsection{Equivariant quantum $K$-theory of flag varieties}
Let us denote subgroups of $GL_n$ as follows
$$
T=\left[\begin{matrix}
* \\
&\ddots\\
&& *
\end{matrix}\right]\subset 
B=\left[\begin{matrix}
* & \cdots&*\\
&\ddots & \vdots\\
&& *
\end{matrix}\right]\subset G=GL_n
\supset
\left[\begin{matrix}
* \\
\vdots&\ddots\\
*&\cdots&*
\end{matrix}\right]=B^-.$$

\subsubsection*{Equivariant parameters} Let $\epsilon_i:T\to \mathbb{C}^\times$ be the $i$-th projection. Then we can identify the weight lattice 
$X^*(T) =\mathbb{Z}\epsilon_1\oplus \cdots\oplus \mathbb{Z}\epsilon_n$, and obtain the set $\{\alpha_i=\epsilon_i-\epsilon_{i+1}\}_{i=1}^{n-1}$ of simple roots. Let $\{\varpi_i\}_i$ denote the fundamental weights, and set $\rho:=\sum_{i=1}^{n-1}\varpi_i$. 
Identifying the $T$-equivariant $K$-theory of a point with the representation ring of $T$ via $[\mathbb{C}_{\epsilon_i}]=e^{\epsilon_i}$ and denoting $t_i:=1-e^{\epsilon_i}$, we have 
$$K_T(\pt)=
\mathbb{Z}[(1-t_1)^{\pm 1},\ldots,(1-t_n)^{\pm 1}].$$

\subsubsection*{Equivariant $K$-theory}
Consider the (complete) flag variety $G/B$,  
$$G/B = \left\{
 0=V_0\subset V_1\subset \cdots \subset V_{n-1}\subset V_n=\mathbb{C}^n: \dim V_j=j\right\}=:F\ell_n.$$
There are  (opposite) Schubert varieties $Y(w)=\overline{B^-wB/B}\subset G/B$ of codimension $\ell(w)$, indexed by permutations $w\in S_n=W$. Here we have taken a canonical isomorphism that sends each transposition $(i, i+1)$ to     the simple reflection $s_i=s_{\alpha_i}$. The Schubert varieties $Y(w)$ are $T$-invariant, and the classes of structure sheaves on them form a basis of the $T$-equivariant $K$-theory of $G/B$:
 $$K_T(G/B)=\bigoplus_{w\in S_n}K_T(\pt)\mathcal{O}_w,\qquad\mbox{where}\quad\mathcal{O}_w:=[\mathcal{O}_{Y(w)}].$$ 
The Schubert structure coefficients $N_{u, v}^w$ satisfy the positivity property \cite{AGM}: for any $u, v\in S_n$
$$\mathcal{O}_u\cdot \mathcal{O}_v=\sum_{w\in S_n}N_{u, v}^w\mathcal{O}_w,\qquad\mbox{with}\quad (-1)^{\ell(u)+\ell(v)+\ell(w)}N_{u, v}^w\in \mathbb{Z}_{\geq 0}[e^{-\alpha_1}-1, \cdots, e^{-\alpha_{n-1}}-1].$$

\subsubsection*{Special classes}
 Denote the $j$-th tautological sub-bundle on $G/B$
by $\mathcal{V}_j$, whose fiber at a point $V_{\bullet}$ is the vector space $V_j$. The sub-bundles $\mathcal{V}_j$ and the line bundles   $\mathcal{L}_j:=\mathcal{V}_j/\mathcal{V}_{j-1}$ are all $G$-bundles on $G/B$. We  consider the $G$-equivariant $K$-theoretic Chern polynomial for the dual bundle (with a convention slightly different from  the Hirzebruch $\lambda_y$-class):  
$$\lambda_y^G(\mathcal{V}_k)=\sum_{r=0}^kc_r^G(\mathcal V_k^\vee)y^{k-r}=\prod_{i=1}^k\,
(y+1-[\mathcal L_i]_G)
=\sum_{r=0}^k (-1)^r(1+y)^{k-r}[\Lambda^r\mathcal{V}_k]
\in K_G(G/B)[y],$$
where   $[\mathcal{L}_1]_G,\ldots,[\mathcal{L}_k]_G$ are the $K$-theoretic Chern roots of $\mathcal{V}_k$. 
Its restriction gives the $T$-equivariant class $\lambda_y(\mathcal{V}_k)\in K_T(G/B)[y]$ in exactly the same form by replacing $G$ with $T$.
We remark that $1-[\mathcal{L}_i]$ is usually denoted by $x_i$ in classical Schubert calculus, so the coefficients of $\lambda_y(\mathcal{V}_k)$ are the elementary symmetric polynomials in $x_1,\ldots,x_k$. 
By \cite[Corollary 7.4]{FGSX}, we have a precise relation between the Chern classes 
$c_r^T(\mathcal V_k^\vee)$ and the special Schubert classes $\mathcal{O}_{s_is_{i+1}\cdots s_k}$, $1\leq i\leq k\leq n-1$, which is equivalently given by 
\begin{equation}\label{eqn:lambday}
\lambda_y(\mathcal V_k)
=\prod_{i=1}^k(y+t_i) +\sum_{i=1}^k(1-t_i)(1+y)^{k-i}
\prod_{a<i}(y+t_a)\,
\mathcal{O}_{s_is_{i+1}\cdots s_k}.
\end{equation}

\subsubsection*{Equivariant quantum $K$-theory}  Consider  the coroot lattice $Q^\vee=\bigoplus_{i=1}^{n-1}\mathbb{Z}\alpha_i^\vee$. For $d=\sum_id_i\alpha_i^\vee$, we denote  the Laurent monomial $q^d:=q_1^{d_1}\cdots q_{n-1}^{d_{n-1}}$  in the Novikov variables $q_i$. For $1\leq a<b\leq n$ and  subsets $I_1, I_2$ of $[n]=\{1, \cdots, n\}$, we define   $q_{aa}:=1$,  
\begin{equation}
  (q_{ba})^{-1}=q_{ab}:=q_aq_{a+1}\cdots q_{b-1}\quad\mbox{ and }\quad  q_{I_1I_2}:={\prod_{a\in I_1} q_{an}\over \prod_{b\in I_2} q_{bn}},
\end{equation}
where empty products are understood to be 1; in particular,   $q_{I_1I_2}=1$ if $I_1=I_2=\varnothing$.

 The (localized) $T$-equivariant quantum $K$-theory of $G/B$ is a deformation of $K_T(G/B)$,  
$$QK_T(G/B) = \big(K_T(G/B)\otimes_{\mathbb{Z}}\mathbb{C}[q_1^{\pm 1},\cdots, q_{n-1}^{\pm 1}], \star\big),$$
with the Schubert structure coefficients $N_{u, v}^{w, d}$ in the quantum product  $$\mathcal{O}_u \star\mathcal{O}_v=\sum_{w\in S_n, d\in Q^\vee} N_{u, v}^{w, d}\mathcal{O}_w {q}^d$$
defined by a combination of genus-zero, two- and three-point $T$-equivariant $K$-theoretic Gromov--Witten invariants   (see \cite{Giv00, Lee, BCMP} for more details). Here  $Q^\vee$ is canonically identified with $H_2(G/B, \mathbb{Z})$; $N_{u, v}^{w, 0}=N_{u, v}^w$;   $N_{u, v}^{w, d}=0$ unless $d$ is effective. Moreover, 
$\sum_{d}N_{u, v}^{w, d}q^d$ is a polynomial in $q_i$ by  the finiteness property proved in \cite{ACT, Kato1}. 
We invert the Novikov variables to use the quantum-affine correspondence in the proof of
Theorem  \ref{thm:full-flag-Demazure}. The following positivity property was predicted in \cite[Conjecture 2.2]{BCMP}:
\begin{equation}\label{eqn:conjPos}
    (-1)^{\ell(u)+\ell(v)+\ell(w)+\deg q^d}N_{u, v}^{w, d}\in \mathbb{Z}_{\geq 0}[e^{-\alpha_1}-1, \cdots, e^{-\alpha_{n-1}}-1]
\end{equation}
    where $\deg q^d=\sum_i 2d_i$ and therefore never contributes a sign in the case of $F\ell_n$.

The  $G$-equivariant quantum  $K$-theory $QK_G(G/B)=\big(K_G(G/B)\otimes_{\mathbb{Z}}\mathbb{C}[q_1^{\pm 1},\cdots, q_{n-1}^{\pm 1}], \star\big)$ is defined similarly, and we can view 
$QK_T(G/B)$ as a $QK_G(G/B)$-module, induced from the $K_G(G/B)$-module structure on $K_T(G/B)$.

\subsection{Affine Weyl group and semi-infinite Bruhat order}
\label{subsec:affine-conventions}
The affine Weyl group is
\(
\widehat{W}=W\ltimes Q^\vee.
\)
That is, every element of $
\widehat{W}$ is uniquely written as $wt_ \mu$, where $w\in W$ and $t_ \mu$ is the symbol for an element
$ \mu\in Q^\vee$ and  
the multiplication is characterized by $w t_{ \mu} w^{-1}=t_{w \mu}$.
Define the affine simple reflection $s_0:=s_{1n}\cdot t_{\epsilon_n-\epsilon_1}$, where $s_{1n}\in S_n=W$ denotes the transposition $(1,n)$. 
We  regard the indices $0,1,\ldots,n-1$ of $s_0,s_1,\ldots,s_{n-1}$ as elements of $\mathbb{Z}/n\mathbb{Z}$. 

We   use the semi-infinite order on $\widehat{W}$,  defined by Lusztig  \cite{Lus80} in order to study the periodic module of affine Hecke algebras. This order  was also studied under the name stable Bruhat order in Peterson's notes \cite{Pet97} and the name quantum Bruhat order in \cite{BFP99,Pos05,LaSh, LNSSS,INS16}. 
The definition we use here is equivalent to the definitions in these references. 

For a coroot $\nu=\sum_{i=1}^{n-1} d_i\alpha_i^\vee\in Q^\vee$, by  
$\nu\gg0$ we mean that $\nu$ is sufficiently dominant, that is, $\langle \alpha_i, \nu\rangle \gg0$ for all $i$. 
The \emph{semi-infinite length} is defined by
$$\ell_{\semi}(wt_ \mu) 
= \ell(wt_{ \mu+ \nu})-\ell(t_ \nu)\qquad \text{when $ \nu \gg0$},$$
where $\ell:\widehat{W}\to \mathbb{Z}_{\geq0}$ is the Coxeter length of $\widehat{W}$.
This is well-defined and can be computed explicitly:  
$\ell_{\semi}(wt_ \mu)
   =\ell(w)+\langle 2\rho, \mu\rangle$. 
The semi-infinite order $\leq_{\semi}$ on $\widehat{W}$ can be characterized by 
$$
wt_ \mu \leq_{\semi} w't_{ \mu'}\iff 
wt_{ \mu+ \nu} \leq w't_{ \mu'+ \nu}$$
where $\nu\gg 0$ and $\leq$ is the Bruhat order on $\widehat{W}$. 
In this paper, we will only use the following case.
For $i\in \mathbb{Z}/n\mathbb{Z}$
$$s_iwt_ \mu >_{\semi} wt_ \mu
\iff \ell_{\semi}(s_iwt_ \mu)=\ell_{\semi}(wt_ \mu)+1.$$
One can check the following
\[
 \begin{aligned}
 s_iwt_ \mu>_{\semi}wt_ \mu
   &\iff s_iw>w\iff w^{-1}(i)<w^{-1}(i+1)
       &&(1\leq i<n),\\
 s_0wt_ \mu>_{\semi}wt_ \mu
   &\iff s_{1n}w<w\iff w^{-1}(n)<w^{-1}(1).
 \end{aligned} 
\]

\begin{Lemma}[{\cite[Lemma 6.12]{LNSSS}}]\label{lem:chain}
For any $u,w\in W$, there exist $i_1,\ldots,i_\ell\in \mathbb{Z}/n\mathbb{Z}$ and $ \mu\in Q^\vee$ such that 
$$u<_\semi 
s_{i_1}u<_\semi
s_{i_2}s_{i_1}u<_\semi\cdots<_\semi
s_{i_\ell}\cdots s_{i_2}s_{i_1}u=wt_{ \mu}.$$
\end{Lemma}

\section{Recursions for the coefficients $c_{u, w}^q$}

\subsection{Classical coefficients}
We regard \(K_T(G/B)\) as a \(K_G(G/B)\)-algebra via the natural induced  homomorphism \(K_G(G/B)\to K_T(G/B)\). All products involving classes from these two rings are understood in this sense. In  $K_T(G/B)[y]$, we write 
\begin{equation}\label{eq:clscuw}
\calO_{u}\cdot \lambda_y(\mathcal{V}_k)
=\sum_{w\in S_n} c_{u,w}  \calO_{w},\qquad c_{u,w}\in K_T(\pt)[y].
\end{equation}
Following \cite[equation (10)]{MNS},  the \emph{left Demazure operator} on $K_T(G/B)$ is defined as below.  For $1\leq i\leq n-1$, 
\begin{align}
D_i^\vee
=-\frac{1-t_{i}}{t_i-t_{i+1}}\left(1-\frac{1-t_{i+1}}{1-t_{i}}s_i^L\right). 
\end{align}
Here $s_i^L$
denotes the left Weyl group action. By \cite[Remark 5.6]{MNS}, we have 
\[D_i^\vee(\calO_{u})
=\begin{cases}
\calO_{s_iu}, & s_iu<u,\\
\calO_{u}, & s_iu>u.
\end{cases}
\]
For $\gamma\in K_T(G/B)$ and $f\in K_T(\pt)$, we have
\[
D_i^\vee(f\cdot \gamma)=f|_{t_i\leftrightarrow t_{i+1}}\cdot D_i^\vee(\gamma)-\frac{1-t_i}{t_i-t_{i+1}}(f-f|_{t_i\leftrightarrow t_{i+1}})\cdot \gamma.
\]
 {The left Weyl-group action fixes every class coming from \(K_G(G/B)\). It follows that $D_i^\vee$ is $K_G(G/B)$-linear.} That is,  for any $\gamma'\in K_G(G/B)\subset K_T(G/B)\ni \gamma$, we have 
$$D_i^\vee(\gamma'\cdot \gamma)=\gamma'\cdot D_i^\vee(\gamma).$$

\begin{Prop}\label{indonw}\label{thm:clsrec}
If $s_iw>w$ for some $1\leq i\leq n-1$, then 
\begin{align}\label{PURW-1}
c_{u,s_iw} & 
=\begin{cases}
  -\dfrac{1-t_{i+1}}{t_i-t_{i+1}}c_{u,w}|_{t_i\leftrightarrow t_{i+1}}
+\dfrac{1-t_{i}}{t_i-t_{i+1}}c_{u,w}
+ c_{s_iu,w}|_{t_i\leftrightarrow t_{i+1}}, 
& s_iu<u,\\[15pt]
  -\dfrac{1-t_{i}}{t_i-t_{i+1}}c_{u,w}|_{t_i\leftrightarrow t_{i+1}}
+\dfrac{1-t_{i}}{t_i-t_{i+1}}c_{u,w},
& s_iu>u. 
\end{cases}
\end{align}
\end{Prop}

\begin{proof}
Note   $s_iw>w\iff w^{-1}(i)<w^{-1}(i+1).$ Applying $D_i^\vee$ to the left-hand side of \eqref{eq:clscuw}, 
we find
$$
D_i^\vee(\calO_{u}\cdot \lambda_y(\mathcal{V}_k))
=D_i^\vee(\calO_{u})\cdot \lambda_y(\mathcal{V}_k)
=\begin{cases}
    \calO_{s_iu}\cdot \lambda_y(\mathcal{V}_k)=\displaystyle{\sum_w} c_{s_iu,w}\cdot \calO_{w}, & s_iu<u,\\[5pt]
    \calO_{u}\cdot \lambda_y(\mathcal{V}_k)=\displaystyle{\sum_w} c_{u, w}\cdot \calO_{w}, & s_iu>u.\\
\end{cases}$$
Applying $D_i^\vee$ to the right-hand side of \eqref{eq:clscuw}, we find
\begin{align*}
& \sum_w D_i^\vee\big(c_{u, w}\cdot \calO_{w}\big)\\
& = \sum_{w}\left(c_{u, w}|_{t_i\leftrightarrow t_{i+1}}\cdot D_i^\vee(\calO_w)
    -\frac{1- t_i}{t_i-t_{i+1}}(c_{u, w}-c_{u, w}|_{t_i\leftrightarrow t_{i+1}})\cdot\calO_w\right) \\
& = \sum_{s_iw<w}\left(c_{u, w}|_{t_i\leftrightarrow t_{i+1}}\cdot\calO_{s_iw}
    -\frac{1- t_i}{t_i-t_{i+1}}(c_{u, w}-c_{u, w}|_{t_i\leftrightarrow t_{i+1}})\cdot\calO_w\right)\\ 
&\quad +\sum_{s_iw>w}\left(c_{u, w}|_{t_i\leftrightarrow t_{i+1}}\cdot\calO_w
    -\frac{1- t_i}{t_i-t_{i+1}}(c_{u, w}-c_{u, w}|_{t_i\leftrightarrow t_{i+1}})\cdot\calO_w\right) \\
    & = -\sum_{s_iw<w} \frac{1- t_i}{t_i-t_{i+1}}(c_{u, w}-c_{u, w}|_{t_i\leftrightarrow t_{i+1}})\cdot\calO_w\\
    &\quad + \sum_{s_iw>w}
    \left(c_{u,s_iw}|_{t_i\leftrightarrow t_{i+1}}
    + c_{u,w}|_{t_i\leftrightarrow t_{i+1}}
    -\frac{1- t_i}{t_i-t_{i+1}}(c_{u,w}-c_{u,w}|_{t_i\leftrightarrow t_{i+1}})\right) \cdot\calO_{w}\\
    & = -\sum_{s_iw<w} \frac{1- t_i}{t_i-t_{i+1}}(c_{u,w}-c_{u,w}|_{t_i\leftrightarrow t_{i+1}})\cdot\calO_{w}\\
    &\quad + \sum_{s_iw>w}
    \left(c_{u,s_iw}|_{t_i\leftrightarrow t_{i+1}}
    -\frac{1-t_i}{t_i-t_{i+1}}c_{u,w}
    +\frac{1-t_{i+1}}{t_i-t_{i+1}}c_{u,w}|_{t_i\leftrightarrow t_{i+1}}\right)\cdot \calO_{w}.
\end{align*}

Extracting  the coefficients of $\calO_w$ with $s_iw>w$ on both sides, we deduce  that 
\begin{equation*}
    c_{u,s_iw}|_{t_i\leftrightarrow t_{i+1}}
=\frac{1- t_i}{t_i-t_{i+1}}c_{u,w}
-\frac{1- t_{i+1}}{t_i-t_{i+1}}c_{u,w}|_{t_i\leftrightarrow t_{i+1}}
+\begin{cases}
c_{s_iu,w}, & s_iu<u,\\[5pt]
c_{u,w}, & s_iu>u,
\end{cases}
\end{equation*}
which coincides with \eqref{PURW-1} after the variable exchange $t_i\leftrightarrow t_{i+1}$. 
\end{proof}

\begin{Prop}\label{prop:clasinitial}
Let $\delta$ denote the  Kronecker symbol. For any $u\in S_n$, we have 
\begin{equation}\label{eq:clsini}
c_{u,\id}(t,y)=\delta_{u,\id}\cdot\prod_{i=1}^k (y+t_i).
\end{equation} 
\end{Prop}
\begin{proof}
Let us denote the localization at $\id$ by 
$$\cdot |_{\id}:K_T(G/B)\to K_T(\pt).$$
Since 
$\calO_{u}|_{\id} = \delta_{u,\id}$ and 
$\lambda_y(\mathcal{V}_k)|_{\id}
=\prod_{i=1}^k (y+t_i)$, localizing \eqref{eq:clscuw} at $\id$ gives \eqref{eq:clsini}. 
\end{proof}

\begin{Th}\label{thm:classical-uniqueness}
The recursion \eqref{PURW-1} and the initial condition
\eqref{eq:clsini} determine the coefficients
$$\{c_{u,w}:u,w\in S_n\}.$$
\end{Th}

\begin{proof}
The column indexed by $\id$ is prescribed by \eqref{eq:clsini}.
For any nonidentity $w$, choose $i$ with $s_iw<w$ and apply
\eqref{PURW-1} to the already known column indexed by $s_iw$.
It expresses every entry in the column indexed by $w$ in terms of entries of
that shorter column. Induction on $\ell(w)$ proves uniqueness.
\end{proof}

\subsection{Quantum coefficients}
In $QK_T(G/B)[y]$, we write 
\begin{equation}\label{eq:qmtcuw}
\calO_{u}\star \lambda_y(\mathcal{V}_k)
=\sum_{w\in S_n} c^q_{u,w}  \calO_{w},\qquad c^q_{u,w}\in K_T(\pt)[q_1, \cdots, q_{n-1}][y].
\end{equation}
The operator $D_i^\vee$ can be extended to quantum $K$-theory by \cite{MNS}. 
We have the following property.

\begin{Prop}\label{qmtindonw}
For $1\le i\le n-1$ and $s_iu<u$ (equivalently $s_iu<_\semi u$), we have 
\begin{align}\label{qmtPURW-1}
c^q_{s_iu, w} & 
=\begin{cases}
-\dfrac{1- t_i}{t_i-t_{i+1}}(c^q_{u,w}-c^q_{u,w}|_{t_i\leftrightarrow t_{i+1}}), 
& s_iw<w,\\[15pt]
c^q_{u,s_iw}|_{t_i\leftrightarrow t_{i+1}}
    -\dfrac{1-t_i}{t_i-t_{i+1}}c^q_{u,w}
    +\dfrac{1-t_{i+1}}{t_i-t_{i+1}}c^q_{u,w}|_{t_i\leftrightarrow t_{i+1}},
& s_iw>w.\end{cases}
\end{align}
\end{Prop}
\begin{proof}
In the proof of Proposition \ref{indonw}, the same computation applies with $c_{u,w}$ replaced by $c_{u,w}^q$ and the usual product replaced by the quantum product. 
Comparing the coefficients when $s_iu<u$, we get \eqref{qmtPURW-1} immediately. 
\end{proof}

We still need the $i=0$ version of Proposition \ref{qmtindonw}. 
This is guaranteed by the following, which is essentially a reformulation of  Kato's results in \cite{Kato1}. 

\begin{Th}\label{thm:full-flag-Demazure}
There is a $\mathbb C[q_1^{\pm1},\ldots,q_{n-1}^{\pm1}]$-linear operator
$$D_0^\vee:QK_T(G/B)\longrightarrow QK_T(G/B)$$
with
\[
 D_0^\vee\calO_w=
 \begin{cases}
\calO_{s_{1n} w}\cdot  q_{w^{-1}(n),w^{-1}(1)},&s_0 w<_\semi w,\\
 \calO_w,&s_0 w>_\semi w,
 \end{cases}
\]
For $\gamma\in QK_T(G/B)$ and $f\in K_T(\pt)$, we have
\[
 D_0^\vee(f\gamma)
 =f|_{t_1\leftrightarrow t_n}D_0^\vee(\gamma)
 -\frac{1-t_n}{t_n-t_1}\bigl(f-f|_{t_1\leftrightarrow t_n}\bigr)\gamma.
\]
For $\gamma'\in QK_G(G/B)$,
\[
 D_0^\vee(\gamma' \star \gamma)=\gamma' \star D_0^\vee(\gamma).
\]
\end{Th}

\begin{proof}
Put $\theta=\epsilon_1-\epsilon_n$ and $s_\theta=s_{1n}$.
Let $w_0$ be the longest element of $S_n$. 

\subsubsection*{Step 1}
We first translate Kato's conventions into ours.
For $x\in S_n$, his finite Schubert variety is
\[
 \mathscr B(x)=\overline{Bxw_0B/B}=X(xw_0)=\overline{Bxw_0B/B},
\]
and his line bundle is
\[
 \mathcal O_{\mathscr B}(\lambda)
 =G\times_B\mathbb C_{w_0\lambda};
\]
see \cite[\S1.1]{Kato1}.
Thus the curve class $\mu\in Q^\vee$ in our notation has degree
$-w_0\mu$ in Kato's notation. In particular, his Novikov variable
$Q_i$ corresponds to our $q_{n-i}$.

Let $\mathbf Q^{\mathrm{rat}}$ be the semi-infinite flag variety. We use the completed definition of 
its K-theory $K_{T}(\mathbf Q^{\mathrm{rat}})$
in \cite[\S1.5]{Kato1}.
Write $\mathbf Q(x)$ for its Schubert variety indexed by
$x\in\widehat W$ in Kato's convention.
By \cite[Theorem 3.11, Corollary 3.13, and Theorem 4.17]{Kato1},
there is a $K_T(\pt)$-linear isomorphism
\[
 \Psi:QK_T(G/B)_{\mathrm{loc}}\longrightarrow K_{T}(\mathbf Q^{\mathrm{rat}})
\]
such that, for $w\in S_n$ and $\mu\in Q^\vee$,
\[
 \Psi\bigl([\mathcal O_{X(w)}]q^\mu\bigr)
 =\bigl[\mathcal O_{\mathbf Q((ww_0)t_{-w_0\mu})}\bigr],
\]
where $t_\nu$ is affine translation by the coroot $\nu$.
Moreover, for $1\leq i<n$ and
$\gamma\in QK_T(G/B)_{\mathrm{loc}}$,
\[
 \Psi\bigl([\mathcal O_{\mathscr B}(-\varpi_i)]\star\gamma\bigr)
 =\Xi(-\varpi_i)\Psi(\gamma),
\]
where $\varpi_i$ is the $i$-th fundamental weight and
$\Xi(-\varpi_i): K_{T}(\mathbf Q^{\mathrm{rat}})\to  K_{T}(\mathbf Q^{\mathrm{rat}})$ is the corresponding line-bundle twist.

\subsubsection*{Step 2}
Let $\leq_K$ denote Kato's semi-infinite order, defined using
sufficiently antidominant right translations. Its relation to our
order is
\[
 x\leq_K y\iff xw_0\leq_\semi yw_0
 \qquad(x,y\in\widehat W).
\]
In particular, Kato's order on $S_n$ is the
reverse of finite Bruhat order.

By \cite[Theorem 1.24]{Kato1}, $K_{T}(\mathbf Q^{\mathrm{rat}})$ carries Demazure operators
$\mathsf D_i$, $0\leq i<n$, with
\[
 \mathsf D_i[\mathcal O_{\mathbf Q(x)}]
 =\begin{cases}
 [\mathcal O_{\mathbf Q(s_ix)}],&s_ix>_Kx,\\
 [\mathcal O_{\mathbf Q(x)}],&s_ix<_Kx.
 \end{cases}
\]
They commute with right translations. Moreover,
\cite[Theorem 1.26]{Kato1} gives
\[
 \mathsf D_i\Xi(-\varpi_j)=\Xi(-\varpi_j)\mathsf D_i
 \qquad(0\leq i<n,\ 1\leq j<n).
\]
Define
\[
 D_0=\Psi^{-1}\circ\mathsf D_0\circ\Psi
 :QK_T(G/B)_{\mathrm{loc}}
 \longrightarrow QK_T(G/B)_{\mathrm{loc}}.
\]
It commutes with Novikov monomials and satisfies
\[
 D_0\bigl([\mathcal O_{\mathscr B}(-\varpi_j)]\star\gamma\bigr)
 =[\mathcal O_{\mathscr B}(-\varpi_j)]\star D_0(\gamma).
\]

By the order comparison and Section~\ref{subsec:affine-conventions},
\[
 s_0(ww_0)>_Kww_0
 \iff s_0w>_\semi w
 \iff s_\theta w<w.
\]
Also, the semidirect-product relation gives
\[
 s_0(ww_0)t_{-w_0\mu}
 =(s_\theta w)w_0t_{-w_0(\mu+w^{-1}\theta^\vee)}.
\]
Combining these identities with the formula for $\Psi$, we obtain
\[
 D_0[\mathcal O_{X(w)}]
 =\begin{cases}
 q^{w^{-1}\theta^\vee}[\mathcal O_{X(s_\theta w)}],
      &s_\theta w<w,\\
 [\mathcal O_{X(w)}],&s_\theta w>w.
 \end{cases}
\]
For $f\in K_T(\pt)$ and $\gamma\in QK_T(G/B)_{\mathrm{loc}}$,
the nil-DAHA relation in \cite[Definition 1.1(5)]{Kato1} gives
\[
 D_0(f\gamma)
 =f|_{t_1\leftrightarrow t_n}D_0(\gamma)
  +\frac{f-f|_{t_1\leftrightarrow t_n}}{1-e^{-\theta}}\gamma.
\]
The quotient belongs to $K_T(\pt)$ by the finite geometric-series
identity applied to characters. These formulas act term by term on
the completed Schubert expansions. Since there are finitely many
Novikov shifts, $D_0$ preserves finite Laurent sums.

\subsubsection*{Step 3}
Put
\[
 D_0^\vee=w_0^L\circ D_0\circ w_0^L,
\]
where $w_0^L$ is the left Weyl action of \cite{MNS}.
This action is a quantum-ring automorphism, fixes Novikov variables,
and fixes $G$-equivariant classes. Hence
\[
 D_0^\vee(a\star\gamma)=a\star D_0^\vee(\gamma)
 \qquad\bigl(a\in QK_G(G/B)_{\mathrm{loc}}\bigr).
\]
Moreover,
\[
 w_0^L[\mathcal O_{X(w)}]=\mathcal O_{w_0w},\qquad
 w_0\theta=-\theta,\qquad
 w_0s_\theta w_0=s_\theta.
\]
So we have 
$$
 D_0^\vee\mathcal O_w
 =\begin{cases}
 q^{-w^{-1}\theta^\vee}\mathcal O_{s_\theta w},
      &s_\theta w>w,\\
 \mathcal O_w,&s_\theta w<w.
 \end{cases}$$
Here
$q^{-w^{-1}\theta^\vee}=q_{w^{-1}(n),w^{-1}(1)}$, and
$s_\theta w>w$ is equivalent to $s_0w<_\semi w$ by
Section~\ref{subsec:affine-conventions}.
Conjugating the scalar rule gives
$$ D_0^\vee(f\gamma)
 =f|_{t_1\leftrightarrow t_n}D_0^\vee(\gamma)
  +\frac{f-f|_{t_1\leftrightarrow t_n}}{1-e^\theta}\gamma.$$

Note that $\mathcal{O}(-\varpi_i)$ generates $QK_G(G/B)$ after localization and completion. 
Since $D_0$ commutes with their quantum multiplication
and is $K_G(\pt)$-linear, we obtain
$$D_0^\vee(\gamma' \star \gamma)=\gamma'\star D_0^\vee(\gamma),$$
for $\gamma'\in QK_G(G/B)\subset QK_T(G/B)\ni \gamma$. 
\end{proof}

\begin{Prop}\label{qmtindonw2}
Assume that $s_{1n}u> u$  (equivalently, $s_0u<_\semi u$). Then we have 
\begin{align}\label{qmtPURW-2}
c^q_{s_{1n}u, w} & 
=q_{u^{-1}(1),u^{-1}(n)}\cdot 
\begin{cases}
-\dfrac{1- t_n}{t_n-t_{1}}(c^q_{u,w}-c^q_{u,w}|_{t_n\leftrightarrow t_{1}}), 
& s_{1n}w>w,\\[15pt]
c^q_{u, s_{1n}w}|_{t_n\leftrightarrow t_{1}}\cdot 
q_{w^{-1}(1),w^{-1}(n)}
    -\dfrac{1-t_n}{t_n-t_{1}}c^q_{u,w}
    +\dfrac{1-t_{1}}{t_n-t_{1}}c^q_{u,w}|_{t_n\leftrightarrow t_{1}},
& s_{1n}w<w. 
\end{cases}
\end{align}
\end{Prop}

\begin{proof}
    The statement is  a direct consequence of Theorem \ref{thm:full-flag-Demazure}, by   applying
$D_0^\vee$ to \eqref{eq:qmtcuw} and comparing the structure  coefficients.
\end{proof}

\begin{Prop}\label{prop:quaninitial} We have 
\begin{equation}\label{eq:quanini}
c_{\id,w}=c_{\id,w}^q,\qquad 
\text{for all $w\in S_n$}.
\end{equation}
\end{Prop}
\begin{proof}
This follows because $\mathcal{O}_{\id}=1$ is the multiplicative identity in $QK_T(G/B)$. 
\end{proof}

\begin{Th}\label{thm:quantum-uniqueness}
The recursions \eqref{qmtPURW-1} and \eqref{qmtPURW-2}, together with  the initial condition
\eqref{eq:quanini}, determine the coefficients
$$\{c_{u,w}^q:u,w\in S_n\}.$$
\end{Th}
\begin{proof}
This follows immediately from \cite[Lemma 6.12]{LNSSS}, recalled in Lemma \ref{lem:chain}. 
\end{proof}

\section{Proof of Theorem \ref{thm:main-formula}}
\label{sec:two-columns}
For $k\in [n-1]$ and  $u, w\in S_n$, recall that the  storm diagram $\mathcal{D}_k(u, w)$  has total weight 
$ F_{u, w}^q(\beta):=q_{IJ}\prod_{j=1}^nc_{u, w}^{(j)}(t_j, y)
$. Each local weight  $c_{u, w}^{(j)}(t_j, y)$ is nonzero if and only if the local configuration at interval $j$ is    given by one of the seven local configurations in Figure \ref{fig:local-configurations}, or equivalently, the cases in Table \ref{tab:local-rules} where we denote $I^{(j)}_{<b}=\{p\in I^{(j)}\,:\, p<b\}$.
 Note $$I^{(0)}:= I^{(n)}=I =\{p\leq k\,:\, w(p)<u(p)\},\qquad  J^{(0)}:=J^{(n)}=J =\{p> k\,:\, u(p)<w(p)\}.$$

  \begin{table}[!ht]
\centering
\caption{Non-vanishing local weights with $a_j=u^{-1}(j)$ and $b_j=w^{-1}(j)$}
\label{tab:local-rules}

\renewcommand{\arraystretch}{1.5}
\setlength{\tabcolsep}{5pt}
\small

\begin{tabular}{|c|l@{\(,\quad\)}l|c|c|c|}
\hline
Case & \multicolumn{2}{c|}{Condition} & $c_{u, w}^{(j)}(t_j, y)$
 & \(|I^{(j-1)}_{<b_j}|\)
 & \(|J^{(j-1)}_{<b_j}|\)\\ \hline

1 & \(a_j=b_j\le k\)
  & \(I^{(j-1)}_{<a_j}=\varnothing\)
  & \(y+t_j\)
  & \(0\) & \(0\)\\ \hline

2 & \(a_j<b_j\le k\)
  & \(I^{(j-1)}_{<a_j}=\varnothing\)
  & \(\beta^{m_j}(1+\beta t_j)\)
  & \(m_j\) & \(0\)\\ \hline

3 & \(b_j<a_j\le k\)
  & \(I^{(j-1)}_{<a_j}=\{b_j\}\)
  & \(1-\beta y\)
  & \(0\) & \(0\)\\ \hline

4 & \(k<a_j=b_j\)
  & \(J^{(j-1)}_{>a_j}=\varnothing\)
  & \(1\)
  & \(|I^{(j-1)}|\) & \(|J^{(j-1)}|\)\\ \hline

5 & \(k<b_j<a_j\)
  & \(J^{(j-1)}_{>a_j}=\varnothing\)
  & \(\beta^{m_j+1}\)
  & \(|I^{(j-1)}|\) & \(|J^{(j-1)}|-m_j-1\)\\ \hline

6 & \(b_j\le k<a_j\)
  & \(J^{(j-1)}_{>a_j}=\varnothing\)
  & \(\beta^{m_j}\)
  & \(m_j\) & \(0\)\\ \hline

7 & \(a_j\le k<b_j\)
  & \(I^{(j-1)}_{<a_j}=\varnothing\)
  & \(\beta^{m_j}(1+\beta t_j)\)
  & \(|I^{(j-1)}|\) & \(|J^{(j-1)}|-m_j\)\\ \hline
\end{tabular}
\end{table}

\noindent Recall $\deg q_i=2$ and assign  $\deg \beta=-1, \deg y=\deg t_j=1$.   Then it   follows immediately from Table \ref{tab:local-rules} that 
  each nonzero local weight $c_{u, w}^{(j)}(t_j, y)$ is a homogeneous polynomial in $t_j, y, \beta$  with 
\begin{equation}\label{eqn:localWT}
    \deg c_{u, w}^{(j)}(t_j, y)+|I^{(j-1)}_{<b_j}|-|J^{(j-1)}_{<b_j}|
    =\begin{cases}
1, &\text{if}\ a_j=b_j\leq k,\\
0,& otherwise.
\end{cases}
\end{equation}
In this section, we prove Theorem \ref{thm:main-formula} by showing that the total   weights $F_{u, w}^q(-1)$ satisfy the same recursive relation and initial condition as the structure coefficients $c_{u, w}^q$.  

\subsection{A formula for  $F_{u, w}^q(\beta)$}
Consider the following sets 
\begin{align*}
\mathsf F&=\{j\in u([k]):u^{-1}(j)=w^{-1}(j)\},\\
\mathsf U&=\{j\in u([k]):u^{-1}(j)<w^{-1}(j)\}.
\end{align*}
 \begin{Prop} 
    \label{prop:cformula}
 If $F_{u, w}^q(\beta)\neq 0$, then it is homogeneous  of degree $\ell(u)+k-\ell(w)$ and is given by  
\begin{equation}\label{lr:weight}
F_{u, w}^q(\beta)
 =q_{IJ}\beta^{-\ell(u)-k+\ell(w)+\deg q_{IJ}+|\mathsf F|}(1-\beta y)^{k-|\mathsf F|-|\mathsf U|}
   \prod_{j\in\mathsf F}(y+t_j)\prod_{j\in\mathsf U}(1+\beta t_j).
\end{equation}
 \end{Prop}

 \begin{proof}
    By definition, $F_{u, w}^q(\beta)\neq 0$   if and only if the storm  diagram $\mathcal{D}_k(u, w)$ is balanced and every interval $j$ satisfies one of the 7 conditions in Table \ref{tab:local-rules}. Since every $c_{u, w}^{(j)}(t_j, y)$ is a homogeneous polynomial, it follows directly from Table \ref{tab:local-rules} and the definitions of $\mathsf F$ and $\mathsf U$ that $$F_{u, w}^q(\beta)=q_{IJ}\beta^{M}(1-\beta y)^{k-|\mathsf F|-|\mathsf U|}
   \prod_{j\in\mathsf F}(y+t_j)\prod_{j\in\mathsf U}(1+\beta t_j).$$
 It suffices to prove that $M$ equals the desired exponent of $\beta$ in \eqref{lr:weight}.

  Since $q_{r-1}$ is contained in  $q_{pn}$  exactly when $p<r$, its exponent in $q_{IJ}
=\frac{\prod_{p\in I}q_{p,n}}{\prod_{p\in J}q_{p,n}}$ equals $|I_{<r}|-|J_{<r}|$. We first compare this difference with that for interval $j$. Note that  path $p$ increases the counting  by $1$ if $u(p)\le j<w(p)$, or by $-1$ if $w(p)\le j<u(p)$, and by $0$ otherwise. We have 
$$  
|I^{(j)}_{<r}|-|J^{(j)}_{<r}|
=|I_{<r}|-|J_{<r}|
+\#\{p<r:u(p)\le j\}
-\#\{p<r:w(p)\le j\}.
 $$
Also note $|I^{(j)}_{<a_j}|=|J^{(j)}_{<a_j}|$ for all $j$. (Indeed, if $a_j>k$, then $I^{(j)}_{<a_j}=I^{(j)}$ by definition, and $J^{(j)}_{\geq a_j}=\varnothing$ by Table \ref{tab:local-rules} so that $J^{(j)}_{< a_j}=J^{(j)}$. If $a_j\leq k$, then both sets are empty.)
 Applying the above equalities, we obtain
\begin{align*}
&\quad \sum_{j=1}^n
 \left(|I^{(j-1)}_{<b_j}|-|J^{(j-1)}_{<b_j}|\right)\\
&=\sum_{r=1}^n
 \Bigl(
 |I^{(w(r)-1)}_{<r}|-|J^{(w(r)-1)}_{<r}|
 +|I^{(u(r))}_{<r}|-|J^{(u(r))}_{<r}|
 \Bigr)\\
&=\deg q_{IJ}+\ell(w)-\ell(u)
 + \sum_{r=1}^n\big(|\{p<r\,:\,
  w(p)>u(r)\}|-|\{p<r\,:\, u(p)\ge w(r)\}|
 \big) \\
 &=\deg q_{IJ}+\ell(w)-\ell(u)
 +  |\mathsf F|-k.
\end{align*} 
Here the second equality holds since 
\begin{align*}
\ell(w)-\ell(u)&=\sum_{r=1}^n(|\{p<r\,:\, w(p)>w(r)\}|-|\{p<r\,:\, u(p)>u(r)\}|)\\
&=\sum_{r=1}^n(-|\{p<r\,:\, w(p)<w(r)\}|+|\{p<r\,:\, u(p)<u(r)\}|).
\end{align*}
The sum in the third line is equal to
$$
 |\{(p,r)\in[n]^2:p\ne r,\ u(p)<w(r)\}| -\binom n2 
 =-|\{r:u(r)<w(r)\}|.$$
Using $|I|=|J|$, this becomes $-k+|\mathsf F|$, which proves the third equality. Equation \eqref{eqn:localWT} now yields the claim. 
 \end{proof}

\subsection{A relation among four diagrams}

Let us consider two adjacent intervals. 
For simplicity, let us consider the first and the second interval. 
Recall their equivariant parameters are $t_{1}$ and $t_{2}$ respectively. 
Fix all paths outside these two intervals. 
In particular, the two boundary walls are fixed. 
Let $a<b$ be the two source labels and $c<d$ the two target labels.
We use the following four arrangements
$$
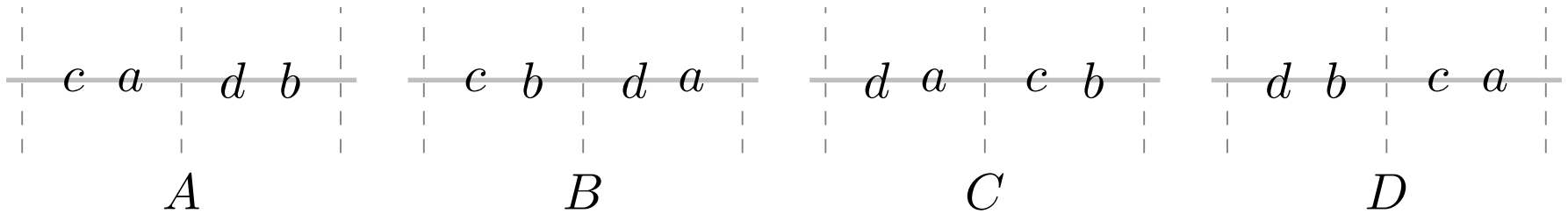
\def\twocolumnframe#1#2#3#4#5{
\begin{eye}{1}{2}{2}
  \wind{#2}{#3}{0}{1}{1}
  \wind{#4}{#5}{0}{2}{2}
  \label{1.5}{$#1$}
\end{eye}}
\begin{matrix}
 \twocolumnframe{A}{c}{a}{d}{b}&
 \twocolumnframe{B}{c}{b}{d}{a}&
 \twocolumnframe{C}{d}{a}{c}{b}&
 \twocolumnframe{D}{d}{b}{c}{a}
\end{matrix}
$$
where the left
and right boundary walls of the four arrangements are the same. 
Precisely, we take $u, w\in S_n$ such that $u(a)=1=w(c)$ and $u(b)=2=w(d)$. Then the  products
of the two interval weights are respectively given by 
\begin{align*}
&A=c_{u, w}^{(1)}(t_1,y)c_{u, w}^{(2)}(t_2,y), 
&&B=c_{s_1u, w}^{(1)}(t_1,y)c_{s_1u, w}^{(2)}(t_2,y),\\
&C=c_{u, s_1w}^{(1)}(t_1,y)c_{u, s_1w}^{(2)}(t_2,y), 
&&D=c_{s_1u, s_1w}^{(1)}(t_1,y)c_{s_1u, s_1w}^{(2)}(t_2,y).
\end{align*}

The weights \(A,B,C,D\) are understood as contributions with the prescribed data on the two outer boundary walls. An arrangement contributes only if it can be realized by the seven allowed local configurations while matching these fixed boundary data; otherwise, its contribution is defined to be zero. In particular, an arrangement may admit a realization with nonzero local weights for some boundary data, yet have zero contribution when the outer walls are fixed differently.

For example, consider Case \ref{app:case-4-9} in Appendix \ref{sec:twocolproof} with \(a=c=1\), \(k=1\), \(b=d=2\), and \(m_1=\ell_1=0\). In the first row, both outer walls are empty, and \(B=0\). In contrast, the \(B\)-arrangement in the second row is realized by Cases 6 and 7 of Table \ref{tab:local-rules}, whose weights are \(1\) and \(1+\beta t_2\), respectively. This realization requires both outer walls to carry the red path labeled \(1\) and the blue path labeled \(2\). Hence the nonzero product \(1+\beta t_2\) is admissible for the boundary data in the second row, but contributes zero under the empty outer-wall condition in the first row.

\begin{Prop}\label{prop:two-column}
The four weights satisfy
\begin{align}
 (t_1-t_{2})C&=(1+\beta t_{1})(A-A|_{t_1\leftrightarrow t_{2}}),
 \label{eq:local-ascending-column}\\
 (t_1-t_2)D&=(1+\beta t_1)B-(1+\beta t_2)B|_{t_1\leftrightarrow t_2}+(t_1-t_2)A|_{t_1\leftrightarrow t_2}.
 \label{eq:local-descending-column}
\end{align}
Equivalently,
\begin{align}
 (t_1-t_2)A&=(t_1-t_2)D|_{t_1\leftrightarrow t_2}-(1+\beta t_1)B+(1+\beta t_2)B|_{t_1\leftrightarrow t_2},
 \label{eq:local-ascending-target-row}\\
 (t_1-t_2)C&=(1+\beta t_1)(D|_{t_1\leftrightarrow t_2}-D).
 \label{eq:local-descending-target-row}
\end{align}
\end{Prop}

\begin{proof}
The proof is by case-by-case verification. 
All the cases are displayed in Appendix \ref{sec:twocolproof}.
We illustrate the verification with Case~\ref{app:case-4-4}, where $c<a\leq k<d<b$. 
The four diagrams below are arranged in the order $A,B,C,D$.
The multiplicities $m_1,\ell_1,\ell_2$ count paths crossing both intervals
with labels in $(a,k]$, $(k,d)$, and $(d,b)$, respectively.
\eyesetup{interval width=12mm,level distance=2.5mm,first level=3.5mm,
  endpoint gap=1.5mm,label offset=.18,label font=\scriptsize,
  wall label sep=1.5pt,wall padding=1.5mm,path width=.55pt,
  bundle width=1.5pt,arrow length=1.4mm,arrow width=1mm,
  bundle arrow length=2mm,bundle arrow width=1.5mm}
\setlength{\arraycolsep}{2pt}
\[
\begin{matrix}
\begin{eye}{3}{7}{2}
  \wind{a}{}{2}{1}{1.5}
  \wind{}{c}{1}{0.5}{1}
  \wind{}{}{7}{1.5}{0.5}
  \wind{b}{}{7}{2}{1.5}
  \wind{}{d}{5}{2.5}{2}
  \wind{}{}{2}{1.5}{2.5}
  \wind[thick]{}{m_1}{3}{0.5}{2.5}
  \wind[thick]{\ell_1}{}{4}{2.5}{0.5}
  \wind[thick]{\ell_2}{}{6}{2.5}{0.5}
  \label{1.5}{\tiny $\beta^{1+\ell_2}(1-\beta y)$}
\end{eye} &
\begin{eye}{3}{7}{2}
  \wind{b}{}{7}{1}{0.5}
  \wind{}{c}{1}{0.5}{1}
  \wind{a}{}{2}{2}{2.5}
  \wind{}{d}{5}{2.5}{2}
  \wind[thick]{}{m_1}{3}{0.5}{2.5}
  \wind[thick]{\ell_1}{}{4}{2.5}{0.5}
  \wind[thick]{\ell_2}{}{6}{2.5}{0.5}
  \label{1.5}{\tiny $\beta^{\ell_2}(1+\beta t_2)$}
\end{eye} &
\begin{eye}{3}{7}{2}
  \wind{a}{}{2}{1}{1.5}
  \wind{}{d}{5}{1.5}{1}
  \wind{}{}{1}{0.5}{1.5}
  \wind{}{}{7}{1.5}{0.5}
  \wind{b}{}{7}{2}{1.5}
  \wind{}{c}{1}{1.5}{2}
  \wind{}{}{2}{1.5}{2.5}
  \wind{}{}{5}{2.5}{1.5}
  \wind[thick]{}{m_1}{3}{0.5}{2.5}
  \wind[thick]{\ell_1}{}{4}{2.5}{0.5}
  \wind[thick]{\ell_2}{}{6}{2.5}{0.5}
  \label{1}{$0$}
\end{eye} &
\begin{eye}{3}{7}{2}
  \wind{b}{}{7}{1}{0.5}
  \wind{}{d}{5}{1.5}{1}
  \wind{}{}{1}{0.5}{1.5}
  \wind{a}{}{2}{2}{2.5}
  \wind{}{c}{1}{1.5}{2}
  \wind{}{}{5}{2.5}{1.5}
  \wind[thick]{}{m_1}{3}{0.5}{2.5}
  \wind[thick]{\ell_1}{}{4}{2.5}{0.5}
  \wind[thick]{\ell_2}{}{6}{2.5}{0.5}
  \label{1.5}{\tiny $\beta^{1+\ell_2}(1-\beta y)$}
\end{eye}
\end{matrix}
\]
A direct calculation verifies  \eqref{eq:local-ascending-column} and
\eqref{eq:local-descending-column}.
\end{proof}

\subsection{The classical case} 
Let $F_{u,w}=F^{q}_{u, w}(-1)|_{q_1=\cdots =q_{n-1}=0}$, which is the classical version of $F^q_{u, w}(-1)$. 
Then $F_{u,w}=0$ whenever there is a path
crossing the cut between intervals $n$ and $1$. 
With Proposition \ref{prop:two-column}, we can prove the classical version of Theorem \ref{thm:main-formula}. 

 \begin{Th}\label{thm:classical-model}
The weights $F_{u,w}$ satisfy the recursion and the initial condition in Theorem \ref{thm:classical-uniqueness}. In particular, 
$F_{u,w}=c_{u,w}$ for all $u,w\in S_n$.
\end{Th}

\begin{proof}
Fix $1\leq i<n$. In the four arrangements of
Proposition~\ref{prop:two-column}, choose $u,w$ with
\[
 u^{-1}(i)=a<b=u^{-1}(i+1),\qquad
 w^{-1}(i)=c<d=w^{-1}(i+1).
\]
In Proposition~\ref{prop:two-column}, replace the local parameter
pair $(t_1,t_2)$ by $(t_i,t_{i+1})$. 
After setting the Novikov variables to zero, all contributing diagrams have the same exterior configuration: interchanging the adjacent values \(i\) and \(i+1\) changes no exterior wall data for diagrams that do not cross the cut. Apply Proposition \ref{prop:two-column} with this exterior configuration fixed and let \(Z\) be the product of the weights of the remaining intervals.
It is easy to check that
$$F_{u,w}=Z\cdot A,\qquad 
F_{s_i u,w} = Z\cdot B,\qquad 
F_{u,s_iw}= Z\cdot C,\qquad 
F_{s_i u,s_iw} = Z\cdot D. $$
Since $s_iw>w$, equations
\eqref{eq:local-ascending-column} and
\eqref{eq:local-descending-column} give, respectively,
\begin{align*}
 F_{u,s_iw}
 &=\frac{1-t_i}{t_i-t_{i+1}}
       (F_{u,w}-F_{u,w}|_{t_i\leftrightarrow t_{i+1}}),\\
 F_{s_i u,s_iw}
 &=\frac{(1-t_i)F_{s_i u,w}-(1-t_{i+1})F_{s_i u,w}|_{t_i\leftrightarrow t_{i+1}}}
         {t_i-t_{i+1}}
   +F_{u,w}|_{t_i\leftrightarrow t_{i+1}}.
\end{align*}
The first formula is the case $s_i u>u$ of \eqref{PURW-1}.
Replacing $s_i u$ by $u$ in the second formula gives its case
$s_i u<u$. Thus both branches of the recursion hold.
This proves the recursion \eqref{PURW-1} 
 in Proposition \ref{thm:clsrec}. 

Now let us check the initial condition. 
Assume $w=\id$, and $F_{u,w}\neq 0$. 
The following diagram shows the case when $n=6$, where $\bar{u}_i=u^{-1}(i)$. 
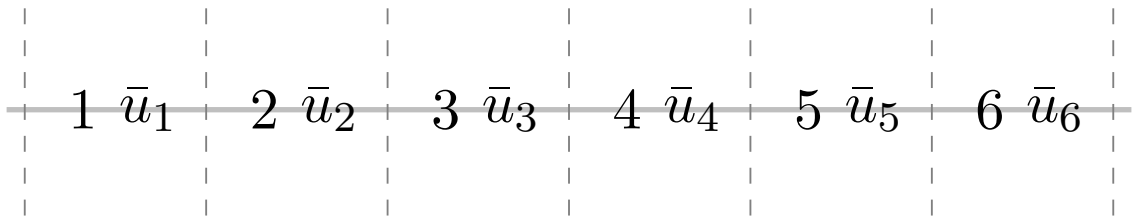
$$\begin{eye}{2}{4}{6}
\wind{1}{\bar{u}_1}{0}{1}{1}
\wind{2}{\bar{u}_2}{0}{2}{2}
\wind{3}{\bar{u}_3}{0}{3}{3}
\wind{4}{\bar{u}_4}{0}{4}{4}
\wind{5}{\bar{u}_5}{0}{5}{5}
\wind{6}{\bar{u}_6}{0}{6}{6}
\end{eye}$$
We claim that the $i$-th interval has to be 
$$
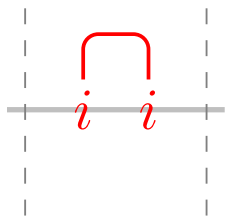
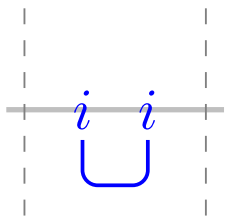
\begin{eye}{2}{4}{1}
\wind{i}{i}{2}{1}{1}
\label{1}{\(y+t_i\)}
\end{eye}
\quad \text{ if $1\leq i\leq k$ }\qquad 
\begin{eye}{2}{4}{1}
\wind{i}{i}{3}{1}{1}
\label{1}{\(1\)}
\end{eye}
\quad \text{ if $k< i\leq n$}.
$$
In the first interval, since we set $q_1=\cdots=q_{n-1}=0$, the left boundary is empty. 
By checking the nonzero weights, it is easy to see the claim for $i=1$. 
Similarly, the claim is true for the second interval, the third interval, etc. 
Thus $F_{u,\id}=0$ unless $u=\id$, and 
$F_{\id,\id}=\prod_{i=1}^k(y+t_i)$. This proves the initial condition 
\eqref{eq:clsini} in Theorem \ref{thm:classical-uniqueness}. 
\end{proof}

\subsection{The quantum case}
With Proposition \ref{prop:two-column} and Theorem \ref{thm:classical-model}, we can prove Theorem \ref{thm:main-formula}. 

\begin{Th}
The weights $F^q_{u,w}$ satisfy the recursion and the initial condition in Theorem \ref{thm:quantum-uniqueness}. In particular, 
$F^q_{u,w}=c^q_{u,w}$ for all $u,w\in S_n$.
\end{Th}
\begin{proof}
For $1\leq i< n$, the argument is similar to the proof of  Theorem \ref{thm:classical-model}. 
Choose $u,w$ with
\[u^{-1}(i)=a<b=u^{-1}(i+1),\qquad
 w^{-1}(i)=c<d=w^{-1}(i+1).\]
Then we can write
$$F^q_{u,w}=Z^q\cdot A,\qquad 
F^q_{s_i u,w} = Z^q\cdot B,\qquad 
F^q_{u,s_iw}= Z^q\cdot C,\qquad 
F^q_{s_i u,s_iw} = Z^q\cdot D$$
where $Z^q$ is the product of weights for the remaining $n-2$ intervals and the Novikov variable. 
Equations \eqref{eq:local-ascending-target-row} and
\eqref{eq:local-descending-target-row} give
\begin{align*}
 F^q_{u,w}
 &=F^q_{s_i u,s_iw}|_{t_i\leftrightarrow t_{i+1}}
   -\frac{(1-t_i)F^q_{s_i u,w}-(1-t_{i+1})F^q_{s_i u,w}|_{t_i\leftrightarrow t_{i+1}}}
          {t_i-t_{i+1}},\\
 F^q_{u,s_iw}
 &=\frac{1-t_i}{t_i-t_{i+1}}
      \left(F^q_{s_i u,s_iw}|_{t_i\leftrightarrow t_{i+1}}-F^q_{s_i u,s_iw}\right).
\end{align*}
Here $s_i u>u$ and $s_iw>w$.
Replacing $s_i u$ by the source index in
\eqref{qmtPURW-1}, these are its two branches.
This proves the first recursion \eqref{qmtPURW-1} in Theorem \ref{thm:quantum-uniqueness}. 

For the second recursion, the proof is very similar, with the convention $i=n$ and $i+1\equiv 1\bmod n$. The only thing we need to worry about is the Novikov parameters. 
We use the same setup as above and keep track of the additional Novikov factors. 
Recall that $q_{a,b}=q^{\epsilon_a-\epsilon_b}$ for any $1\leq a,b\leq n$ (not necessarily $a<b$). 
Assume that 
$F_{u,w}^q=Z^q\cdot A$, then by Lemma \ref{lem:shiftQpara} below, we have
$$F_{s_i u,w}^q = q_{a,b}^{-1}\cdot Z^q\cdot B,\qquad 
F_{u,s_iw}^q= q_{c,d}\cdot Z^q\cdot C,\qquad 
F_{s_i u,s_iw}^q = q_{a,b}^{-1}q_{c,d}\cdot Z^q\cdot D.$$
Combining these identities we get the second recursion \eqref{qmtPURW-2}. 

Now let us check the initial condition. 
Assume $u=\id$, and $F^q_{u,w}\neq 0$. 
The following diagram shows the case when $n=6$, where $\bar{w}_i=w^{-1}(i)$. 
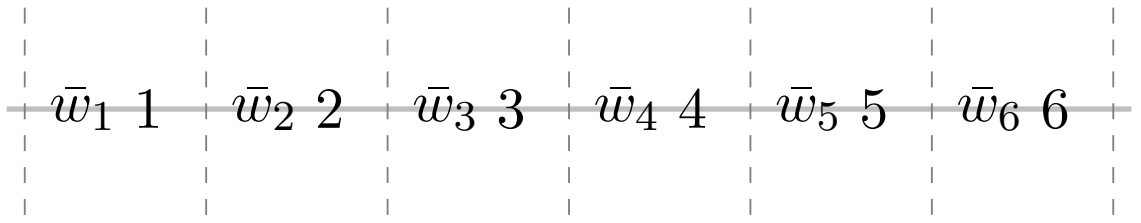
$$\begin{eye}{2}{4}{6}
\wind{\bar{w}_1}{1}{0}{1}{1}
\wind{\bar{w}_2}{2}{0}{2}{2}
\wind{\bar{w}_3}{3}{0}{3}{3}
\wind{\bar{w}_4}{4}{0}{4}{4}
\wind{\bar{w}_5}{5}{0}{5}{5}
\wind{\bar{w}_6}{6}{0}{6}{6}
\end{eye}$$
Since neither of the following local configurations is allowed
$$
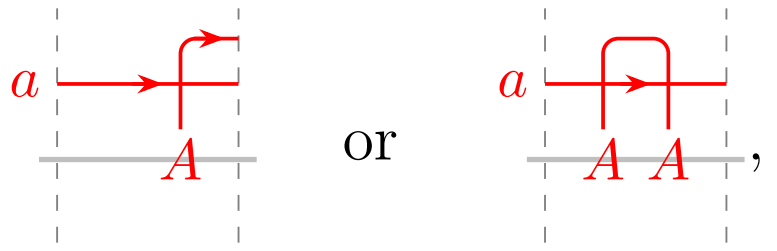
\begin{eye}{4}{5}{1}
\wind{A}{}{4}{1}{1.5}
\wind{a}{}{2}{0.5}{1.5}
\end{eye}
\quad\text{ or }\quad 
\begin{eye}{4}{5}{1}
\wind{A}{A}{4}{1}{1}
\wind{a}{}{2}{0.5}{1.5}
\end{eye},$$
the path labeled by $a\leq k-1$ cannot go beyond $a+1\leq k$. 
We can conclude that $1,\ldots,k-1\in \{\bar{w}_1,\bar{w}_2,\ldots,\bar{w}_k\}$. 
In particular, $1,\ldots,k-1$ will not contribute to the Novikov parameters.
If the path with label $k$
contributes a nontrivial Novikov factor, it must pass through interval $k+1$. 
This  contradicts  Lemma \ref{lem:onewallnopath} below. 
This proves $F_{\id,w}=F_{\id,w}^q$. 
\end{proof}

\begin{Lemma}\label{lem:shiftQpara}
If the following interval has nonzero weight, 
$$\begin{eye}[interval width=15mm]{2}{4}{1}
\node at ({0.5*\eyeintervalwidth},0.5) [left] {\(A\)}; 
\node at ({0.5*\eyeintervalwidth},-0.5) [left] {\(B\)}; 
\node at ({1.5*\eyeintervalwidth},0.5) [right]{\(A'\)}; 
\node at ({1.5*\eyeintervalwidth},-0.5) [right] {\(B'\)}; 
\wind{a}{b}{0}{1}{1}
\end{eye}$$
then 
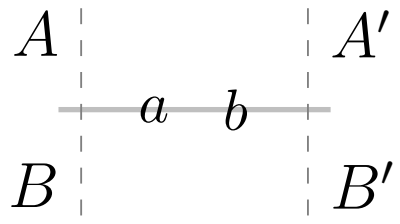
$$q_{A',B'} = q_{b,a}
\cdot q_{A,B}.$$
\end{Lemma}
\begin{proof}
This follows directly by examining the seven nonzero weights. 
\end{proof}

\begin{Lemma}\label{lem:onewallnopath}
Let $u,w\in S_n$ with $F_{u,w}^q\neq 0$. 
Then the right wall of the $u(k+1)$-th interval is empty. 
\end{Lemma}
\begin{proof}
The path with label $k+1$
is the lowest blue path. If $u(k+1)=w(k+1)$,   then the nonzero local configuration is given by the fourth picture in Figure \ref{fig:local-configurations}  with the right wall empty. If $u(k+1)\neq w(k+1)$, then  the only nonzero possibility is given by the sixth picture in Figure \ref{fig:local-configurations}: 
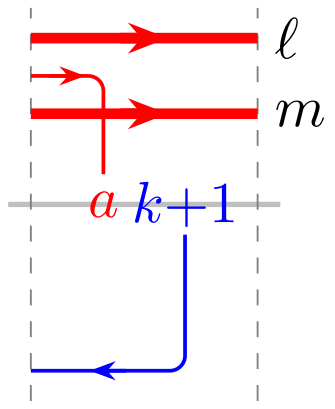
$$\begin{eye}[interval width=15mm,level distance=2.5mm]{4}{8}{1}
  \wind{}{a}{3}{.5}{1}
  \wind[thick]{}{m}{2}{.5}{1.5}
  \wind[thick]{}{\ell}{4}{.5}{1.5}
  \wind{k\text{+}1}{}{5}{1}{.5}
\end{eye}\qquad $$
By the balanced condition, $m=\ell=0$. 
\end{proof}

\section{Applications}
In this section, we establish Pieri rules for $QK_T(G/P)$ of partial flag varieties, as applications of  Theorem \ref{thm:main-formula}.
\subsection{Equivariant quantum Pieri rule for $QK_T(F\ell_n)$} 
For any $u, w\in S_n$,
 we consider the normalized localization of $\mathcal{O}_u$ at the $T$-fixed point indexed by $w$: 
\begin{equation}\label{eqn:localization}
 \psi_u(w):=(-1)^{\ell(u)}N_{u,w}^{w}.
\end{equation}
Moreover,  for $A\subsetneq[n]$, we consider the Grassmannian permutation  $w_A\in S_n$ with  the one-line notation $[w_A(1), \cdots, w_A(n)]$   listing $A=\{a_1, a_2, \cdots, a_r\}$ increasingly, followed by its
complement increasingly. Then for $0\leq b\leq r$, we simply denote 
\begin{equation}\label{lr:localization}
 \psi_b(A):=\psi_{s_{r-b+1}s_{r-b+2}\cdots s_r}(w_A),\quad \mbox{where}\quad \psi_0(A)=\psi_{\rm id}(w_A)=1.
\end{equation}

\begin{Rmk}\label{rmk:loccoef}
Several  manifestly positive formulae for these localization   coefficients are known \cite{KoKu, McN, GrKr}. Here we provide one derived from \cite{McN}.  Assume $b\geq 1$. Write
\[
 [r-b+1]\setminus A=\{j_1<\cdots<j_m\},\qquad
 A\setminus[r-b+1]=\{p_{1}<\cdots<p_{{m+b-1}}\}.
\]
Then we have 
\[
 \psi_b(A)
 =
 \sum_{S\in{[m+b-1]\choose b}}
 \prod_{\ell\in S} (e^{\epsilon_{p_{\ell}}-\epsilon_{j_S(\ell)}}-1) 
 \prod_{\substack{1\le\ell<\max S\\\ell\notin S}}
e^{\epsilon_{p_{\ell}}-\epsilon_{j_S(\ell)}}\in \mathbb{Z}_{\ge 0}[e^{-\alpha_1}-1, \cdots, e^{-\alpha_{n-1}}-1],
\]
where $j_S(\ell):=j_{1+|[\ell-1]\setminus S|}<p_\ell$, and the sum is defined to be $0$ when $m=0$.
\end{Rmk}

Denote $\check{\mathbb{P}}^{i-1}:= \mathbb P\left(\left(\bigoplus_{j=1}^i\mathbb{C}_{\epsilon_j}\right)^\vee\right)$ and $L:=\calO_{\check{\mathbb{P}}^{i-1}}(-1)$. 
\begin{Lemma}\label{lemma:Olambda}
Let $1\leq i\leq k<n$ and  $y=L^{-1}-1$. In $K_T(G/B)$, we have   
\begin{equation}\label{eqn:Os} 
 \calO_{s_is_{i+1}\cdots s_k}
 = \chi_T\left(\check{\mathbb P}^{i-1},L^k \lambda_y(\mathcal V_k)\right),   
\end{equation}
 where we view $L^k \lambda_y(\mathcal V_k)\in K_T(\check{\mathbb P}^{i-1})\otimes_{R(T)}K_T(G/B)$ and apply $\chi_T\left(\check{\mathbb P}^{i-1}, -\right)$ to the first factor.  
 
 Moreover, for any     $A\subsetneq [n]$,   we have 
 \begin{equation*}
     \chi_T\left(\check{\mathbb P}^{i-1}, \prod_{a\in A}(1-(1-t_a)L)\right)=\begin{cases}
(-1)^{|A|-i+1}\psi_{|A|-i+1}(A), &\text{if}\ |A|\geq i-1,\\[5pt]
1,& otherwise.
\end{cases}
 \end{equation*}
\end{Lemma}

\begin{proof}
    Multiplying  both sides of  \eqref{eqn:lambday} by $L^k$, we obtain 
\begin{equation}\label{eqn:Lklamb}
    L^k\lambda_{y}(\mathcal V_k)
=\prod_{h=1}^{k}\bigl(1-(1-t_h)L\bigr)
+\sum_{j=1}^{k}(1-t_j)L
\prod_{h<j}\bigl(1-(1-t_h)L\bigr)
\mathcal O_{s_j\cdots s_k}.
\end{equation}

The evaluation map $\mathbb C^i\otimes L\to\mathcal O$  is surjective.
Thus its Koszul complex is exact, and consequently we have  
$$
0=\sum\nolimits_{m=0}^{i}(-1)^m
[\textstyle\bigwedge^m\mathbb C^i]L^m
=\prod_{h=1}^{i}\bigl(1-(1-t_h)L\bigr).
$$
Then the first product in \eqref{eqn:Lklamb} is zero since $k\ge i$.
Note that if $m=0$, then
$$
\chi_T\left(\check{\mathbb P}^{i-1}, L^m\right)
=\sum_{j=0}^{i-1}(-1)^j
\bigl[H^j(\check{\mathbb P}^{i-1},L^{\otimes m})\bigr]
=1$$
If  $1\leq m<i$, then $\chi_T\left(\check{\mathbb P}^{i-1}, L^m\right)
=0.$
Thus it follows from 
 \eqref{eqn:Lklamb} that
\begin{align*}
 \chi_T\left(\check{\mathbb P}^{i-1},L^k \lambda_y(\mathcal V_k)\right)&=\chi_T\left(\check{\mathbb P}^{i-1}, (1-t_i)L\prod_{h<i}\bigl(1-(1-t_h)L\bigr) \mathcal O_{s_i\cdots s_k}\right)\\
&=\chi_T\left(\check{\mathbb P}^{i-1}, \prod_{h<i}\bigl(1-(1-t_h)L\bigr)\right)\mathcal O_{s_i\cdots s_k}=\mathcal O_{s_i\cdots s_k}.
\end{align*}
Therefore  \eqref{eqn:Os} holds.

Note that $\mathcal V_r|_{w_A}
=\bigoplus_{a\in A}\mathbb C_{\epsilon_a}$. Thus 
$
\left.\lambda_y(\mathcal V_r)\right|_{w_A}
=\prod_{a\in A}(y+t_a).
$
For  $r:=|A|\geq i-1$, we have 
\begin{align*}
\left. N_{s_i\cdots s_r, w_A}^{w_A}=\mathcal O_{s_i\cdots s_r}\right|_{w_A}
&=\chi_T\left(
\check{\mathbb P}^{i-1},
L^r\prod_{a\in A}(y+t_a)
\right)\\
&=\chi_T\left(
\check{\mathbb P}^{i-1},
\prod_{a\in A}\bigl(1-(1-t_a)L\bigr)
\right).
\end{align*}
Then the second half of the statement follows from the definition of $\psi_{r-i+1}(A)$ when $r\geq i-1$, and follows from a direct expansion when $r<i-1$.
\end{proof}

We use the conventions: ${a \choose b}=0$ unless $0\leq b\leq a$;   ${0 \choose b}=1$ if $b=0$, or $0$ otherwise.
 Recall \[
\mathsf F=\{j\in u([k]):u^{-1}(j)=w^{-1}(j)\}, \qquad  \mathsf U=\{j\in u([k]):u^{-1}(j)<w^{-1}(j)\}.
\]

\begin{Th}\label{thm:equivPieri}
Let $u,w\in S_n$, $1\leq i\leq k<n$, and $d\in Q^\vee$. Then   $N_{u, s_{i}s_{i+1}\cdots s_k}^{w,d}=0$ unless  $c_{u, w}^q(t, y)\neq 0$ and $q^d=q_{IJ}$.  
When both conditions hold, we have
\begin{equation*} 
 (-1)^{\ell(u)+k-i+1-\ell(w)}
 N_{u,s_{i}s_{i+1}\cdots s_k}^{w,d}
 =\binom{|\mathsf U|-1}{i-1-|\mathsf F|}
+ \sum_{\substack{S\subseteq\mathsf U~:~|S|\ge i-|\mathsf F|}}
\psi_{|\mathsf F|+|S|-i+1}({\mathsf F\cup S}).
\end{equation*}
In particular, 
$$(-1)^{\ell(u)+k-i+1-\ell(w)}
 N_{u,s_{i}s_{i+1}\cdots s_k}^{w,d}\in \mathbb{Z}_{\ge 0}[e^{-\alpha_1}-1, \cdots, e^{-\alpha_{n-1}}-1].
 $$
\end{Th}

 \begin{proof}
By \eqref{eqn:Os} we have
\begin{align*}
\mathcal O_u\star\mathcal O_{s_i\cdots s_k}
&=\chi_T\left(
\check{\mathbb P}^{i-1},
L^k 
\bigl(\mathcal O_u\star\lambda_y(\mathcal V_k)|_{y=L^{-1}-1}\bigr)
\right)\\
&=\sum_{w\in S_n}
\chi_T\left(
\check{\mathbb P}^{i-1},
L^k c^q_{u,w}(t,L^{-1}-1)
\right)\mathcal O_w.
\end{align*}
Therefore, we have  
$$
\sum_d N_{u,s_i\cdots s_k}^{w,d}q^d
=
\chi_T\!\left(
\check{\mathbb P}^{i-1},
L^k c^q_{u,w}(t,L^{-1}-1)
\right).
$$
Hence, $\sum_d N_{u,s_i\cdots s_k}^{w,d}q^d$ does not vanish only if   $c^q_{u,w}(t,y)\neq 0$; when this holds, 
it is a scalar (possibly zero) multiple of the unique  $q^d=q_{IJ}$. Moreover, by   Proposition~\ref{prop:cformula}, we have   
\begin{equation}\label{eqn:PieriEuler}
\begin{aligned}
&(-1)^{\ell(u)+k-i+1-\ell(w)}
N_{u,s_i\cdots s_k}^{w,d}=
(-1)^{i-1+|\mathsf F|}
\chi_T\left(
\check{\mathbb P}^{i-1},
\prod_{a\in\mathsf F}\bigl(1-(1-t_a)L\bigr)
\prod_{a\in\mathsf U}\bigl((1-t_a)L\bigr)
\right).
\end{aligned}
\end{equation}

Suppose first that $\mathsf U=\varnothing$.
Then by definition,  $w^{-1}u(a)\le a$ for $1\le a\le k$, implying $w^{-1}u(a)=a$ for all $a\le k$. Thus $|\mathsf F|=k$.
By Lemma \ref{lemma:Olambda}, the right-hand side of
\eqref{eqn:PieriEuler} is therefore
$\psi_{k-i+1}(\mathsf F)$. Thus the statement  holds.

Now suppose that $\mathsf U\ne\varnothing$, and fix
$b\in\mathsf U$.   Expanding
$(1-t_a)L=1-\bigl(1-(1-t_a)L\bigr)$, \,\,$a\in\mathsf U\setminus\{b\}$, we obtain
\begin{align*}
&(-1)^{i-1+|\mathsf F|}
\chi_T\left(
\check{\mathbb P}^{i-1},
\prod_{a\in\mathsf F}\bigl(1-(1-t_a)L\bigr)
\prod_{a\in\mathsf U}\bigl((1-t_a)L\bigr)
\right)\\
&\quad= \sum_{S\subseteq\mathsf U\setminus\{b\}}
(-1)^{|\mathsf F|+|S|-i+1}
\chi_T\left(
\check{\mathbb P}^{i-1},
(1-t_b)L
\prod_{a\in\mathsf F\cup S}\bigl(1-(1-t_a)L\bigr)
\right)\\
&\quad=\sum_{S\subseteq\mathsf U\setminus\{b\}}
(-1)^{|\mathsf F|+|S|-i+1}
\chi_T\left(
\check{\mathbb P}^{i-1},
\prod_{a\in \mathsf F\cup S}\bigl(1-(1-t_a)L\bigr)
-
\prod_{a\in \mathsf F\cup S\cup\{b\}}\bigl(1-(1-t_a)L\bigr)
\right)\\
&\quad=\sum_{S\subseteq\mathsf U\setminus\{b\}: |\mathsf F\cup S|\geq i-1}\left(\psi_{|\mathsf F\cup S|-i+1}(\mathsf F\cup S)+\psi_{|\mathsf F\cup S\cup \{b\}|-i+1}(\mathsf F\cup S\cup \{b\})\right)\\
&\quad=\sum_{S\subseteq\mathsf U\setminus\{b\}: |\mathsf F\cup S|=i-1}\psi_{|\mathsf F\cup S|-i+1}(\mathsf F\cup S)+\sum_{S\subseteq\mathsf U: |\mathsf F\cup S|\geq i}\psi_{|\mathsf F\cup S|-i+1}(\mathsf F\cup S). 
\end{align*}
 Here the third equality holds by Lemma \ref{lemma:Olambda}, and the first sum on the last row contributes $\binom{|\mathsf U|-1}{i-1-|\mathsf F|}$. Thus the statement follows, by using \eqref{eqn:PieriEuler} and noting  $\mathsf F\cap\mathsf U=\varnothing$.  
\end{proof}

\begin{Rmk}
Taking the non-equivariant limit $\alpha_i=0$ in Theorem \ref{thm:equivPieri}, we obtain the Pieri rule for $QK(F\ell_n)$ given by the binomial coefficient  in that formula. 
A combinatorial argument identifies the degree-zero specialization with the classical Pieri rule in \cite{LeSo}.
A direct combinatorial proof of the equivalence between  the non-equivariant limit and the quantum Pieri rule in \cite{NaSa} is still missing. 
\end{Rmk}

\subsection{Seidel action}
      
   The $K$-theoretic (equivariant) Seidel product concerns  quantum multiplication by a specific class in (equivariant) quantum $K$-theory. 
It has been studied for     
(co)minuscule Grassmannians  in  \cite{BCP} (see also  \cite{BM11, LLSY} for complex Grassmannians) in the non-equivariant setting, and in \cite{BCP0} in  the equivariant setting.   The $K$-theoretic  equivariant Seidel product for any flag variety of general Lie type was recently proved in \cite{IKN}. Here we reprove   the equivariant Seidel product in type $A$  using Theorem \ref{thm:main-formula}, and show further applications in the next subsection.

Let $\eta=s_1s_2\cdots s_{n-1}=(12\cdots n)\in S_n$.    Recall $G=GL_n$ and $F\ell_n=G/B$.    Define a  $\mathbb{C}$-linear map $\mathcal{S}$ by 
$$
\mathcal S:QK_T(F\ell_n)\longrightarrow QK_T(F\ell_n),
\qquad
\mathcal S\!\big(q^df_u\mathcal O_u\big)
=q^d\eta(f_u)\,q_{u^{-1}(n),n}\mathcal O_{\eta u},
$$
where $\eta(t_i)=t_{\eta(i)}$.
Then $
\mathcal S^n
=q_1q_2^2\cdots q_{n-1}^{n-1}\,\mathrm{id}. $

For $v\in S_n$, denote by $v^L$ the left Weyl action on
$QK_T(F\ell_n)$.   
\begin{Prop}\label{prop:Seidelspec}
In $QK_T(F\ell_n)$, for any $u\in S_n$, we have
$$
\mathcal{S}(\mathcal{O}_u)=\mathcal O_\eta\star\eta^L(\mathcal O_u).
$$
\end{Prop}

\begin{proof}
By the cyclic symmetry of the diagram associated to $(u, w)$ and 
Lemma \ref{lem:shiftQpara}, we have 
$$
q_{u^{-1}(n),n}\,c^q_{\eta u,\eta w}(t,y)
=
q_{w^{-1}(n),n}\,\eta\bigl(c^q_{u,w}(t,y)\bigr).
$$ 
Consequently, we have
\begin{align*}
\mathcal S\bigl(\mathcal O_u\star\lambda_y(\mathcal V_k)\bigr)
&=\mathcal S\left(\sum_w c^q_{u,w}(t,y)\mathcal O_w\right)\\
&=\sum_w\eta\bigl(c^q_{u,w}(t,y)\bigr)
       q_{w^{-1}(n),n}\mathcal O_{\eta w}\\
&=q_{u^{-1}(n),n}
       \sum_w c^q_{\eta u,\eta w}(t,y)\mathcal O_{\eta w}\\
&=q_{u^{-1}(n),n}
       \mathcal O_{\eta u}\star\lambda_y(\mathcal V_k)\\
&=\mathcal S(\mathcal O_u)\star\lambda_y(\mathcal V_k).
\end{align*}
It follows that $\mathcal S\bigl(\mathcal O_u\star c_r^T(\mathcal V_k^\vee)\bigr)
 =\mathcal S(\mathcal O_u)\star c_r^T(\mathcal V_k^\vee)$
 for all $1\leq r\leq k<n$. These equivariant Chern classes generate a basis $\{B_\nu\}_\nu$ of $QK_T(F\ell_n)$ over $\operatorname{Frac}(R(T)[q_1^{\pm1},\cdots,q_{n-1}^{\pm1}])$
(after extending the coefficient ring), so we can write $\mathcal{O}_{u}=\sum_\nu f_\nu B_\nu$. 
 Consequently, we have
 \begin{align*}
 \mathcal S(\mathcal O_u)&=\mathcal S(\mathcal{O}_{\rm id}\star \mathcal O_u)=\mathcal S\left(\mathcal O_{\rm id}\star \sum_\nu f_\nu B_\nu\right)\\
 &=\mathcal S(\mathcal O_{\rm id})\star \sum_\nu \eta^L(f_\nu) B_\nu=\mathcal O_\eta\star  \eta^L(\mathcal{O}_{u}).
 \end{align*}
\end{proof}
\noindent We remark that    $\mathcal{S}$ is $QK_G(F\ell_n)$-linear by the above proposition, after viewing $QK_T(F\ell_n)$ as a $QK_G(F\ell_n)$-module.

Let  $\pi_j: F\ell_n\longrightarrow Gr(j, n)$ denote the    natural projection to  complex Grassmannians. We have $\mathcal{O}_{v_j}=\pi_j^*[\mathcal{O}_{\rm pt}]\in K(F\ell_n)$, where $v_j=[n-j+1, \cdots, n, 1, 2, \cdots, n-j]=\eta^{n-j}$.

\begin{Coro}\label{cor:equivSeidel}
For any $1\leq j<n$ and $u\in S_n$,  in $QK_T(F\ell_n)$  we have
$$
\mathcal O_{v_j}\star v_j^L(\mathcal O_u)
=
q^{\varpi_j^\vee-u^{-1}(\varpi_j^\vee)}
\mathcal O_{v_ju},
$$
where $\varpi_j^\vee$ is the $j$-th fundamental coweight.
\end{Coro}

\begin{proof}
By Proposition \ref{prop:Seidelspec}, we have 
\begin{align*}
\mathcal S^{n-j}(\mathcal O_u)&=\mathcal{S}^{n-j-1}\big(\mathcal{S}(\mathcal O_{\rm id}) \star \eta^L(\mathcal{O}_u)\big)\\
&=\cdots\\
&=
\mathcal S^{n-j}(\mathcal O_{\rm id})\star (\eta^{n-j})^L(\mathcal O_u).
\end{align*}
Thus we have 
\begin{equation}\label{eqn:degome}
    \mathcal O_{v_j}\star v_j^L(\mathcal O_u)
=
\frac{\prod_{a=j+1}^nq_{u^{-1}(a),n}}
     {\prod_{a=j+1}^nq_{a,n}}\,
\mathcal O_{v_ju}.
\end{equation}

The exponent of $q_r$ in this quotient is $
\#\{a>j:u^{-1}(a)\le r\}-\#\{a:j<a\le r\}$, equal to
$$
 \#\{b\le r:u(b)>j\}-\max(r-j,0)
 =\min(j,r)-\#\{b\le r:u(b)\le j\}\ge0.
$$
 Note that
$\varpi_j^\vee=\sum_{b=1}^j e_b-\frac jn\sum_{b=1}^n e_b$ and $\alpha_r^\vee=e_r-e_{r+1}$.
Therefore,
\begin{align*}
 \varpi_j^\vee-u^{-1}(\varpi_j^\vee)
 &=\sum_{b=1}^j (e_b-e_n)-\sum_{u(b)\leq j} (e_b-e_n)\\
 &=\sum_{r=1}^{n-1}
   \left(\min(j,r)-\#\{b\le r:u(b)\le j\}\right)\alpha_r^\vee.
\end{align*}
\end{proof}

\subsection{Equivariant quantum Pieri rule for $QK_T(F\ell_{n_1, \cdots, n_r; n})$}
The following interpretation is a typical application of the Seidel product: 
 $$\mbox{quantum Pieri rule } = \mbox{ classical Pieri rule } + \mbox{ Seidel action}.$$
    This has been shown for the quantum $K$-theory of complex Grassmannians in \cite{LLSY} and minuscule Grassmannians of classical Lie type in \cite{BCP}. 
In this subsection, we provide the Pieri rule for the equivariant quantum $K$-theory of partial flag varieties, also  using the principle ``quantum=classical+ Seidel". Its specialization to $QK(F\ell_n)$   
          resolves \cite[Conjecture 3.1 and Problem 3.2]{LiSo} of the third author and Song. 
 
Fix any $1\leq n_1<n_2<\cdots <n_r< n$. We consider the $r$-step partial flag variety   
$$X=F\ell_{n_1, \cdots, n_r; n} = \left\{
 V_{n_1}\subset \cdots \subset V_{n_r}\subset \mathbb{C}^n: \dim V_{n_j}=n_j\right\}.$$
 It is a quotient $G/P$ of $G=GL_n$ by a   parabolic subgroup $P$.
There are  (opposite) Schubert varieties $X^w=\overline{B^-wP/P}\subset G/P$ of codimension $\ell(w)$, indexed by permutations $w\in W^P\subset S_n$ that have no descents except possibly at $n_1, \cdots, n_r$. Denote 
 $\mathcal{O}_w^X:=[\mathcal{O}_{X^w}]\in K_T(G/P)$, and consider the  (localized) $T$-equivariant quantum $K$-theory $QK_T(G/P) = \big(K_T(G/P)\otimes_{\mathbb{Z}}\mathbb{C}[Q_1^{\pm 1},\cdots, Q_r^{\pm 1}], \star \big)$ with the quantum product    
 $$\mathcal{O}_u^X \star \mathcal{O}_v^X=\sum_{w\in W^P, d\in H_2(G/P, \mathbb{Z})} \bar N_{u, v}^{w, d}\mathcal{O}_w^X {Q}^d.$$
Set $n_0:=0$, $n_{r+1}:=n$. The first Chern class and  the degree of $Q_i$   are respectively given by 
$$c_1(TX)=\sum_{i=1}^r (n_{i+1}-n_{i-1}){\rm PD}[X^{s_{n_i}}];\quad \deg Q_i= n_{i+1}-n_{i-1},\qquad i=1, \cdots, r.$$

 Let $\pi: G/B=F\ell_n\longrightarrow F\ell_{n_1, \cdots, n_r; n}=G/P$ denote the natural projection.  The following is a special property of quantum $K$-theory, proved by Kato and generalizing the special case $G/B\to {\rm pt}$ in \cite{BCLM}.

\begin{Prop}[\protect{\cite[Theorem 2.19]{Kato3}}]\label{prop:Phimorph}
    There is a morphism $\Phi$ of $R(T)$-algebras,
     $$\Phi: QK_T(G/B)\longrightarrow QK_T(G/P),$$
     defined by $\Phi(\mathcal O_u)=\pi_*\mathcal O_u$, and $\Phi(q_j)=Q_{i}$ if $j=n_i$, or 1 otherwise.
\end{Prop}

For $u\in W$, we denote by $\bar u\in W^P$ the minimal length representative of the coset $uW_P$, where $W_P \leq S_n$ is the subgroup generated by those $s_j$ with $j\neq n_i$ for any $i$. Equivalently, $\bar u$ is obtained from $u$ by sorting $\{u(n_{i-1}+1), u(n_{i-1}+2), \cdots, u(n_i)\}$ increasingly for all $1\leq i\leq r+1$.

\begin{Lemma}\label{lemma:SeiPhi}
    For any $u, v\in S_n$, the following hold:
   {\upshape $$ \mbox{(i) }\, \overline{u \bar  v}=\overline{uv};\qquad\mbox{(ii) }\, \Phi(q_{\bar u^{-1}(n), n})=\Phi(q_{u^{-1}(n), n});\qquad \mbox{(iii) }\, \Phi(v^L(\mathcal{O}_u))= v^L(\mathcal O_{\bar u}^X). $$ }
\end{Lemma}

\begin{proof}
Since $\bar v W_P=v W_P$,   $u\bar v W_P=uv W_P$. Thus (i) holds. Write $u=\bar u u'$ with $u'\in W_P$. Then $u^{-1}(n)=(u')^{-1}\bar u^{-1}(n)$ and $\bar u^{-1}(n)$ belong to the same block $(n_{i-1}, n_i]$ for some $i$. Thus (ii) holds.
Note that 
$$
\Phi(v^L(\mathcal{O}_u))=\pi_*v^L(\mathcal O_u)
=[R\pi_*\mathcal O_{ vX^u}]
=[\mathcal O_{ vX^{\bar u}}]
=v^L(\mathcal O_{\bar u}^X),$$
which implies (iii).
\end{proof}
Recall $\eta=(12\cdots n)$.  Define the Seidel operator $\mathcal S_X:QK_T(G/P)\longrightarrow QK_T(G/P)$ by 
\begin{align}
\mathcal S_X\!\big(f_u\mathcal O_u^X\big)
=\eta(f_u)\Phi(q_{u^{-1}(n),n})\mathcal O_{\overline{\eta u}}^X.
\end{align}

\begin{Prop}\label{prop:Seidpara}  We have $\mathcal{S}_X\Phi=\Phi \mathcal S$. Moreover, for any $m\in \mathbb{Z}_{>0}$ and $u, v\in W^P$, 
\begin{align*}
\mathcal S_X^m(\mathcal O_v^X\star \mathcal O_u^X)
&=\Phi(\mathcal S^m(\mathcal O_v)\star (\eta^m)^L(\mathcal O_u))\\
&=\mathcal S_X^m(\mathcal O_v^X)\star (\eta^m)^L(\mathcal O_u^X).
\end{align*}    
\end{Prop}  
\begin{proof}
For any $w\in S_n$, by Lemma \ref{lemma:SeiPhi}, we have  
\begin{align*}
\mathcal{S}_X\Phi(\mathcal{O}_w)&=\mathcal S_X(\mathcal O_{\bar w}^X)=\Phi(q_{\bar w^{-1}(n), n}) \mathcal O_{\overline{\eta \bar w}}^X\\
&=\Phi(q_{w^{-1}(n), n}) \mathcal O_{\overline{\eta w}}^X\\
&=\Phi \mathcal S(\mathcal{O}_w).
\end{align*}
Combining it with Proposition \ref{prop:Phimorph}, we obtain
\begin{align*}
\mathcal S_X^m(\mathcal O_v^X\star \mathcal O_u^X)&=\mathcal S_X^m\Phi(\mathcal O_v\star \mathcal O_u)=
    \Phi\mathcal S^m(\mathcal O_v\star \mathcal O_u)\\
&=\Phi\left(\mathcal S^m(\mathcal O_v)\star (\eta^m)^L(\mathcal O_u)\right). 
\end{align*}
    The last equality also   follows in the same way, using Lemma   \ref{lemma:SeiPhi} (iii).  
\end{proof}

\begin{Lemma}\label{lemma:noquantum} Let $1\leq k<n$ and   $u\in S_n$ with $u(k+1)=n$. Then for any $1\leq i\leq k$ and $v\in S_n$, in $QK_T(F\ell_n)$, we have
 $$ \mathcal{O}_u\star v^L(\mathcal{O}_{s_is_{i+1}\cdots s_k})=\mathcal{O}_u \cdot v^L(\mathcal{O}_{s_is_{i+1}\cdots s_k}).$$    
\end{Lemma}

\begin{proof}
Let  $w\in S_n$ with $c^q_{u,w}(t,y)\neq 0$. 
By Lemma \ref{lem:onewallnopath},  $I=J=\varnothing$.   Thus by Theorem \ref{thm:main-formula}, we have $
 \mathcal O_u\star\lambda_y(\mathcal V_k)
 =
 \mathcal O_u\cdot\lambda_y(\mathcal V_k).
$
Comparing coefficients of powers of $1 + y$ in 
$ \lambda_y(\mathcal V_k)
 =
 \sum_{j=0}^k(-1)^j(1+y)^{k-j}
 [\textstyle\bigwedge^j\mathcal V_k]$,
 gives
 $
 \mathcal O_u\star[\textstyle\bigwedge^j\mathcal V_k]
  =
  \mathcal O_u\cdot[\textstyle\bigwedge^j\mathcal V_k] $ for each $0\le j\le k$.
Since $\mathcal{V}_k$ is $G$-equivariant, $v^L$ fixes  $[\textstyle\bigwedge^j\mathcal V_k]$.
Hence, by Lemma \ref{lemma:Olambda},
$$v^L(\mathcal O_{s_is_{i+1}\cdots s_k}) =
 \sum_{j=0}^k(-1)^j
 v\left(\chi_T(\check{\mathbb P}^{i-1},L^j)\right)
 [\textstyle\bigwedge^j\mathcal V_k].
$$ Thus the statement holds.  
\end{proof}

For any $u, w\in S_n$ and $1\leq k<n$, we denote
\begin{equation}
    \delta_k(u, w):=\overline{v_hw}^{-1}v_h(\varpi^\vee_{h})-\overline{v_hu}^{-1}v_h(\varpi^\vee_{h})\in Q^\vee\in H_2(G/B, \mathbb{Z}),
\end{equation}
where $h=u(k+1)$ and $\varpi_n^\vee:=0$ by convention, and we recall $v_h=\eta^{n-h}$. Note $Q^{\pi_*(\delta_k(u, w))}=\Phi(q^{\delta_k(u, w)})$. The following theorem shows that $\mathcal O_w^X$ occurs in at most one quantum term in the quantum product $\mathcal{O}_u^X\star \mathcal O_{s_is_{i+1}\cdots s_{n_j}}^X$, and the  coefficient equals a specific coefficient in the Pieri rule for $QK_T(F\ell_n)$. 
\begin{Th}\label{thm:parttofull} 
Let $k\in \{n_1, \cdots, n_r\}$ and $u, w\in W^P$. Denote $h=u(k+1)$. Set $\hat u=v_h^{-1}\overline{v_hu}$ and $\hat w=v_h^{-1}\overline{v_hw}$.
For any $1\leq i\leq k$, we have 
 $$ \sum_{d\in H_2(G/P, \mathbb{Z})}\bar N_{u, s_is_{i+1}\cdots s_k}^{w, d}Q^d= N_{\hat u, s_is_{i+1}\cdots s_k}^{\hat w, \delta_k(u, w)} Q^{\pi_*(\delta_k(u, w))}. $$ 
\end{Th}
\begin{proof}
Denote $\mathsf s:=s_is_{i+1}\cdots s_k$. Note $v_hu(k+1)=n$. By Lemma \ref{lemma:noquantum}, we have 
\begin{align}\label{eqn:redformain}
        \Phi(\mathcal O_{{v_h u}}\star v_h^L(\mathcal O_{\mathsf s}))= \Phi(\mathcal O_{{v_h u}}\cdot v_h^L(\mathcal O_{\mathsf s})) 
       =
 \mathcal O_{\overline{v_hu}}^X
 \cdot v_h^L(\mathcal O_{\mathsf s}^X).
\end{align}
Then by Proposition \ref{prop:Seidpara}, we have
\begin{align*}
\mathcal{S}_X^{n-h}(\mathcal{O}_u^X\star \mathcal O_{\mathsf s}^X)
    &=   \Phi\left(\mathcal{S}^{n-h}(\mathcal{O}_u)\star v_h^L(\mathcal O_{\mathsf s})\right)\nonumber\\
      & =
 \Phi\left(\prod_{a=h+1}^n q_{u^{-1}(a),n}\right)
 \mathcal O_{\overline{v_hu}}^X
 \cdot v_h^L(\mathcal O_{\mathsf s}^X).
\end{align*}
The left-hand side of the above formula,
after  applying $\mathcal{S}_X^{n-h}$ on the expansion, is given by  
\begin{align}\label{eqn:redbarN}
     \sum_{w, d} v_h(\bar N_{u, \mathsf s}^{w, d})Q^d\Phi\left(\prod_{a=h+1}^n q_{w^{-1}(a), n}\right) \mathcal O_{\overline{v_h w}}^X .
\end{align}
The above right-hand side can be read off from $ \mathcal O_{\overline{v_hu}}
 \cdot v_h^L(\mathcal O_{\mathsf s})$, since the  $R(T)$-algebra morphism $\pi^*: K_T(G/P)\to K_T(G/B)$ satisfies $\pi^*(\mathcal{O}_w^X)=\mathcal{O}_w$ and $v^L\circ \pi^*=\pi^*\circ v^L$.  
 
 Applying the Seidel product again, we have  
\begin{equation}\label{eqn:SOhatu}
  \mathcal S^{n-h}(\mathcal O_{\hat u}\star \mathcal O_{\mathsf s})=\mathcal S^{n-h}(\mathcal O_{\hat u})\star v_h^L(\mathcal O_{\mathsf s})=\prod_{a=h+1}^n q_{\hat u^{-1}(a),n}\mathcal O_{\overline{v_hu}}
\star v_h^L(\mathcal O_{\mathsf s}).
\end{equation}
 Hence,
$\mathcal O_{\overline{v_hu}}
 \cdot v_h^L(\mathcal O_{\mathsf s})$
is obtained by taking the coefficient of
$\prod_{a=h+1}^n q_{\hat u^{-1}(a),n}$
in the left-hand side of \eqref{eqn:SOhatu}. 
Since $v_h\hat w=\overline{v_hw}$, the coefficient of $\mathcal O_{\overline{v_hw}}$ in the classical product is read off from the left-hand side of \eqref{eqn:SOhatu}, which is 
$
v_h\big(N_{\hat u,\mathsf s}^{\hat w,\delta}\big),
$ where 
$$q^\delta=
\frac{\prod_{a=h+1}^n q_{\hat u^{-1}(a),n}}
     {\prod_{a=h+1}^n q_{\hat w^{-1}(a),n}}
=
q^{\hat w^{-1}(\varpi_h^\vee)
   -\hat u^{-1}(\varpi_h^\vee)}
=
q^{\delta_k(u,w)}$$
by Corollary \ref{cor:equivSeidel}.
Noting  that $w\mapsto\overline{v_hw}$ is a bijection on $W^P$, we have 
\begin{equation}\label{eqn:redtoful}
\mathcal O_{\overline{v_hu}}^X
 \cdot v_h^L(\mathcal O_{\mathsf s}^X)
=
\sum_{w\in W^P}
v_h\bigl(
 N_{\hat u,\mathsf s}^{\hat w,\delta_k(u,w)}
\bigr)
\mathcal O_{\overline{v_hw}}^X.
\end{equation} 
Comparing \eqref{eqn:redbarN} with
\eqref{eqn:redtoful} multiplied by
$\Phi\left(\prod_{a=h+1}^n q_{u^{-1}(a),n}\right)$ and applying $v_h^{-1}$ on both sides, we obtain
\begin{equation*}
\sum_{d\ge0}
\bar N_{u,\mathsf s}^{w,d}Q^d
=
 N_{\hat u,\mathsf s}^{\hat w,\delta_k(u,w)}
\Phi\left(
\frac{\prod_{a=h+1}^n q_{u^{-1}(a),n}}
     {\prod_{a=h+1}^n q_{w^{-1}(a),n}}
\right).
\end{equation*}
 Since $\hat uW_P=v_h^{-1}\overline{v_h u} W_P=v_h^{-1}{v_h u} W_P= uW_P$ and similarly $\hat wW_P= wW_P$, by  Lemma \ref{lemma:SeiPhi} (ii),  the statement follows from
$$
\Phi\left(
\frac{\prod_{a=h+1}^n q_{u^{-1}(a),n}}
     {\prod_{a=h+1}^n q_{w^{-1}(a),n}}
\right)
=
\Phi\left(
\frac{\prod_{a=h+1}^n q_{\hat u^{-1}(a),n}}
     {\prod_{a=h+1}^n q_{\hat w^{-1}(a),n}}
\right)
=Q^{\pi_*(\delta_k(u,w))}.
$$ 
\end{proof}

\begin{Lemma}\label{lemma:signcorr} With the notation in Theorem \ref{thm:parttofull} and assuming $u(k+1)<n$, we have  
    $$\ell(u)+\ell(w)+\deg Q^{\pi_*(\delta_k(u,w))} \equiv \ell(\hat u)+\ell(\hat w) \quad \mod 2.$$
\end{Lemma}

 \begin{proof}
Let $R$ denote the root system, and  $R_P$ (resp. $R_h$) denote the root subsystem generated by the complement of $\{\alpha_{n_1}, \cdots, \alpha_{n_r}\}$ (resp. $\{\alpha_h\}$) in the set of simple roots. Use the superscript $+$ for the subset of positive roots, and set
 $\rho_P={1\over 2}\sum_{\alpha\in R_P^+}\alpha$.
Note that, for $\alpha\in R^+$, $v_h(\alpha)<0$ if and only if $\alpha\in R^+\setminus R_h^+$.   

Now take $\alpha\in R_P^+$. Since 
$v_h\hat w=\overline{v_hw}\in W^P$, $v_h\hat w(\alpha)>0$.
If $\hat w(\alpha)>0$,  then $\hat w(\alpha)\in R_h^+$, implying $\langle \hat w\alpha, \varpi_h^\vee\rangle
=0$. If $\hat w(\alpha)<0$, then   $ v_h\bigl(-\hat w(\alpha)\bigr)<0$, implying  $\langle \hat w\alpha, \varpi_h^\vee\rangle
=-1$. 
 Since $w$ is the minimal representative of $\hat wW_P$,
the following holds (for any element $w$ of $W^P$)
$$
\ell(\hat w)-\ell(w)
=\sharp\{\alpha\in R_P^+ \mid \hat w(\alpha)<0\}=
-\langle 2\rho_P,\hat w^{-1}(\varpi_h^\vee)\rangle.
$$
By \cite[Lemma 3.5]{FuWo}, we have  
$
\deg Q^{\pi_*(\delta_k(u,w))}
=
\langle 2\rho-2\rho_P,\delta_k(u,w)\rangle$. 
Hence, we have 
\begin{align*}
&\ell(\hat u)-\ell(u)-\ell(\hat w)+\ell(w)+\deg Q^{\pi_*(\delta_k(u,w))}\\
&\quad= \langle 2\rho_P,\hat w^{-1}(\varpi_h^\vee)-\hat u^{-1}(\varpi_h^\vee)\rangle+
\langle 2\rho-2\rho_P,\delta_k(u,w)\rangle\\
&\quad=\langle 2\rho,\delta_k(u,w)\rangle \equiv 0 \mod 2.
\end{align*} 
\end{proof}
With $\hat u$ and $\hat w$ as in Theorem \ref{thm:parttofull}, define  $\hat I=I_k(\hat u, \hat w), \hat J=J_k(\hat u, \hat w)$, and  
\begin{align*}
\hat{\mathsf F}=\{j\in \hat u([k]): \hat u^{-1}(j)=\hat w^{-1}(j)\}, \qquad  \hat{\mathsf U}=\{j\in \hat u([k]): \hat u^{-1}(j)< \hat w^{-1}(j)\}.
\end{align*}
\begin{Th}\label{thm:equiPierPara}
    Let $k\in\{n_1, \cdots, n_r\}$ and $u, w\in W^P$. For any $1\leq i\leq k$, all $\bar N_{u, s_is_{i+1}\cdots s_k}^{w, d}$ vanish unless  $c_{\hat u, \hat w}^q(t, y)\neq 0$ and $q_{\hat I\hat J}=q^{\delta_k(u, w)}$.    When both conditions hold, only $d=\pi_*(\delta_k(u, w))$ can contribute, and   
\begin{align*} 
& (-1)^{\ell(u)+k-i+1+\ell(w)+\deg Q^{\pi_*(\delta_k(u, w))}}
 \bar N_{u,s_{i}s_{i+1}\cdots s_k}^{w, \pi_*(\delta_k(u, w))}\\
&\quad =\binom{|\hat{\mathsf U}|-1}{i-1-|\hat{\mathsf F}|}
+ \sum_{\substack{S\subseteq \hat{\mathsf U}~:~|S|\ge i-|\hat{\mathsf F}|}}
\psi_{|\hat{\mathsf F}|+|S|-i+1}({\hat{\mathsf F}\cup S}).
\end{align*}
\end{Th}

\begin{proof}
    When $u(k+1)=n$, we have $v_h=\mbox{id}$, so that $\hat u=u$ and $\hat w=w$. The quantum product in \eqref{eqn:redformain} has no quantum terms, so $N_{\hat u, s_is_{i+1}\cdots s_k}^{\hat w, \delta_k(u, w)}\neq 0$   only if $\delta_k(u, w)=0$. Thus there is no sign issue in this case.

    When $u(k+1)<n$, the statement follows immediately from Theorem \ref{thm:equivPieri}, Theorem  \ref{thm:parttofull} and Lemma \ref{lemma:signcorr}.
\end{proof}

\appendix
\section{The two-interval cases}\label{app:two-column}\label{sec:twocolproof}

We omit the balance conditions. All diagrams are understood to be balanced.
\newcounter{appcase}
\newcommand{\appcaseanchor}[2]{\renewcommand{\theappcase}{#1}\refstepcounter{appcase}\label{app:case-#2}}
\eyesetup{interval width=12mm,level distance=2.5mm,first level=3.5mm,
  endpoint gap=1.5mm,label offset=.18,label font=\scriptsize,
  wall label sep=1.5pt,wall padding=1.5mm,path width=.55pt,
  bundle width=1.5pt,arrow length=1.4mm,arrow width=1mm,
  bundle arrow length=2mm,bundle arrow width=1.5mm}

\subsection{Case 1} When $a<b\leq k$ and $c<d\leq k$, there are 13 cases.

\appcaseanchor{1.1}{1-1}
\subsubsection*{Case 1.1} When $a<b<c<d\leq k$, only one nonzero case:
\[

\]

\appcaseanchor{1.2}{1-2}
\subsubsection*{Case 1.2} When $a<c<b<d\leq k$, only one nonzero case:
\[
%
\]

\appcaseanchor{1.3}{1-3}
\subsubsection*{Case 1.3} When $a<c<d<b\leq k$, all four weights are zero.
\appcaseanchor{1.4}{1-4}
\subsubsection*{Case 1.4} When $c<a<b<d\leq k$, only one nonzero case:
\[
%
\]

\appcaseanchor{1.5}{1-5}
\subsubsection*{Case 1.5} When $c<a<d<b\leq k$, all four weights are zero.
\appcaseanchor{1.6}{1-6}
\subsubsection*{Case 1.6} When $c<d<a<b\leq k$, all four weights are zero.
\appcaseanchor{1.7}{1-7}
\subsubsection*{Case 1.7} When $a=c<b<d\leq k$, only one nonzero case:
\[
%
\]

\appcaseanchor{1.8}{1-8}
\subsubsection*{Case 1.8} When $a=c<d<b\leq k$, only one nonzero case:
\[
%
\]

\appcaseanchor{1.9}{1-9}
\subsubsection*{Case 1.9} When $a<c<b=d\leq k$, only one nonzero case:
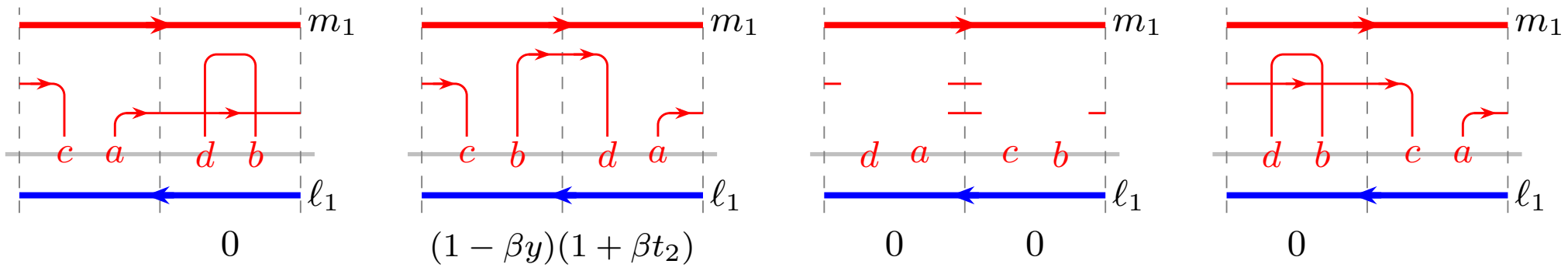
\[
%
\]

\appcaseanchor{1.10}{1-10}
\subsubsection*{Case 1.10} When $c<a<b=d\leq k$, only one nonzero case:
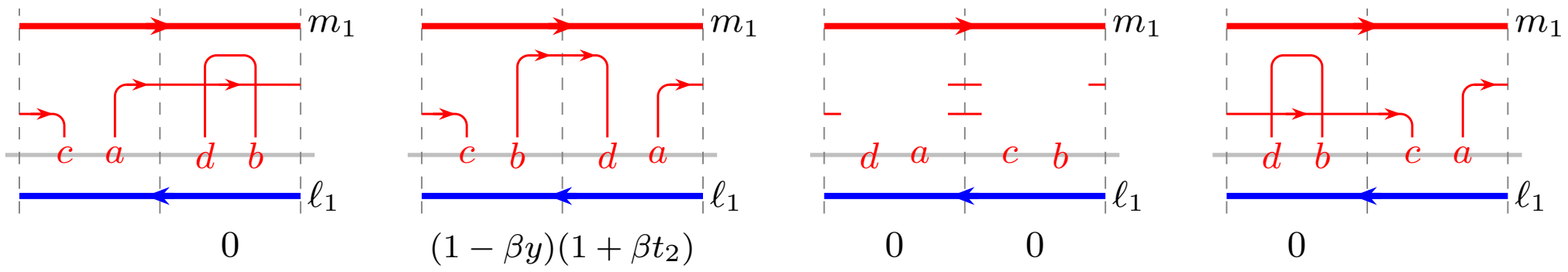
\[
%
\]

\appcaseanchor{1.11}{1-11}
\subsubsection*{Case 1.11} When $a<b=c<d\leq k$, only one nonzero case:
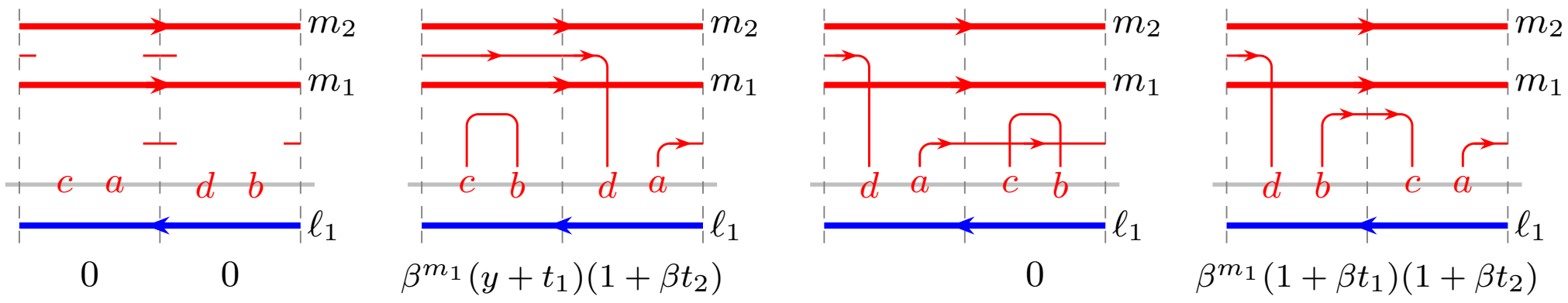
\[
%
\]

\appcaseanchor{1.12}{1-12}
\subsubsection*{Case 1.12} When $c<a=d<b\leq k$, only one nonzero case:
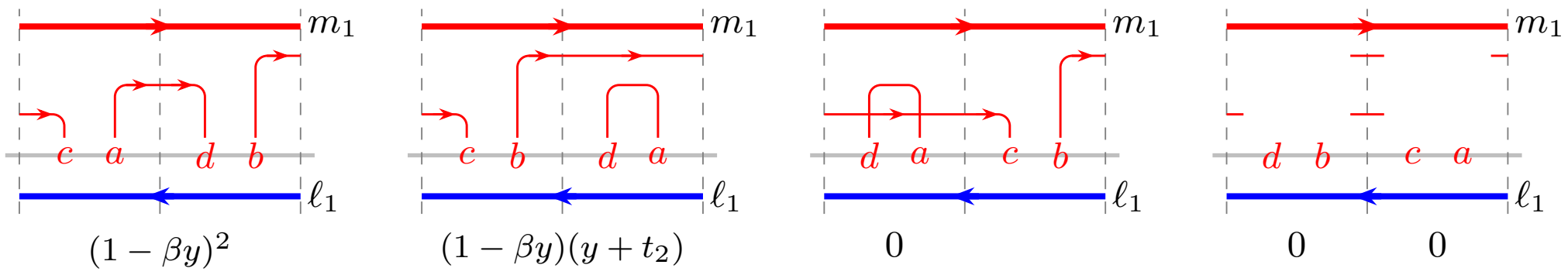
\[
%
\]

\appcaseanchor{1.13}{1-13}
\subsubsection*{Case 1.13} When $a=c<b=d\leq k$, only one nonzero case:
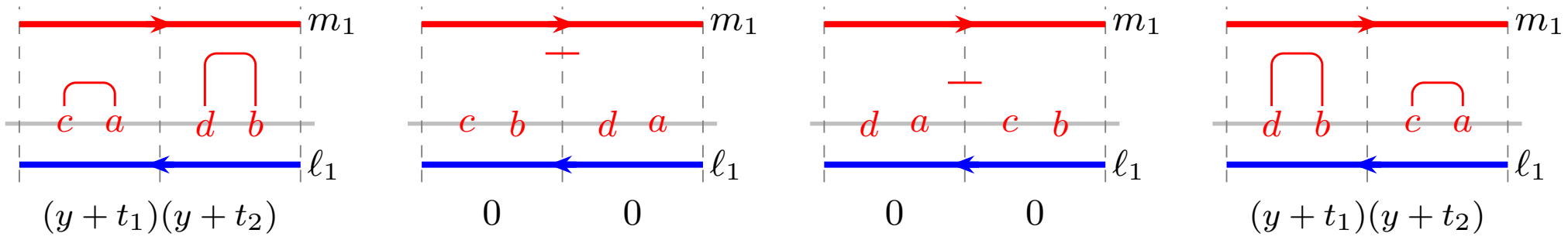
\[
%
\]

\subsection{Case 2} When $a<b\leq k$ and $c\leq k<d$, there are 5 cases.

\appcaseanchor{2.1}{2-1}
\subsubsection*{Case 2.1} When $c<a<b\leq k<d$, only one nonzero case:
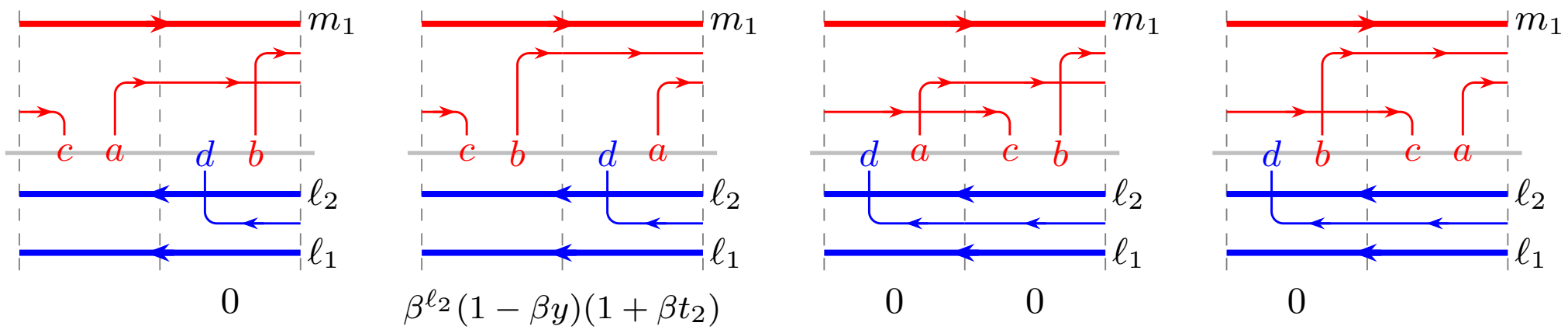
\[
%
\]

\appcaseanchor{2.2}{2-2}
\subsubsection*{Case 2.2} When $a<c<b\leq k<d$, only one nonzero case:
\[
%
\]

\appcaseanchor{2.3}{2-3}
\subsubsection*{Case 2.3} When $a<b<c\leq k<d$, only one nonzero case:
\[
%
\]

\appcaseanchor{2.4}{2-4}
\subsubsection*{Case 2.4} When $a=c<b\leq k<d$, only one nonzero case:
\[
%
\]

\appcaseanchor{2.5}{2-5}
\subsubsection*{Case 2.5} When $a<b=c\leq k<d$, only one nonzero case:
\[
%
\]

\subsection{Case 3} When $a\leq k<b$ and $c<d\leq k$, there are 5 cases.

\appcaseanchor{3.1}{3-1}
\subsubsection*{Case 3.1} When $a<c<d\leq k<b$, only one nonzero case:
\[
%
\]

\appcaseanchor{3.2}{3-2}
\subsubsection*{Case 3.2} When $c<a<d\leq k<b$, only one nonzero case:
\[
%
\]

\appcaseanchor{3.3}{3-3}
\subsubsection*{Case 3.3} When $c<d<a\leq k<b$, only one nonzero case:
\[
%
\]

\appcaseanchor{3.4}{3-4}
\subsubsection*{Case 3.4} When $a=c<d\leq k<b$, only one nonzero case:
\[
%
\]

\appcaseanchor{3.5}{3-5}
\subsubsection*{Case 3.5} When $c<a=d\leq k<b$, only one nonzero case:
\[
%
\]

\subsection{Case 4} When $a\leq k<b$ and $c\leq k<d$, there are 9 cases.
\appcaseanchor{4.1}{4-1}
\subsubsection*{Case 4.1} When $a<c\leq k<b<d$, only one nonzero case:
\[
%
\]

\appcaseanchor{4.2}{4-2}
\subsubsection*{Case 4.2} When $a<c\leq k<d<b$, only one nonzero case:
\[
%
\]

\appcaseanchor{4.3}{4-3}
\subsubsection*{Case 4.3} When $c<a\leq k<b<d$, only one nonzero case:
\[
%
\]

\appcaseanchor{4.4}{4-4}
\subsubsection*{Case 4.4} When $c<a\leq k<d<b$, only one nonzero case:
\[
%
\]

\appcaseanchor{4.5}{4-5}
\subsubsection*{Case 4.5} When $a=c\leq k<b<d$, only one nonzero case:
\[
%
\]

\appcaseanchor{4.6}{4-6}
\subsubsection*{Case 4.6} When $a=c\leq k<d<b$, two nonzero cases:
\[
%
\]

\appcaseanchor{4.7}{4-7}
\subsubsection*{Case 4.7} When $a<c\leq k<b=d$, two nonzero cases:
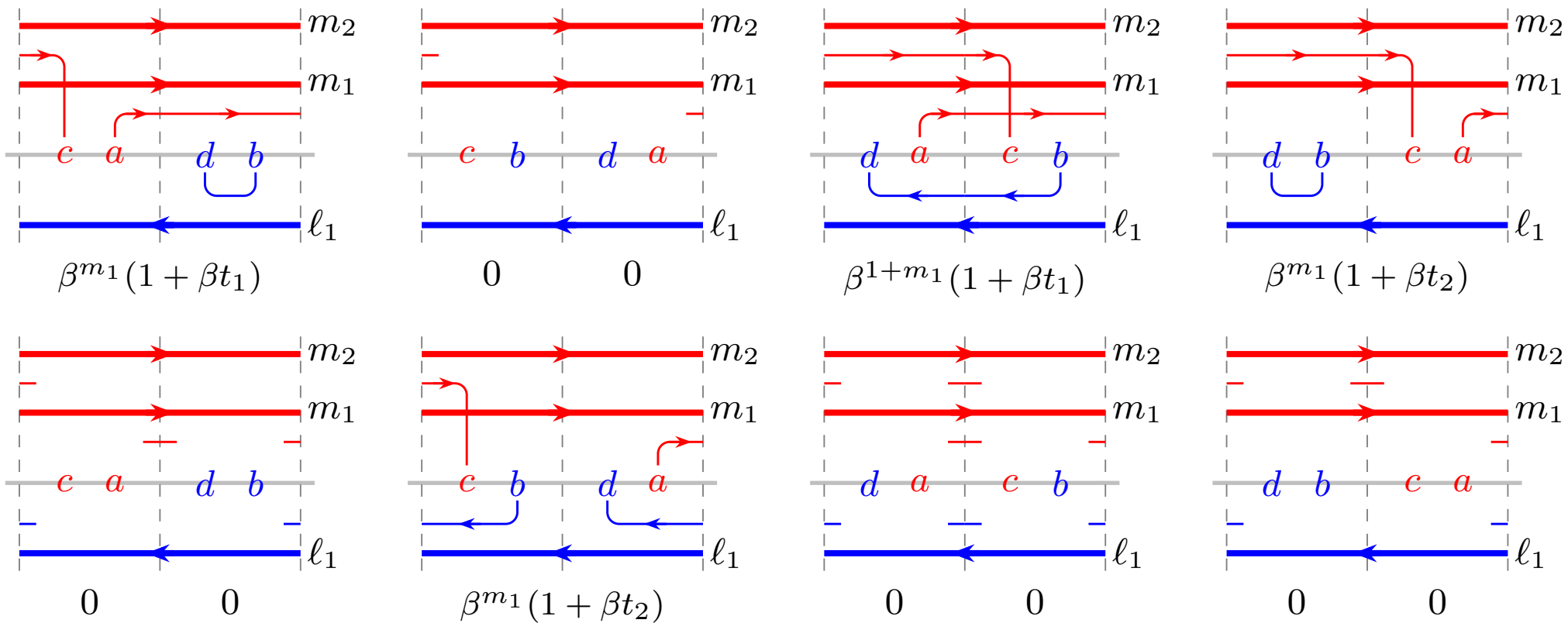
\[
%
\]

\appcaseanchor{4.8}{4-8}
\subsubsection*{Case 4.8} When $c<a\leq k<b=d$, two nonzero cases:
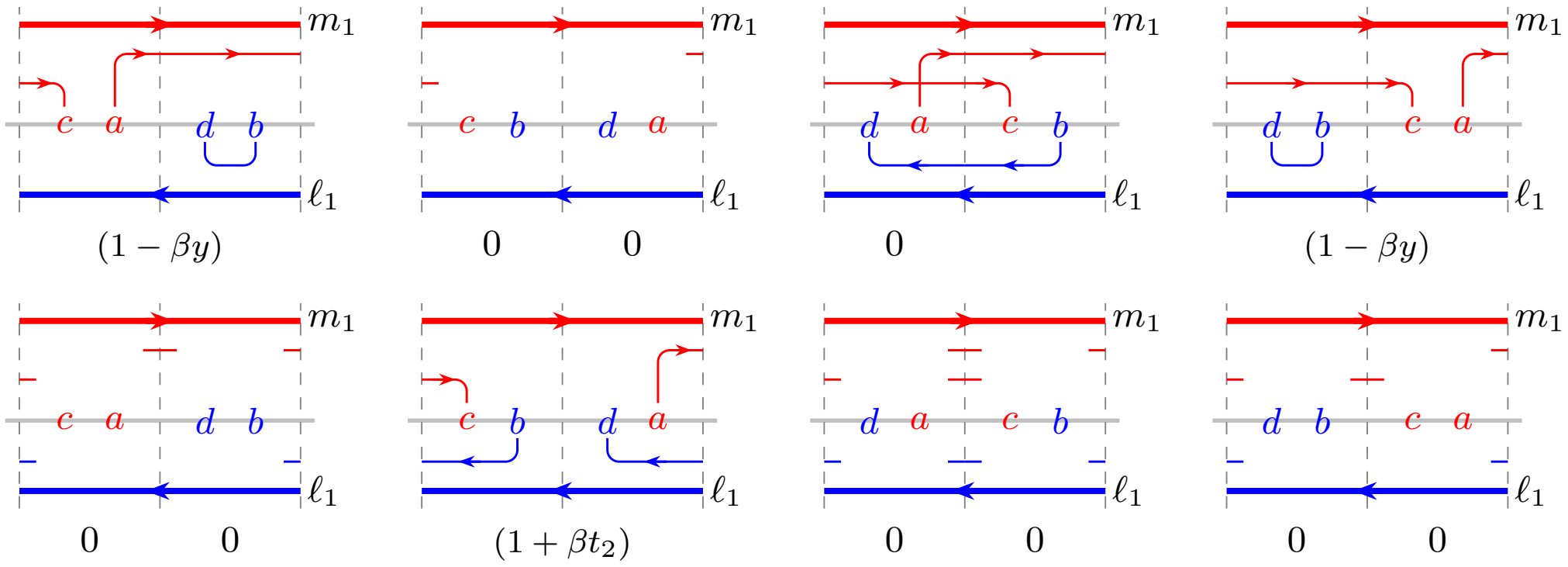
\[
%
\]

\appcaseanchor{4.9}{4-9}
\subsubsection*{Case 4.9} When $a=c\leq k<b=d$, two nonzero cases:
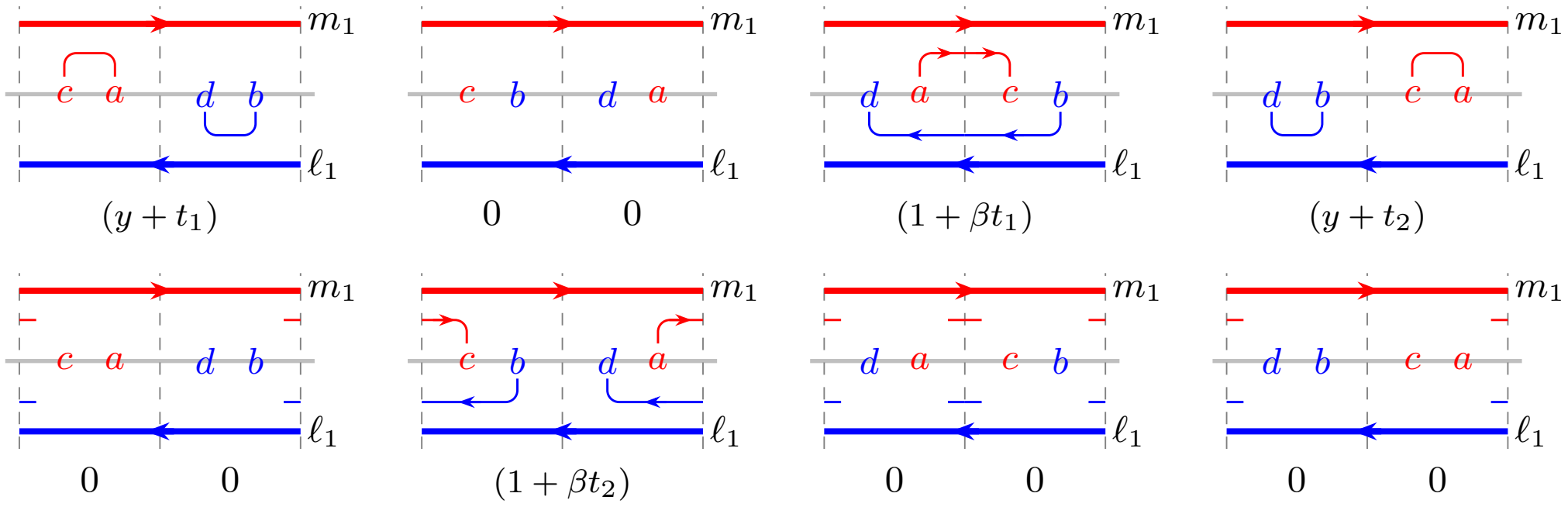
\[
%
\]

\subsection{Case 5} When $c<d\leq k<a<b$, there is one case.

\appcaseanchor{5.1}{5-1}
\subsubsection*{Case 5.1} When $c<d\leq k<a<b$, only one nonzero case:
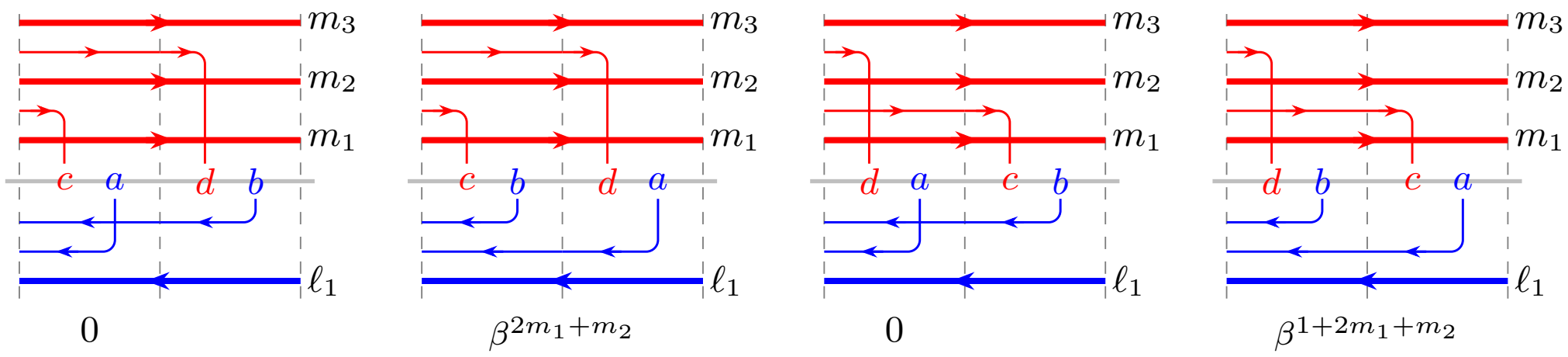
\[
%
\]

\subsection{Case 6} When $k<a<b$ and $c\leq k<d$, there are 5 cases.

\appcaseanchor{6.1}{6-1}
\subsubsection*{Case 6.1} When $c\leq k<d<a<b$, only one nonzero case:
\[
%
\]

\appcaseanchor{6.2}{6-2}
\subsubsection*{Case 6.2} When $c\leq k<a<d<b$, all four weights are zero.

\appcaseanchor{6.3}{6-3}
\subsubsection*{Case 6.3} When $c\leq k<a<b<d$, all four weights are zero.

\appcaseanchor{6.4}{6-4}
\subsubsection*{Case 6.4} When $c\leq k<a=d<b$, only one nonzero case:
\[
%
\]

\appcaseanchor{6.5}{6-5}
\subsubsection*{Case 6.5} When $c\leq k<a<b=d$, only one nonzero case:
\[
%
\]

\subsection{Case 7} When $a<b\leq k<c<d$, there is one case.

\appcaseanchor{7.1}{7-1}
\subsubsection*{Case 7.1} When $a<b\leq k<c<d$, only one nonzero case:
\[
%
\]

\subsection{Case 8} When $a\leq k<b$ and $k<c<d$, there are 5 cases.

\appcaseanchor{8.1}{8-1}
\subsubsection*{Case 8.1} When $a\leq k<c<d<b$, only one nonzero case:
\[
%
\]

\appcaseanchor{8.2}{8-2}
\subsubsection*{Case 8.2} When $a\leq k<c<b<d$, only one nonzero case:
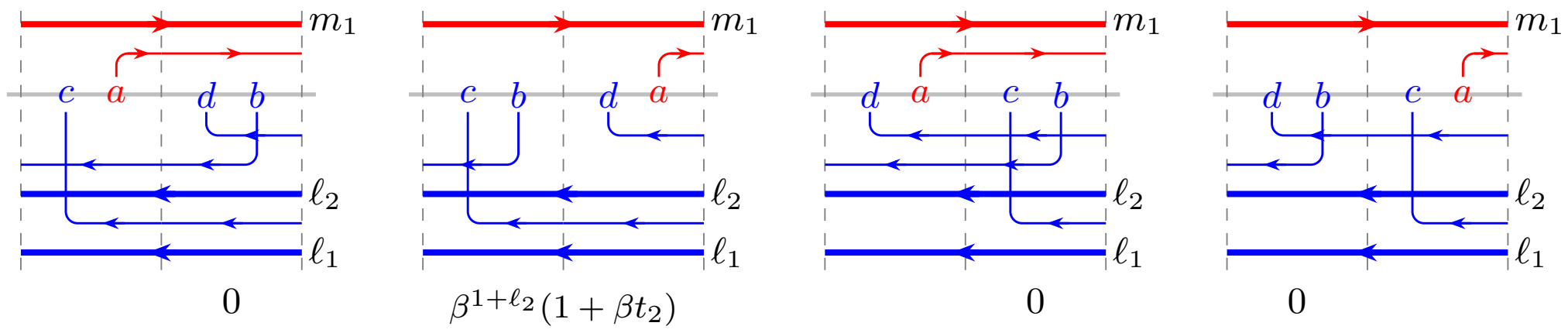
\[
%
\]

\appcaseanchor{8.3}{8-3}
\subsubsection*{Case 8.3} When $a\leq k<b<c<d$, all four weights are zero.

\appcaseanchor{8.4}{8-4}
\subsubsection*{Case 8.4} When $a\leq k<b=c<d$, only one nonzero case:
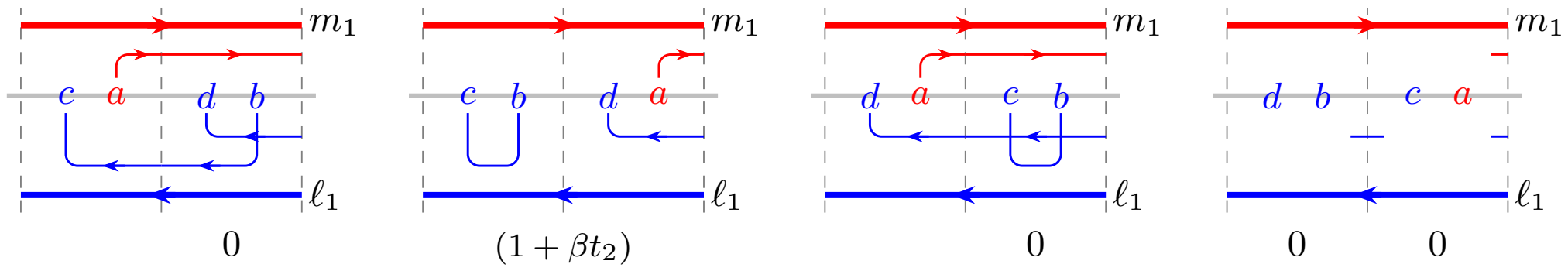
\[
%
\]

\appcaseanchor{8.5}{8-5}
\subsubsection*{Case 8.5} When $a\leq k<c<b=d$, two nonzero cases:
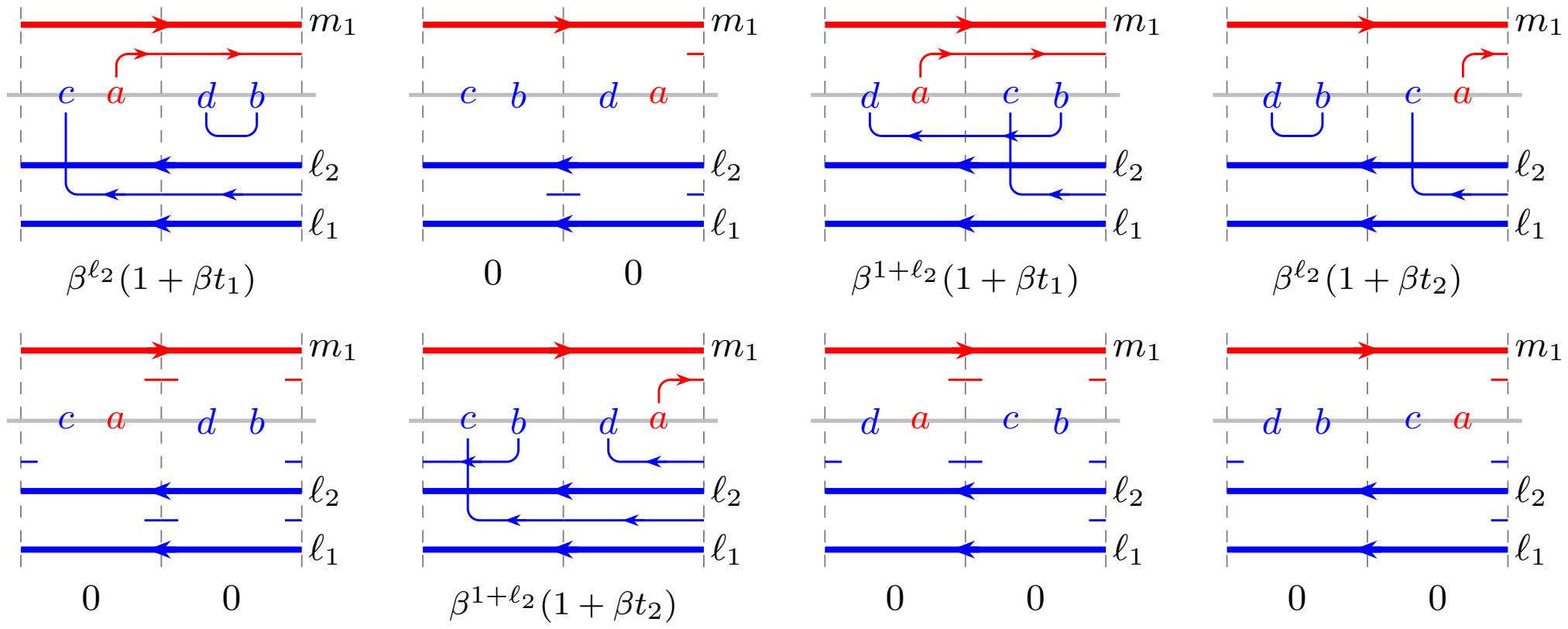
\[
%
\]

\subsection{Case 9} When $k<a<b$ and $k<c<d$, there are 13 cases.

\appcaseanchor{9.1}{9-1}
\subsubsection*{Case 9.1} When $k<a<b<c<d$, all four weights are zero.

\appcaseanchor{9.2}{9-2}
\subsubsection*{Case 9.2} When $k<a<c<b<d$, all four weights are zero.

\appcaseanchor{9.3}{9-3}
\subsubsection*{Case 9.3} When $k<a<c<d<b$, all four weights are zero.

\appcaseanchor{9.4}{9-4}
\subsubsection*{Case 9.4} When $k<c<a<b<d$, all four weights are zero.

\appcaseanchor{9.5}{9-5}
\subsubsection*{Case 9.5} When $k<c<a<d<b$, all four weights are zero.

\appcaseanchor{9.6}{9-6}
\subsubsection*{Case 9.6} When $k<c<d<a<b$, only one nonzero case:
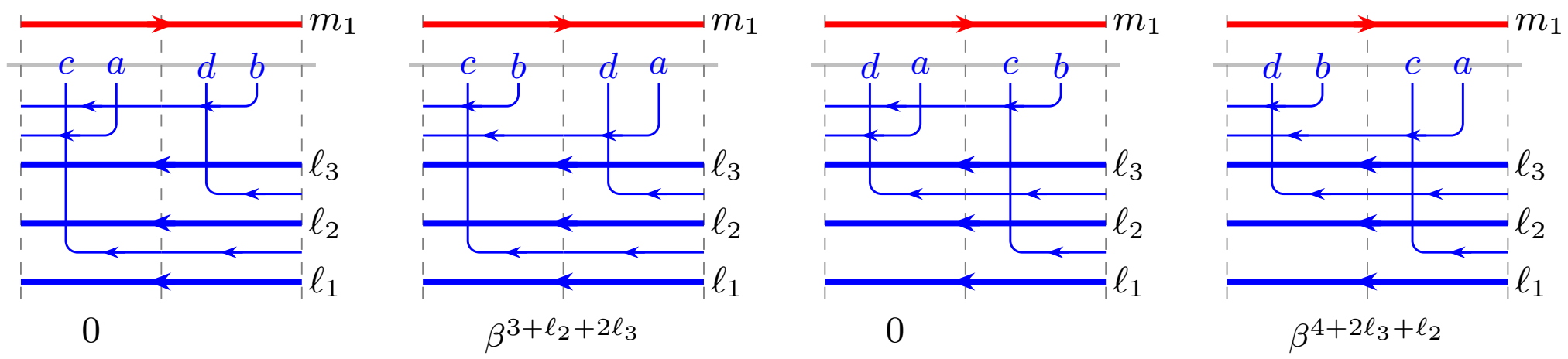
\[
\begin{matrix}
\begin{eye}{1}{8}{2}
  \wind{a}{}{7}{1}{0.5}
  \wind{}{c}{3}{1.5}{1}
  \wind{}{}{8}{1.5}{0.5}
  \wind{b}{}{8}{2}{1.5}
  \wind{}{d}{5}{2.5}{2}
  \wind{}{}{3}{2.5}{1.5}
  \wind[thick]{}{m_1}{1}{0.5}{2.5}
  \wind[thick]{\ell_1}{}{2}{2.5}{0.5}
  \wind[thick]{\ell_2}{}{4}{2.5}{0.5}
  \wind[thick]{\ell_3}{}{6}{2.5}{0.5}
  \label{1}{$0$}
\end{eye} &
\begin{eye}{1}{8}{2}
  \wind{b}{}{8}{1}{0.5}
  \wind{}{c}{3}{1.5}{1}
  \wind{}{}{7}{1.5}{0.5}
  \wind{a}{}{7}{2}{1.5}
  \wind{}{d}{5}{2.5}{2}
  \wind{}{}{3}{2.5}{1.5}
  \wind[thick]{}{m_1}{1}{0.5}{2.5}
  \wind[thick]{\ell_1}{}{2}{2.5}{0.5}
  \wind[thick]{\ell_2}{}{4}{2.5}{0.5}
  \wind[thick]{\ell_3}{}{6}{2.5}{0.5}
  \label{1.5}{\tiny $\beta^{3+\ell_2+2\ell_3}$}
\end{eye} &
\begin{eye}{1}{8}{2}
  \wind{a}{}{7}{1}{0.5}
  \wind{}{d}{5}{1.5}{1}
  \wind{}{}{8}{1.5}{0.5}
  \wind{b}{}{8}{2}{1.5}
  \wind{}{c}{3}{2.5}{2}
  \wind{}{}{5}{2.5}{1.5}
  \wind[thick]{}{m_1}{1}{0.5}{2.5}
  \wind[thick]{\ell_1}{}{2}{2.5}{0.5}
  \wind[thick]{\ell_2}{}{4}{2.5}{0.5}
  \wind[thick]{\ell_3}{}{6}{2.5}{0.5}
  \label{1}{$0$}
\end{eye} &
\begin{eye}{1}{8}{2}
  \wind{b}{}{8}{1}{0.5}
  \wind{}{d}{5}{1.5}{1}
  \wind{}{}{7}{1.5}{0.5}
  \wind{a}{}{7}{2}{1.5}
  \wind{}{c}{3}{2.5}{2}
  \wind{}{}{5}{2.5}{1.5}
  \wind[thick]{}{m_1}{1}{0.5}{2.5}
  \wind[thick]{\ell_1}{}{2}{2.5}{0.5}
  \wind[thick]{\ell_2}{}{4}{2.5}{0.5}
  \wind[thick]{\ell_3}{}{6}{2.5}{0.5}
  \label{1.5}{\tiny $\beta^{4+2\ell_3+\ell_2}$}
\end{eye}
\end{matrix}
\]

\appcaseanchor{9.7}{9-7}
\subsubsection*{Case 9.7} When $k<a=c<b<d$, all four weights are zero.

\appcaseanchor{9.8}{9-8}
\subsubsection*{Case 9.8} When $k<a=c<d<b$, all four weights are zero.

\appcaseanchor{9.9}{9-9}
\subsubsection*{Case 9.9} When $k<a<c<b=d$, all four weights are zero.

\appcaseanchor{9.10}{9-10}
\subsubsection*{Case 9.10} When $k<c<a<b=d$, only one nonzero case:
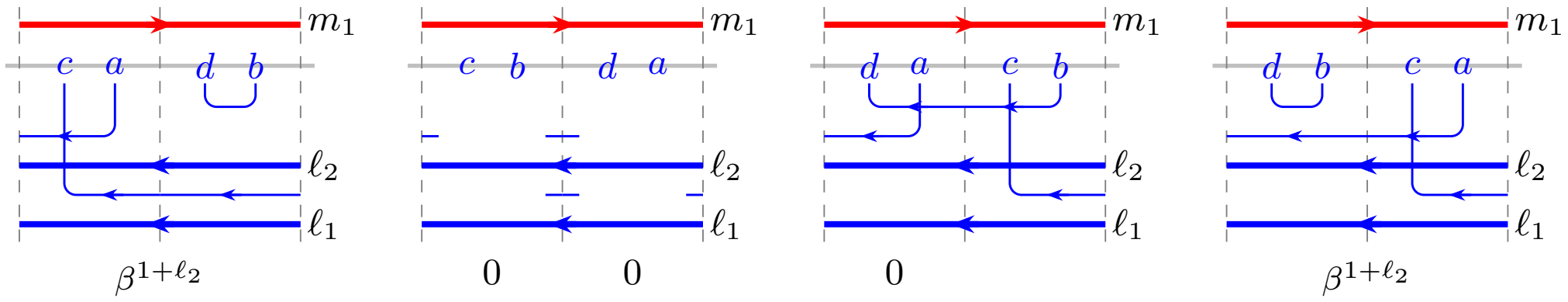
\[
\begin{matrix}
\begin{eye}{1}{6}{2}
  \wind{a}{}{5}{1}{0.5}
  \wind{}{c}{3}{1.5}{1}
  \wind{d}{b}{6}{2}{2}
  \wind{}{}{3}{2.5}{1.5}
  \wind[thick]{}{m_1}{1}{0.5}{2.5}
  \wind[thick]{\ell_1}{}{2}{2.5}{0.5}
  \wind[thick]{\ell_2}{}{4}{2.5}{0.5}
  \label{1.5}{\tiny $\beta^{1+\ell_2}$}
\end{eye} &
\begin{eye}{1}{6}{2}
  \wind[labels only]{c}{b}{2}{1}{1}
  \wind[labels only]{d}{a}{2}{2}{2}
  \brokenwind{}{}{5}{0.5}{1.5}
  \brokenwind[ends=right]{}{}{3}{0.5}{1.5}
  \brokenwind{}{}{3}{1.5}{2.5}
  \brokenwind[ends=left]{}{}{5}{1.5}{2.5}

  \wind[thick]{}{m_1}{1}{0.5}{2.5}
  \wind[thick]{\ell_1}{}{2}{2.5}{0.5}
  \wind[thick]{\ell_2}{}{4}{2.5}{0.5}
  \label{1}{$0$}
  \label{2}{$0$}
\end{eye} &
\begin{eye}{1}{6}{2}
  \wind{a}{}{5}{1}{0.5}
  \wind{}{d}{6}{1.5}{1}
  \wind{b}{}{6}{2}{1.5}
  \wind{}{c}{3}{2.5}{2}
  \wind[thick]{}{m_1}{1}{0.5}{2.5}
  \wind[thick]{\ell_1}{}{2}{2.5}{0.5}
  \wind[thick]{\ell_2}{}{4}{2.5}{0.5}
  \label{1}{$0$}
\end{eye} &
\begin{eye}{1}{6}{2}
  \wind{d}{b}{6}{1}{1}
  \wind{}{}{5}{1.5}{0.5}
  \wind{a}{}{5}{2}{1.5}
  \wind{}{c}{3}{2.5}{2}
  \wind[thick]{}{m_1}{1}{0.5}{2.5}
  \wind[thick]{\ell_1}{}{2}{2.5}{0.5}
  \wind[thick]{\ell_2}{}{4}{2.5}{0.5}
  \label{1.5}{\tiny $\beta^{1+\ell_2}$}
\end{eye}
\end{matrix}
\]

\appcaseanchor{9.11}{9-11}
\subsubsection*{Case 9.11} When $k<a<b=c<d$, all four weights are zero.

\appcaseanchor{9.12}{9-12}
\subsubsection*{Case 9.12} When $k<c<a=d<b$, only one nonzero case:
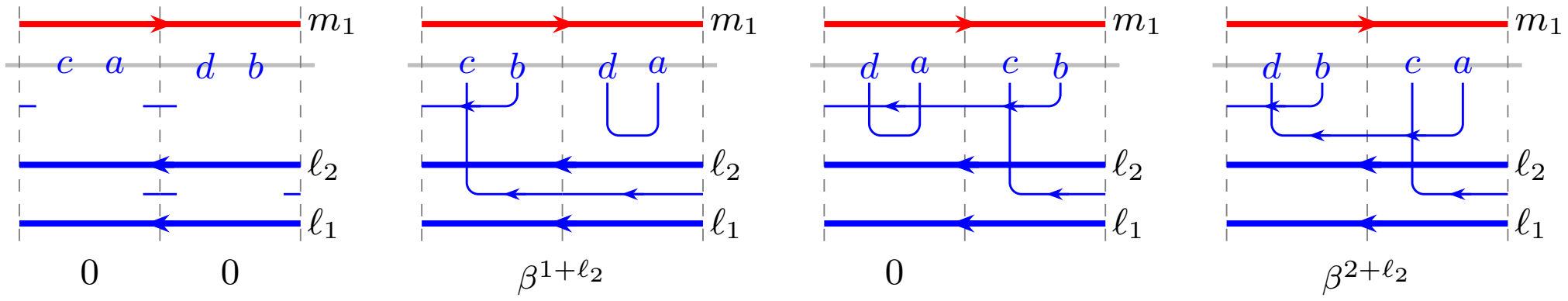
\[
\begin{matrix}
\begin{eye}{1}{6}{2}
  \wind[labels only]{c}{a}{2}{1}{1}
  \wind[labels only]{d}{b}{2}{2}{2}
  \brokenwind{}{}{6}{0.5}{1.5}
  \brokenwind[ends=right]{}{}{3}{0.5}{1.5}
  \brokenwind{}{}{3}{1.5}{2.5}
  \brokenwind[ends=left]{}{}{6}{1.5}{2.5}

  \wind[thick]{}{m_1}{1}{0.5}{2.5}
  \wind[thick]{\ell_1}{}{2}{2.5}{0.5}
  \wind[thick]{\ell_2}{}{4}{2.5}{0.5}
  \label{1}{$0$}
  \label{2}{$0$}
\end{eye} &
\begin{eye}{1}{6}{2}
  \wind{b}{}{6}{1}{0.5}
  \wind{}{c}{3}{1.5}{1}
  \wind{d}{a}{5}{2}{2}
  \wind{}{}{3}{2.5}{1.5}
  \wind[thick]{}{m_1}{1}{0.5}{2.5}
  \wind[thick]{\ell_1}{}{2}{2.5}{0.5}
  \wind[thick]{\ell_2}{}{4}{2.5}{0.5}
  \label{1.5}{\tiny $\beta^{1+\ell_2}$}
\end{eye} &
\begin{eye}{1}{6}{2}
  \wind{d}{a}{5}{1}{1}
  \wind{}{}{6}{1.5}{0.5}
  \wind{b}{}{6}{2}{1.5}
  \wind{}{c}{3}{2.5}{2}
  \wind[thick]{}{m_1}{1}{0.5}{2.5}
  \wind[thick]{\ell_1}{}{2}{2.5}{0.5}
  \wind[thick]{\ell_2}{}{4}{2.5}{0.5}
  \label{1}{$0$}
\end{eye} &
\begin{eye}{1}{6}{2}
  \wind{b}{}{6}{1}{0.5}
  \wind{}{d}{5}{1.5}{1}
  \wind{a}{}{5}{2}{1.5}
  \wind{}{c}{3}{2.5}{2}
  \wind[thick]{}{m_1}{1}{0.5}{2.5}
  \wind[thick]{\ell_1}{}{2}{2.5}{0.5}
  \wind[thick]{\ell_2}{}{4}{2.5}{0.5}
  \label{1.5}{\tiny $\beta^{2+\ell_2}$}
\end{eye}
\end{matrix}
\]

\appcaseanchor{9.13}{9-13}
\subsubsection*{Case 9.13} When $k<a=c<b=d$, only one nonzero case:
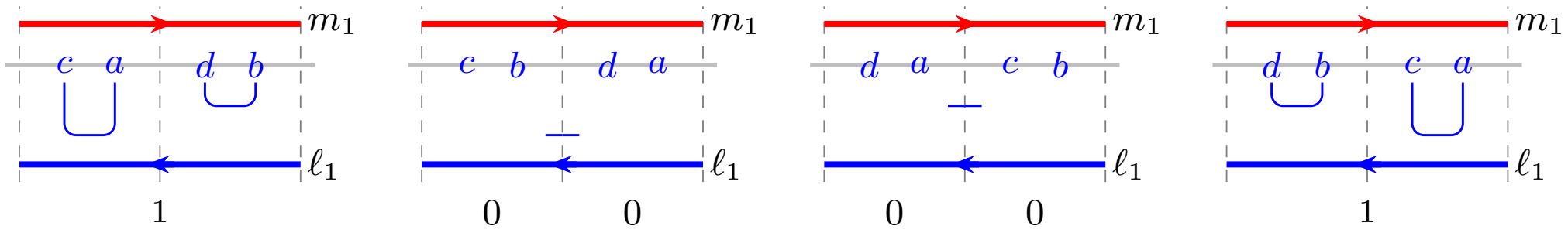
\[
\begin{matrix}
\begin{eye}{1}{4}{2}
  \wind{c}{a}{3}{1}{1}
  \wind{d}{b}{4}{2}{2}
  \wind[thick]{}{m_1}{1}{0.5}{2.5}
  \wind[thick]{\ell_1}{}{2}{2.5}{0.5}
  \label{1.5}{\tiny $1$}
\end{eye} &
\begin{eye}{1}{4}{2}
  \wind[labels only]{c}{b}{2}{1}{1}
  \wind[labels only]{d}{a}{2}{2}{2}
  \brokenwind[ends=right]{}{}{3}{0.5}{1.5}
  \brokenwind[ends=left]{}{}{3}{1.5}{2.5}

  \wind[thick]{}{m_1}{1}{0.5}{2.5}
  \wind[thick]{\ell_1}{}{2}{2.5}{0.5}
  \label{1}{$0$}
  \label{2}{$0$}
\end{eye} &
\begin{eye}{1}{4}{2}
  \wind[labels only]{d}{a}{2}{1}{1}
  \wind[labels only]{c}{b}{2}{2}{2}
  \brokenwind[ends=right]{}{}{4}{0.5}{1.5}
  \brokenwind[ends=left]{}{}{4}{1.5}{2.5}

  \wind[thick]{}{m_1}{1}{0.5}{2.5}
  \wind[thick]{\ell_1}{}{2}{2.5}{0.5}
  \label{1}{$0$}
  \label{2}{$0$}
\end{eye} &
\begin{eye}{1}{4}{2}
  \wind{d}{b}{4}{1}{1}
  \wind{c}{a}{3}{2}{2}
  \wind[thick]{}{m_1}{1}{0.5}{2.5}
  \wind[thick]{\ell_1}{}{2}{2.5}{0.5}
  \label{1.5}{\tiny $1$}
\end{eye}
\end{matrix}
\]


\begin{thebibliography}{99}


\bibitem{ACT}  D. Anderson, L. Chen and H.-H. Tseng, On the Finiteness of Quantum $K$-Theory of a Homogeneous Space, with Appendix B by H. Iritani. Int. Math. Res. Not. IMRN 2022, no. 2, 1313–1349.

\bibitem{AGM} D. Anderson, S. Griffeth and E. Miller,  
Positivity and Kleiman transversality in equivariant $K$-theory of homogeneous spaces.
J. Eur. Math. Soc. (JEMS) 13 (2011), no. 1, 57–84.





\bibitem{BFP99}
F. Brenti, S. Fomin and A. Postnikov,
Mixed Bruhat operators and Yang--Baxter equations for Weyl groups.
Internat. Math. Res. Notices, 1999(8)(1999), 419--441.


\bibitem{Buch} A.S. Buch, A Littlewood--Richardson rule for the  $K$-theory of Grassmannians.
Acta Math. 189 (2002), no. 1, 37–78.

\bibitem{BCMP} A.S. Buch, P.-E. Chaput, L.C. Mihalcea and N. Perrin, A Chevalley formula for the equivariant quantum $K$-theory of cominuscule varieties. Algebraic Geometry 5 (2018), no. 5, 568–595.

\bibitem{BCMP22} A.S. Buch, P.-E. Chaput, L. C. Mihalcea and N. Perrin, Positivity of minuscule quantum $K$ theory. arxiv: math.AG/2205.08630.

\bibitem{BCP0} A.S. Buch, P.-E. Chaput and N. Perrin, Equivariant rigidity of Richardson varieties. 
Pacific J. Math. 338 (2025), no. 2, 209–229.

\bibitem{BCP} A.S. Buch, P.-E. Chaput and N. Perrin, Seidel and Pieri products in cominuscule quantum  $K$-theory.
Math. Ann. 395 (2026), no. 4, Paper No. 81, 60 pp.

\bibitem{BCLM} A.S. Buch, S. Chung, C. Li and L.C. Mihalcea, Euler characteristics in the quantum $K$-theory of flag varieties.  Selecta Math. (N.S.) 26 (2020), no. 2, 29.

\bibitem{BM11} A.S. Buch and  L.C. Mihalcea,  Quantum  $K$-theory of Grassmannians. 
Duke Math. J. 156 (2011), no. 3, 501–538.



\bibitem{FGSX}
Neil J.Y. Fan, Peter L. Guo, C. Su and  R. Xiong,
 A Pieri type formula for motivic Chern classes of Schubert cells in Grassmannians,
arXiv: {math.CO/2402.04500}, 2024.

\bibitem{FGX}
N.J.Y. Fan, P.L. Guo, Rui Xiong, Bumpless pipe dreams meet puzzles, Advances in Mathematics 463 (2025), 110113.

\bibitem{FuLa} W. Fulton and  A. Lascoux, A Pieri formula in the Grothendieck ring of a flag bundle. Duke Math. J. 76 (1994), no. 3, 711–729.

\bibitem{FuWo} W. Fulton and C. Woodward,  
On the quantum product of Schubert classes. J. Algebraic Geom. 13 (2004), no. 4, 641–661.



\bibitem{Giv00} A. Givental, On the WDVV equation in quantum $K$-theory, dedicated to William Fulton on the occasion of his 60th birthday. Michigan Math. J. 48 (2000), 295–304.









\bibitem{GrKr} W. Graham and V. Kreiman, Excited Young diagrams, equivariant  $K$-theory, and Schubert varieties.
Trans. Amer. Math. Soc. 367 (2015), no. 9, 6597–6645.

\bibitem{GrRa} S. Griffeth and A. Ram, Affine Hecke algebras and the Schubert calculus. European J. Combin. 25 (2004), no. 8, 1263–1283.

\bibitem{GKM25} V. Gorbounov, C. Korff and L. C. Mihalcea, 
Quantum $K$-theory of Grassmannians from a Yang-Baxter algebra.
  arXiv: math.AG/2503.08602.


\bibitem{GKM} V. Gorbounov, C. Korff and L. C. Mihalcea,
Quantum $K$-structure constants and cylindric Grothendieck polynomials for Grassmannians. arXiv: math.CO/2610.31607.


\bibitem{HL15} Y. Huang and C. Li, On equivariant quantum Schubert calculus for $G/P$. J. Algebra 441 (2015), 21--56.

\bibitem{IIM} T. Ikeda, S. Iwao and T. Maeno,  
Peterson isomorphism in $K$-theory and relativistic Toda lattice. 
Int. Math. Res. Not. IMRN 2020, no. 19, 6421–6462.

\bibitem{IKN} T. Ikeda, T. Kouno and S. Naito,  
Seidel product formula in equivariant quantum $K$-theory of flag varieties, arXiv: math.AG/2602.24050.




\bibitem{INS16}
M. Ishii, S. Naito and D. Sagaki,
Semi-infinite Lakshmibai--Seshadri path model for level-zero
extremal weight modules over quantum affine algebras.
Adv. Math., 290(2016), 967--1009.


\bibitem{Kato3}
S. Kato. On quantum K-groups of partial flag manifolds. J. Math. Pures Appl. (9) 212 (2026), Paper No. 103916, 20 pp.




\bibitem{Kato1}
S. Kato. Loop structure on equivariant K-theory of semi-infinite flag manifolds. Ann.
of Math. (2) 202 (2025), no. 3, 1001–1075.


\bibitem{Kato8}
S. Kato, S. Naito, and D. Sagaki. Equivariant K-theory of semi-infinite flag manifolds and the Pieri-Chevalley formula. Duke Math. J., 169(13):2421–2500, 2020.


\bibitem{KZJ1} A. Knutson and P.  Zinn-Justin,  Schubert puzzles and integrability I: Invariant trilinear forms.
Commun. Am. Math. Soc. 6 (2026), 1–67.

 
\bibitem{KZJ3} A. Knutson and P.  Zinn-Justin,  Schubert puzzles and integrability III: separated descents.
arXiv: math.CO/2306.13855


\bibitem{KoKu} B. Kostant and S. Kumar, $T$-equivariant  $K$-theory of generalized flag varieties.  J. Differential Geom. 32 (1990), no. 2, 549–603. 

\bibitem{LLMS} T. Lam, C. Li, L.C. Mihalcea and M. Shimozono, A conjectural Peterson isomorphism in $K$-theory, J. Algebra 513 (2018), 326–343. 

\bibitem{LaSh}
T. Lam and M. Shimozono,  Quantum cohomology of $G/P$ and homology of affine Grassmannian,  Acta Math. 204 (2010), no. 1, 49–90.



\bibitem{Lee} Y.-P. Lee, Quantum $K$-theory. I. Foundations, Duke Math. J. 121 (2004), no. 3, 389-424.


\bibitem{LeMa}
C.~Lenart and T.~Maeno,
\newblock Quantum Grothendieck polynomials,
\newblock arXiv:math/0608232 (2006).



\bibitem{LNS}
C. Lenart, S. Naito and  D. Sagaki,
A general Chevalley formula for semi-infinite flag manifolds and quantum $K$-theory. 
Selecta Math. (N.S.) 30 (2024), no. 3, Paper No. 39, 44 pp.


\bibitem{LNSSS}
C. Lenart, S. Naito,   D. Sagaki, A. Schilling and  M. Shimozono, A uniform model for Kirillov--Reshetikhin crystals I: Lifting the parabolic quantum Bruhat graph, Int. Math. Res. Notices, 2015, no. 7, 1848–1901.

\bibitem{LeSo} C. Lenart and F. Sottile, A Pieri-type formula for the  $K$-theory of a flag manifold.
Trans. Amer. Math. Soc. 359 (2007), no. 5, 2317–2342. 

\bibitem{LLSY} C. Li, Z. Liu, J. Song and M. Yang, On  Seidel representation in quantum $K$-theory of   Grassmannians,  Sci. China Math. 68 (2025), no. 7, 1523–1548.

\bibitem{LiSo} C. Li and J. Song, Toward the quantum Pieri rule for $F\ell_n$ via Seidel representation (in Chinese), SCIENTIA SINICA Mathematica, 2024, 54(12): 2009–2022; an English translation is available at arXiv: math.AG/2507.12351.

\bibitem{Lus80}
G. Lusztig,
Hecke algebras and Jantzen's generic decomposition patterns.
Adv. Math., 37(2)(1980), 121--164.

\bibitem{MNSI}
T.~Maeno, S.~Naito, and D.~Sagaki,
\newblock A presentation of the torus-equivariant quantum $K$-theory ring of
flag manifolds of type $A$, Part I: the defining ideal,
\newblock \emph{J. Lond. Math. Soc.} \textbf{111} (2025), no.~3, e70095.

\bibitem{MNSII}
T.~Maeno, S.~Naito, and D.~Sagaki,
\newblock A presentation of the torus-equivariant quantum $K$-theory ring of
flag manifolds of type $A$, Part II: quantum double Grothendieck polynomials,
\newblock \emph{Forum Math. Sigma} \textbf{13} (2025), e19.

\bibitem{McN} P.J. McNamara,  Factorial Grothendieck polynomials. Electron. J. Combin. 13 (2006), no. 1, Research Paper 71, 40 pp.




\bibitem{MNS}
L.C.  Mihalcea, H. Naruse and C. Su, Left Demazure–Lusztig operators on equivariant (quantum)
cohomology and K-theory, Int. Math. Res. Not. IMRN, 16 (2022):12096–12147.


\bibitem{NNS}  S. Naito, D.-K. Nguyen and D. Sagaki, Pieri formula for the torus-equivariant quantum $K$-theory of flag manifolds of type $A$. In preparation. 


\bibitem{NaSa}
S.~Naito and D.~Sagaki, 
\newblock Pieri-type multiplication formula for quantum Grothendieck
polynomials,
\newblock \emph{Adv. Math.} \textbf{460} (2025), 110051. 



\bibitem{PeYo} O. Pechenik and A. Yong, 
Equivariant  $K$-theory of Grassmannians.
Forum Math. Pi 5 (2017), e3, 128 pp.

\bibitem{Pet97}
D. Peterson,
Quantum cohomology of $G/P$.
Lecture notes, Massachusetts Institute of Technology, Spring 1997.


\bibitem{PiRa} H. Pittie and A. Ram,  A Pieri-Chevalley formula in the  $K$-theory of a  $G/B$-bundle.
Electron. Res. Announc. Amer. Math. Soc. 5 (1999), 102–107.

\bibitem{Pos05}
A. Postnikov, Quantum Bruhat graph and Schubert polynomials.
Proc. Amer. Math. Soc., 133(3)(2005), 699--709.

\bibitem{Pos05-2}
A. Postnikov, Affine approach to quantum Schubert calculus, Duke Math. J., 128
(2005), no. 3, 473–509.


\bibitem{Sei} P. Seidel, $\pi_1$ of symplectic automorphism groups and invertibles in quantum homology rings, Geom. Funct. Anal. 7 (1997), 1046–1096.
 
\bibitem{WZJ} M. Wheeler, P. Zinn-Justin,  
Littlewood-Richardson coefficients for Grothendieck polynomials from integrability,
J. Reine Angew. Math. 757 (2019), 159–195.




\bibitem{Wood} C.T. Woodward, On D. Peterson’s comparison formula for Gromov-Witten invariants of $G/P$. Proc. Amer. Math. Soc. 133 (2005), no. 6, 1601–1609.

\bibitem{Xu} W. Xu, Quantum  $K$-theory of incidence varieties. 
Eur. J. Math. 10 (2024), no. 2, Paper No. 22, 47 pp.











\end{thebibliography}
\end{document}